\documentclass[a4paper,10pt]{amsart}
\usepackage[utf8]{inputenc}
\usepackage{amsmath}
\usepackage{amsfonts}
\usepackage{amssymb}
\usepackage{graphicx}
\usepackage[left=2.8cm,right=2.8cm,top=2cm,bottom=3cm]{geometry}
\numberwithin{equation}{section}
\graphicspath{ {images/} }
\usepackage{amsmath,amssymb, amsthm} 
\usepackage{enumerate}
\usepackage{lipsum}

\makeatletter
\def\@settitle{\begin{center}%
  \baselineskip14\p@\relax
  \bfseries
  \uppercasenonmath\@title
  \@title
  \ifx\@subtitle\@empty\else
     \\[1ex]\uppercasenonmath\@subtitle
     \footnotesize\mdseries\@subtitle
  \fi
  \end{center}%
}
\def\subtitle#1{\gdef\@subtitle{#1}}
\def\@subtitle{}
\makeatother

\newcommand\blfootnote[1]{%
  \begingroup
  \renewcommand\thefootnote{}\footnote{#1}%
  \addtocounter{footnote}{-1}%
  \endgroup
}

\usepackage{tikz-cd}
\usepackage{hyperref}
\hypersetup{colorlinks, linkcolor=black}

\usepackage[english]{babel}
\usepackage[colorinlistoftodos]{todonotes}

\theoremstyle{plain}
\newtheorem{thm}{Theorem}[subsection] % reset theorem numbering for each chapter
\newtheorem*{thmin}{Theorem}
\newtheorem*{propin}{Proposition}
\newtheorem*{corin}{Corollary}

\usepackage{mathrsfs}

\theoremstyle{definition}
\newtheorem{defi}[thm]{Definition}
\newtheorem{defiprop}[thm]{Definition/Proposition}
\newtheorem{rmk}[thm]{Remark}

\theoremstyle{definition}
\newtheorem{ex}[thm]{Example}
\theoremstyle{plain}
\newtheorem{prop}[thm]{Proposition}
\theoremstyle{plain}
\newtheorem{lemma}[thm]{Lemma}
\theoremstyle{plain}
\newtheorem{cor}[thm]{Corollary}


\newcounter{parentnumber}

\DeclareMathOperator{\ord}{ord}

\newcommand{\sym}{\operatorname{Sym}}

\newcommand{\Hom}{\operatorname{Hom}}

\newcommand{\res}{\operatorname{Res}}
\newcommand{\ind}{\operatorname{Ind}}

\newcommand{\dm}{\mathcal{DM}_{-}(k)}

\newcommand\SmallMatrix[1]{{%
  \tiny\arraycolsep=0.3\arraycolsep\ensuremath{\begin{pmatrix}#1\end{pmatrix}}}}

\newcommand{\N}{\mathfrak{N}}
\newcommand{\zita}{\mathfrak{z}}
\newcommand{\Zita}{\mathfrak{Z}}
\newcommand{\W}{\mathcal{W}}
\newcommand{\Z}{\mathbb{Z}}
\newcommand{\Q}{\mathbb{Q}}

\newcommand{\R}{\mathbb{R}}
\newcommand{\C}{\mathbb{C}}
\newcommand{\A}{\mathbb{A}}

\newcommand{\Oo}{\mathcal{O}}

\newcommand{\hh}{\mathcal{H}}

\newcommand{\Ss}{\mathcal{S}}
\newcommand{\F}{\mathcal{F}}
\newcommand{\T}{\mathcal{T}}
\newcommand{\p}{\mathfrak{p}}

\newcommand{\X}{\mathfrak{X}}

\newcommand{\mat}{\SmallMatrix{a & b \\ c & d}}

\DeclareMathOperator{\std}{Std}
\DeclareMathOperator{\Ext}{Ext}

\DeclareMathOperator{\ch}{ch}
\DeclareMathOperator{\Las}{L_{As}}

\DeclareMathOperator{\as}{As}
\DeclareMathOperator{\GL}{GL}
\DeclareMathOperator{\frob}{Frob}
\DeclareMathOperator{\gal}{Gal}

\DeclareMathOperator{\mot}{mot}
\DeclareMathOperator{\eis}{Eis}

\DeclareMathOperator{\et}{\text{\'{e}}t}

\DeclareMathOperator{\tr}{Tr}

\DeclareMathOperator{\reg}{reg}

\DeclareMathOperator{\tsym}{TSym}

\DeclareMathOperator{\rhoas}{\rho^{As}}

\DeclareMathOperator{\re}{Re}
\DeclareMathOperator{\id}{Id}
\DeclareMathOperator{\nm}{Nm}
\DeclareMathOperator{\cores}{cores}
\DeclareMathOperator{\stab}{Stab}
\DeclareMathOperator{\vol}{Vol}
\DeclareMathOperator{\AF}{\mathcal{AF}^{k,k',j}_{mot}}
\DeclareMathOperator{\BF}{\mathcal{BF}^{k,k',j}_{mot}}
  \newcommand{\Addresses}{{% additional braces for segregating \footnotesize
  \bigskip
  \footnotesize

\textsc{London School of Geometry and Number Theory, UCL, Department of Mathematics, Gower street, WC1E 6BT, London, UK}\par\nopagebreak
  \textit{E-mail address}, G.~Grossi: \texttt{giada.grossi.16@ucl.ac.uk}

}}

\begin{document}
\author{Giada Grossi}
\title{On norm relations for Asai--Flach classes}\blfootnote{\emph{Date}. 4th October 2018}

\begin{abstract}
We give a new proof of the norm relations for the Asai--Flach Euler system built by Lei--Loeffler--Zerbes. More precisely, we redefine Asai--Flach classes in the language used by Loeffler--Skinner--Zerbes for Lemma--Eisenstein classes and prove both the vertical and the tame norm relations using local zeta integrals. These Euler system norm relations for the Asai representation attached to a Hilbert modular form over a quadratic real field $F$ have been already proved by Lei--Loeffler--Zerbes for primes which are inert in $F$ and for split primes satisfying some assumption; with this technique we are able to remove it and prove tame norm relations for all inert and split primes. 
\end{abstract}
\maketitle

\tableofcontents

\section*{Introduction}\label{intro}
The theory of Euler systems is a powerful tool used in modern number theory for studying Galois representations: they play an important role in bounding Selmer groups of such Galois representations, conditional on the non-vanishing of a value of the associated $L$-function.

An Euler system for a $p$-adic representation $V$ of $G_{\Q}$, the absolute Galois group of $\Q$, (unramified outside a finite set $\Sigma$ of places including $p$) is a collection of Galois cohomology-classes $(z_m)_{m\geq 1}$ with $z_m\in H^1(\Q(\mu_m),V^*(1))$ satisfying the following norm relations:
\begin{equation}\label{eqNR}
\text{cores}^{\Q(\mu_{m\ell})}_{\Q(\mu_{m})}z_{m\ell}=\begin{cases}
z_m &\ell\mid m\quad\text{or}\quad \ell\in \Sigma\\
P_{\ell}(\frob_{\ell}^{-1})z_m &\text{otherwise},
\end{cases}
\end{equation}
where $\frob_{\ell}^{-1}$ is the arithmetic Frobenius and $P_{\ell}(x)=\det(1-\frob_{\ell}^{-1}x|V)$ is its characteristic polynomial. 

Even though such Euler systems are expected to exist for ``representations coming from geometry'', it is very difficult to construct non-trivial examples. In the last few years some new Euler systems have been constructed, e.g. an Euler system for the $p$-adic representation attached to the Rankin-Selberg conveolution of two modular forms  \cite{RSCMF}, for the Asai representation of a quadratic Hilbert modular form \cite{HMS} and more recently for the spin representation of a genus $2$ Siegel modular form  \cite{GSP4}. The common input of these works, following the ideas of \cite{Kato}, are Siegel units, which are invertible elements in $\Oo(Y_H)$, where $H=\GL_2$ and $Y_H$ is the Shimura variety corresponding to the algebraic group $H$. More generally one considers Eisenstein classes, which are elements in the first motivic cohomology group of $Y_H$ with coefficients in some specific motivic sheaves. 

The idea of the aforementioned papers is then to consider embeddings $\GL_2\hookrightarrow G$ (or $\GL_2\times_{\GL_1} \GL_2\hookrightarrow G$ in \cite{GSP4}),  where $G$ is a suitable algebraic group. Pushing forward the Siegel units via the corresponding  embedding at the level of Shimura varieties, one gets classes in a motivic cohomology group of the desired variety $Y_G$. Such embeddings are then suitably ``perturbed'' in order to define classes in the motivic cohomology of the base change over cyclotomic extensions $Y_G \times \mu_m$. Via the \'etale regulator one obtains  classes in the \'{e}tale cohomology of $Y_G\times \mu_m$, which give rise to classes in Galois cohomology via the Hochschild-Serre spectral sequence.  The Galois  representations are the ones attached to (twists) of automorphic representations for the group $G$; they appear in the middle degree \'{e}tale cohomology of $Y_{G,\overline{\Q}}$.

The main difficulty in proving that the classes form an Euler system is the proof of the tame Euler system norm relations, i.e. comparing classes $z_{m\ell}$ and $z_m$ when $\ell\nmid m$. In the Rankin--Selberg (\cite{RSCMF}) and in the Asai case (\cite{HMS}), these relations are proved via some explicit computations in the Hecke algebra. This approach would have been much more difficult  (or even impossible) for the Euler system attached to a genus $2$ Siegel modular form, as the structure of the group $\mathrm{GSp}_4$ is too complicated. In \cite{GSP4}, indeed, the technique used was different: the norm relations were obtained using results from smooth representation theory. In this note we show that it is possible to adapt this strategy to give new proofs of the norm relations in the Rankin--Selberg and Asai case (Theorem \ref{tame}). We will give the details for the Asai case, which requires separate computations for the cases when $\ell$ is split and when $\ell$ is inert. The split prime computations also give a new proof of the tame norm relations in the  Rankin--Selberg case.

In \cite{HMS}, the authors found classes for (twists of) the Asai representation of a Hilbert cuspidal eigenform over a real quadratic field $F$. These Euler system classes are shown to satisfy the norm relations (\ref{eqNR}) for primes $\ell$ which are either inert in $F$ or which are split and satisfy the additional assumption that the two prime ideals of $F$ above it are narrowly principal (i.e. principal ideal whose generator is a totally real element of $F$). Using smooth representation theory, we are able to fully prove norm relations, with no such assumption on $\ell$.

\subsection*{Strategy and organisation of the paper} We consider the embedding of algebraic groups over $\Q$
\begin{equation}\label{eq:embgrps}
H:=\GL_2 \ \hookrightarrow \ G:= \res_{F/\Q}\GL_2.
\end{equation}
We will be working with automorphic representations $\Pi$ of $\GL_2(\A_{F,f})$, where $\A_{F,f}$ denotes the finite ad\`{e}les over $F$. Equivalently we can view $\Pi$ as a representation of $G(\A_f)$, where $\A_f$ are the finite ad\`{e}les over $\Q$. At every place $\ell$ we have a representation of $G(\Q_{\ell})$, which, depending on $\ell$ has the following shape:
\begin{itemize}
\item if $\ell$ is inert and hence $G(\Q_{\ell})=\GL_2(F_{\ell})$ where $F_{\ell}$ is an unramified quadratic extension of $\Q_{\ell}$, the representation is $\Pi_{\ell}$;
\item if $\ell=v_1\cdot v_2$ is split and hence $G(\Q_{\ell})\simeq\GL_2(\Q_{\ell})\times \GL_2(\Q_{\ell})$, the representation is $\Pi_{v_1}\otimes \Pi_{v_2}$.
\end{itemize} 
In order to (re)define the Euler system constructed in \cite{HMS}, we define a special map $\AF$ for $k,k'\geq 0$ integers and $0\leq j \leq \min(k,k')$ with values in $W=H^3_{\mot}(Y_G,\mathcal{D}(2))$, where $Y_G$ is the Shimura variety associated to $G$ and $\mathcal{D}$ is a motivic sheaf depending on $k,k',j$. Such map will be of ``global nature'', more precisely it is a map 
\begin{displaymath}
\AF: \Ss(\mathbb{A}_f^2, \Q) \otimes \hh(G(\A_f),\Q) \longrightarrow H^3_{\mot}(Y_G,\mathcal{D}(2))
\end{displaymath}
satisfying some conditions of $H(\A_f)\times G(\A_f)$-equivariance. Here $\Ss(\mathbb{A}_f^2, \Q)$ denotes the space of Schwartz functions on $\A_f$ and $\hh(G(\A_f),\Q)$ the Hecke algebra over $G(\A_f)$. The Asai--Flach classes will be defined as images via $\AF$ of precise elements in $\Ss(\mathbb{A}_f^2, \Q) \otimes \hh(G(\A_f),\Q)$. Proving norm relations (in motivic cohomology) will turn out to be equivalent to proving relations of such classes locally at a certain prime $\ell$, i.e. we will be looking at a map 
\begin{displaymath}
(\AF)_{\ell}: \Ss(\Q_{\ell}^2, \Q) \otimes \hh(G(\Q_{\ell}),\Q) \longrightarrow H^3_{\mot}(Y_G,\mathcal{D}(2)).
\end{displaymath}
While we will be able to prove $p$-direction norm relations already in motivic cohomology, in order to prove ``tame norm relations'' we will have strong assumptions on the target of such map. We will have to apply the \'{e}tale regulator and Hochschild–-Serre spectral sequence to pass to Galois cohomology and finally take the projection to an automorphic representation of $G(\A_f)$ as above. The local components $\Pi_v$ at a ``good prime'' $v$ of $F$ will be an irreducible spherical principal series representation $I_{\GL_2(F_v)}(\chi,\psi)$, for $\chi,\psi$ unramified characters of $F_v^{\times}$. Hence we will need to work and prove results for maps
\begin{displaymath}
\Zita: \Ss(\Q_{\ell}^2,\C)\otimes \hh(G) \longrightarrow  \sigma=\begin{cases}
I_{\GL_2(F_{\ell})}(\chi,\psi) &\text{    if $\ell$ inert, $[F_{\ell}:\Q_{\ell}]=2$} \\
I_{\GL_2(\Q_{\ell})}(\chi_1,\psi_1)\otimes  I_{\GL_2(\Q_{\ell})}(\chi_2,\psi_2) &\text{    if $\ell$ splits}.
\end{cases}
\end{displaymath} 

%We will also need consider the ``intermediate'' algebraic group $G^*:= G\times_{\res_{F/\Q}\mathbb{G}_m}\mathbb{G}_m$. Considering cohomology of $Y_{G^*}$ instead of the one of $Y_G$ will allow to find classes in Galois cohomology over cyclotomic extensions.

We proceed first by recalling in $\S \ref{localrep}$ some useful tools to work with representations of $\GL_2$ over a local field. The following sections, $\S$\ref{inertzeta} and $\S$\ref{splitzeta}, should be thought in parallel: we move to local representations of $G$ over $\Q_{\ell}$ proving the same results for both the inert and split case, giving explicit descriptions of local $L$-factors $L(\as(\sigma),s)$ of principal series representations $\sigma$ as above as local zeta integrals. 

In $\S$\ref{towards}, we will relate images through maps as above of elements in $\Ss(\Q_{\ell}^2,\C)\otimes \hh(G)$ given by Definitions \ref{defischwt}, \ref{defischinfinity} and \ref{etadef}. The main results of this section are 

\begin{propin}[Proposition \ref{THEprop}]\label{propintro}
For any $\Zita: \Ss(\Q_{\ell}^2,\C)\otimes \hh(G)\to W$, where $W$ is a smooth complex representation of $G(\Q_{\ell})$ we have
\begin{displaymath}
\Zita(\phi_{1,\infty}\otimes \ch(\eta_{m+1}K_{m,n}))=\begin{cases}
\tfrac{1}{\ell}U'(\ell)  \\
\tfrac{1}{\ell -1}(U'(\ell)-1) 
\end{cases} \cdot \Zita(\phi_{1,\infty}\otimes \ch(\eta_{m}K_{m,n})) \ \  \begin{matrix}
\text{if } m\geq 1\\
\text{if } m=0.
\end{matrix}
\end{displaymath}
\end{propin}

\begin{corin}[Corollary \ref{THEcor}]\label{corintro}
Let $W=\sigma^{\vee}$ for $\sigma$ a principal series representation with central character $\chi_{\sigma}$. Let $\chi=|\cdot |^{1/2+k+h} \tau$, for $\tau$ a finite order character and $k,h\geq 0$, and $\psi=|\cdot |^{-1/2+h}$. Assume 
\begin{equation}\label{eqcentral}
\chi\psi\cdot\chi_{\sigma}=1
\end{equation}
and some additional condition on $\sigma$. Let $\Zita:\Ss(\Q_{\ell}^2,\C)\otimes \hh(G)\to \sigma^{\vee}$ such that it factors as $\Zita = \Zita'\circ f$, where $f$ is the Siegel section map defined in $\S$\ref{localrep} and 
\begin{displaymath}
\Zita': I_{\GL_2(\Q_{\ell})}(\chi,\psi)\otimes \hh(G) \to  \sigma^{\vee}.
\end{displaymath}
Then we have
\begin{displaymath}
\Zita(\phi_{1,\infty} \otimes (\ch(K)-\ch(\eta_1K)))=\tfrac{\ell}{\ell-1}L(\as(\sigma),h)^{-1}\cdot \Zita(\phi_0,\ch(K)).
\end{displaymath}
\end{corin}

While the proposition can be proved directly for any such $\Zita$, the corollary follows from Theorem \ref{thmzita}. It states the analogous equality for any function in $\Hom_{\GL_2(\Q_{\ell})}(I_{\GL_2(\Q_{\ell})}(\chi,\psi)\otimes\sigma, \C)$ (which is canonical bijection with the space of functions $\Zita'$ as above). The proof of the theorem follows from an explicit proof of the claimed equality for a specific choice of a nonzero element $\zita_{\chi,\psi}\in \Hom_{\GL_2(\Q_{\ell})}(I_{\GL_2(\Q_{\ell})}(\chi,\psi)\otimes\sigma, \C)$ built using the local zeta integrals of $\S\S$2,3 (see Definition \ref{defibasis}). One then crucially needs the following multiplicity one result in order to prove it for any  $\zita \in \Hom_{\GL_2(\Q_{\ell})}(I_{\GL_2(\Q_{\ell})}(\chi,\psi)\otimes\sigma, \C)$.
\begin{thmin}[Theorem \ref{thmmult1}, Multiplicity one]
Let $\sigma,\chi,\psi$ satisfying $($\ref{eqcentral}$)$ and the same additional condition of the Corollary. Assume that $\chi\psi^{-1}\neq |\cdot|^{-1}$. We have 
\begin{displaymath}
\dim\left( \Hom_{\GL_2(\Q_{\ell})}(I_{\GL_2(\Q_{\ell})}(\chi,\psi)\otimes\sigma, \C) \right)\leq 1.
\end{displaymath}
\end{thmin}
This theorem follows from \cite[Theorem 1.1]{prasad} in the case where $\ell$ splits and $I_{\GL_2(\Q_{\ell})}(\chi,\psi)$ is irreducible and it is proved in Theorem \ref{thmmult1} for the remaining cases. We use tools of Mackey theory following the strategy used by Prasad in \textit{op. cit.} in the proof of the principal series representations case and a result of \cite{shintani} in some degenerate case.

In $\S$\ref{eisenstein} we recall the definition of Eisenstein classes as $H(\A_f)$-equivariant maps
\begin{displaymath}
\Ss(\A_f^2,\Q) \longrightarrow H^1_{\mot}(Y_H, \tsym^k\hh_{\mathbf{Q}}(\mathcal{E})(1)).
\end{displaymath}

We proceed in $\S$\ref{defi} with the definition of the Asai--Flach map $\AF$, constructed using the Eisenstein classes map and the pushforward in motivic cohomology induced by (\ref{eq:embgrps}). We then define the motivic Asai--Flach elements as image via this map of elements in $\Ss(\A_f^2,\Q)\otimes \hh(G(\A_f))$ that are described at every place in terms of the elements considered in the above Proposition and Corollary.
%\begin{displaymath}
%z_{M,m,n}^{[k,k',j]}:= \tfrac{1}{\vol(W)}\AF (\phi_{M,m,n}\otimes\xi_{M,m,n} ) \in H^{3}_{\mot}(Y_G(K_{M,m,n}), \mathcal{D}(2)).
%\end{displaymath}
%In particular the level subgroup $K_{M,m,n}$ will be defined by matrices congruent to $\SmallMatrix{*&*\\0&1}$ modulo $p^n$ and having determinant congruent to 1 modulo $Mp^m$. Writing $K^*_{M,m,n}=K_{M,m,n}\cap G^*(\A_f)$, this will give isomorphisms 
%\begin{displaymath}
%Y_{G^*}(K^*_{M,m,n})\simeq Y_{G^*}(K^*_n)\times_{\Q} \Q(\mu_{Mp^m}),
%\end{displaymath}
%for $K^*_n$ the subgroup of $G^*(\A_f)$ given by matrices congruent to $\SmallMatrix{*&*\\0&1}$ modulo $p^n$.

Finally in $\S$\ref{ES} we prove some pushforward compatibilities of such motivic classes (Theorem \ref{changelevel} and \ref{quasivert}) using the local result given by the above Proposition. We then use these classes to find elements in Galois cohomology of the Asai representation of a Hilbert cuspidal eigenform and prove Euler system norm relations (Theorem \ref{verticalthm} and \ref{tame}); vertical norm relations follow from the $p$-direction compatibility of motivic classes, while tame norm relations rely on the local result of the above Corollary.

\subsubsection*{Acknowledgements.} This work is part of the author's PhD thesis under the supervision of Sarah Zerbes. I would like to thank her for suggesting the problem and for her continuous guidance and careful proofreading. I would also like to thank David Loeffler for many insightful comments. Special thanks to Antonio Cauchi and Joaqu\`{i}n Rodrigues for helpful comments and discussions. Part of this work was conducted while the author was visiting the Bernoulli Center at EPFL, during the semester ``Euler Systems and Special Values of L-functions''. I am grateful to the Bernoulli Center and to Dimitar Jetchev for their hospitality and invitation. 

\section{Local representation theory for \texorpdfstring{$\GL_2$}{GL2}}\label{localrep}
In this section we recall the standard tools of local representation theory that will be useful later in the proof of norm relations. We follow \cite[Chapter 4]{bump}.

We let $E$ be a non-Archimedean local field and denote with $\Oo, \p, \varpi$ respectively the ring of integers in $E$, the maximal ideal and a fixed uniformiser of $\p$. Let $|\cdot|$ be the norm and $q$ such that $|\varpi|=q^{-1}$. We also fix an Haar measure $dx$ on $E$ and $d^{\times}x$ on $E^{\times}$. For a smooth character $\chi$ of $E^{\times}$ we define its local $L-$factor 
\begin{displaymath}
L(\chi,s)=L(\chi|\cdot|^s,0)= \begin{cases}
(1-\chi(\varpi)q^{-s})^{-1} &\text{if }\chi_{|_{\Oo^{\times}}}=1 \\
1 &\text{otherwise}.
\end{cases}
\end{displaymath}

\subsection{Induced representations} We recall here some basics about induced representations of totally disconnected topological groups. See for example \cite[$\S\S$2.21-2.29]{bern}. Let $X$ be a group as above with a right Haar measure $d_R$ on $X$ and a left Haar measure $d_L$. The modular quasicharacter $\delta_X$ of $X$ is defined by $d_R(x)=\delta_X(x)d_L(x)$. 

Recall that $X$ is said to be unimodular if $\delta_X=1$. A trivial example of unimodular group is any abelian group. A less trivial example is $X=\GL_n(E)$. A group which is not unimodular is the Borel subgroup $B$ of $\GL_n(E)$. For $n=2$, its modular quasicharacter is given by
\begin{displaymath}
\delta_B(\SmallMatrix{a&b\\0&d})=|\tfrac{a}{d}|.
\end{displaymath}

Consider now $Y$ a closed subgroup of $X$. We have a restriction functor from the category of smooth representations of $X$ to smooth representations of $Y$. This functor has a left and a right adjoint, given by induction and compact induction. 

\begin{defi}
Let $(V, \tau)$ be a smooth representation of $Y$. We denote with $\ind_Y^X\tau$ the space of smooth functions $f:X\to V$ satisfying the following condition
\begin{center}
 $f(yx)=\delta_X^{-1/2}(y)\delta_Y^{1/2}(y)f(x)$ for every $x\in X,y\in Y$.
\end{center}
We denote with $\text{c-}\ind_Y^X(\tau)$ the subspace of $\ind_Y^X\tau$ consisting of functions which in addition are compactly supported modulo $Y$. This coincides with $\ind_Y^X(\tau)$ when $X/Y$ is compact. These are $X$-representation with action of $X$ given by right multiplication.
\end{defi}

\begin{thm}[Frobenius reciprocity]
Let $(W,\sigma)$ be a smooth representation of $X$ and $(V,\tau)$ a smooth representation of $Y$. Denote with $( \ )^{\vee}$ the smooth dual of a representation. We then have the following isomorphisms:
\begin{displaymath}
\Hom_X(\sigma,\ind_Y^X\tau)\simeq \Hom_Y(\sigma_{|Y},\delta_Y^{1/2}\delta_X^{-1/2}\tau);
\end{displaymath}
\begin{displaymath}
\Hom_X(\text{c-}\ind_Y^X\tau,\sigma^{\vee})\simeq \Hom_Y(\delta_Y^{-1/2}\delta_X^{1/2}\tau,(\sigma_{|Y})^{\vee}).
\end{displaymath}
\end{thm}

%If $X$ is an algebraic group, we can realise induced representations as sections of a vector bundle $\mathcal{B}$ over the algebraic variety $X/Y$ with fibres $V$. Let $\mathcal{B}$ be the space $X\times V$ with an equivalence relation given by
%\begin{displaymath}
%(x,v)\sim (xy, \delta_X^{1/2}(y)\delta_Y^{-1/2}(y)\tau(y)v) \ \ \ \text{for }x\in X,y\in Y,v\in V.
%\end{displaymath}
%It is straightforward to see that, writing $\Gamma(X/Y,\mathcal{B})$ for the smooth sections of $\mathcal{B}$ and $\Gamma_c(X/Y,\mathcal{B})$ for the compactly supported smooth sections, we have
%\begin{displaymath}
%\Gamma(X/Y,\mathcal{B})=\ind_Y^X(\tau) \ \ \ \text{and} \ \ \ \Gamma_c(X/Y,\mathcal{B})= \text{c-}\ind_Y^X(\tau).
%\end{displaymath}
%We now state a general lemma about compactly supported smooth sections of vector bundles on totally disconnected locally compact algebraic groups (see for example \cite[Lemma 5.1]{prasad}).
%\begin{lemma}
%Let $X$ be a totally disconnected locally compact algebraic group, $Z$ a closed subgroup and $B$ a vector bundle over $X$. Then we have an exact sequence 
%\begin{displaymath}
%0\to\Gamma_c(X-Z,{B}_{|X-Z})\to \Gamma_c(X,{B})\to\Gamma_c(Z,{B}_{|Z})\to 0.
%\end{displaymath}
%\end{lemma}
One can actually realise induced representations as sections of a vector bundle. Using such description, one finds the following exact sequence of induced representations that will be useful later. Let $H,J$ be closed subgroups of a totally disconnected locally compact algebraic group $G$ and $\tau$ a smooth representation of $J$. Assume that the quotient $H\backslash G /J$ has two elements. This means that the action of $H$ on the space $G/J$ has two orbits, one open and one closed. We can write these two orbits as $H/H_1,H/H_2$, where $H_1=\stab_H(1\cdot J)$ and $H_2=\stab_H(\epsilon\cdot J)$, where $\epsilon\in G$ such that $\epsilon\cdot J$ is in the open orbit. We can compute these subgroups as follows
\begin{displaymath}
H_1=H\cap J, \ \ \ H_2=H\cap \epsilon^{-1}J\epsilon.
\end{displaymath} 
Applying \cite[Lemma 5.1]{prasad} to the closed subgroup $H/H_1$ and the complement$H/H_2$ and normalising appropriately one finds an exact sequence of $H$-modules
\begin{equation}\label{eq:exseq}
0\to \text{c-}\ind_{H_2}^H\tau_2 \to (\text{c-}\ind_J^G\tau)_{|H} \to \text{c-}\ind_{H_1}^H\tau_1 \to 0,
\end{equation}
where $\tau_1=\delta_J^{1/2}\cdot\delta_{H_1}^{-1/2}\cdot \tau_{|H_1}$ and $\tau_2$ is a representation of $H_2$ given by
\begin{displaymath}
\tau_2(h)=\delta_J^{1/2}(\epsilon h\epsilon^{-1})\delta_{H_2}^{-1/2}(h)\tau(\epsilon h\epsilon^{-1}).
\end{displaymath}

\subsection{Principal series representations}
\begin{defi} Let $H=\GL_2(E)$ and $\xi,\psi$ two quasicharacters of $E^{\times}$. We define a space of functions on $H$ as follows 
\begin{displaymath}
I_H(\chi,\psi):= \lbrace f: H \rightarrow \C \text{ smooth }: f(\SmallMatrix{a & * \\ 0 & d}\cdot h)=  \left|\tfrac{a}{d}\right|^{1/2}\chi(a)\psi(d)f(h) \rbrace
\end{displaymath}
\end{defi}
We will denote $I_H(\chi,\psi)$ simply as $I(\chi,\psi)$. We see $I(\chi,\psi)$ as a $\GL_2(E)$-representation letting $\GL_2(E)$ act by right translation. 
%i.e. for $g\in \GL_2(E)$
%\begin{displaymath}
%g\cdot f(h)=f(hg) \text{   for every $f\in I(\chi,\psi),h\in\GL_2(E)$}.
%\end{displaymath}
Letting $\tau$ be the one dimensional representation of the Borel subgroup (consisting of upper triangular matrices) $B(E)$ given by $\tau(\SmallMatrix{a & * \\ 0 & d})=  \chi(a)\psi(d)$, we have $I(\chi,\psi)=\ind_{B(E)}^{\GL_2(E)}\tau$.
\begin{defi}
The $\GL_2(E)$-representations $I(\chi,\psi)$, for $\chi,\psi$ quasicharacters of $E^{\times}$, are called \textit{principal series representations}.
\end{defi}

To characterise such representations, we recall the definition of the intertwining operator. Fix $\chi,\psi$ and write them as
\begin{displaymath}
\chi=|\cdot|^{s_1}\xi_1, \ \ \psi=|\cdot|^{s_2}\xi_2,
\end{displaymath}
where $s_i\in\C$ and $\xi_i$ are unitary characters. Let $f\in I(\chi,\psi)$. We write, for $h\in\GL_2(E)$,
\begin{displaymath}
Mf(h)=\int_F f(w\cdot m_x \cdot h)dx, 
\end{displaymath}
where $w=\SmallMatrix{0&-1\\1&0}$ and $m_x=\SmallMatrix{1 & x \\ 0 & 1}$.

By \cite[Proposition 4.5.6]{bump}, if $\re(s_1-s_2)>0$ then the above integral is absolute convergent and it defines a nonzero intertwining map
\begin{align*}
M: I(\chi &,\psi) \rightarrow I(\psi,\chi)\\ &f\mapsto Mf.
\end{align*}

In the case where $\re(s_1-s_2)<0$ we clearly have an analogous operator $M$ obtained by switching $\chi$ and $\psi$. The procedure for defining such an operator in the case where $\re(s_1-s_2)=0$ uses \textit{flat sections} and the fact that we can write $\GL_2(E)=B(E)\cdot K$, where $K=\GL_2(\Oo)$ (known as \textit{Iwasawa decomposition}, see \cite[Proposition 4.5.2]{bump}). Indeed one starts with noticing that $f\in I(\chi,\psi)$ is uniquely determined by its restriction $f_0$ to $K$.
%, which satisfies 
%\begin{displaymath}
%f(\SmallMatrix{a & b \\ 0 & d}\cdot k)=  \xi_1(a)\xi_2(d)f(k),
%\end{displaymath}
%for $a,d\in\Oo^{\times}$,$b\in\Oo_F$ and $k\in K$. We denote with $V_0$ the space of smooth functions on $K$ satisfying this condition, having fixed $\xi_1,\xi_2$. Then for any $s_1,s_2\in \C$ and $f_0\in V_0$ t
One shows that there exists a unique extension of $f_0$ to an element $f_{s_1,s_2}$ in $V_{s_1,s_2}:=I(|\cdot|^{s_1}\xi_1,|\cdot|^{s_2}\xi_2)$. The function
\begin{displaymath}
(s_1,s_2)\mapsto f_{s_1,s_2},
\end{displaymath} 
is called a flat section. We then have
\begin{prop}\label{interwop} \cite[Proposition 4.5.7]{bump}.
For a fixed $h\in \GL_2(E)$ the integral $Mf_{s_1,s_2}(h)$ defined as above for $\re(s_1-s_2)>0$ has analytic continuation to all $s_1,s_2$ where $\chi\neq \psi$. We hence have defined an intertwining operator 
\begin{displaymath}
M: I(|\cdot|^{s_1}\xi_1,|\cdot|^{s_2}\xi_2) \rightarrow I(|\cdot|^{s_2}\xi_2, |\cdot|^{s_1}\xi_1).
\end{displaymath}
\end{prop}

We have the following theorem characterising principal series representations.
\begin{thm}\label{thmmult}\cite[Theorem 4.5.1 and 4.5.2]{bump}.
Let $\chi,\psi$ be quasicharacter of $E^{\times}$. Then $I(\chi,\psi)$ is an irreducible $\GL_2(E)$-representation except in the following two cases
\begin{itemize}
\item[(i)] if $\chi\psi^{-1}=|\cdot|^{-1}$, then $I(\chi,\psi)$ has a one-dimensional invariant subspace and the quotient representation is irreducible; 
\item[(ii)] if $\chi\psi^{-1}=|\cdot|$, then $I(\chi,\psi)$ has an irreducible codimension one invariant subspace.
\end{itemize}
If $I(\chi,\psi)$ is irreducible, then it is isomorphic to $I(\psi,\chi)$ via the intertwining operator $M$. Moreover if we have two such representations, for quasicharacters $\chi_1,\psi_1,\chi_2,\psi_2$, and $\Hom_{\GL_2(E)}(I(\chi_1,\psi_1),I(\chi_2,\psi_2))$ is non zero then either $\chi_1=\chi_2$ and $\psi_1=\psi_2$ or $\chi_1=\psi_2$ and $\psi_1=\chi_2$.
\end{thm}

Another tool we need to introduce is a pairing on $I(\chi,\psi) \times I(\chi^{-1},\psi^{-1})$ which identifies $I(\chi^{-1},\psi^{-1})$ with the smooth dual of $I(\chi,\psi)$. See \cite[Proposition 4.5.5]{bump}. The pairing is defined by an integral as follows
\begin{defiprop}\label{pairing}
The following integral defines a perfect pairing $\langle \ , \ \rangle: I(\chi,\psi) \times I(\chi^{-1},\psi^{-1}) \to \C$
\begin{displaymath}
\langle f_1 , f_2 \rangle:= \int_{\GL_2(\Oo)}f_1(h)f_2(h) dh,
\end{displaymath}
for every $f_1\in I(\chi,\psi), f_2\in I(\chi^{-1},\psi^{-1})$.
\end{defiprop}

\subsection{Whittaker models}

Let $\Psi$ be a fixed nontrivial additive character of $E$. 

\begin{defi}
Let $V$ be a smooth representations of $\GL_2(E)$. A \textit{Whittaker functional} on $V$ is a linear functional $\lambda: V \rightarrow \C$ satisfying
\begin{displaymath}
\lambda(m_x \cdot v)=\Psi(x)\lambda(v),
\end{displaymath}
for every $x\in F, v\in V$, where as above $m_x=\SmallMatrix{1 & x \\ 0 & 1}$.
\end{defi}

\begin{defi}\label{defwhitH}
We define $\mu:I(\chi,\psi)\rightarrow \C$ by
\begin{displaymath}
\mu(f)=\int_F f(w\cdot m_x) \Psi(-x)dx,
\end{displaymath}
where $w$ is as defined in the previous section.
\end{defi}
With the above notation, this integral converges if $\re(s_1-s_2)>0$, but we can proceed with an analytic continuation to every $s_1,s_2$ using flat sections as above. It is easy to verify that this is a Whittaker functional for $I(\chi,\psi)$. By \cite[Proposition 4.5.4]{bump}, every other Whittaker functional for this representation is scalar multiple of $\mu$.

\begin{defi}\label{defwhitmod} We attach to $I(\chi,\psi)$ a function $W=W_{\mu}$ defined by
\begin{displaymath}
W: f \mapsto \left(W_f: h \mapsto \mu(h\cdot f)= \int_F f(w\cdot m_x \cdot h) \Psi(-x) dx \right)
\end{displaymath}
It satisfies $W_f (m_x\cdot \SmallMatrix{z&0\\0&z} \cdot h)=\Psi(x)\cdot \chi\psi(z)\cdot W_f(h)$. 
\end{defi}

In the literature, one refers to \emph{Whittaker model} of a representation as the image of any such \emph{Whittaker function}. By abuse of notation, we will refer to the function itself as \emph{Whittaker model} of $I(\chi,\psi)$.

We can associate to every Whittaker functional $\lambda$ a \emph{Whittaker function} $W_{\lambda}$ as in the previous definition. The same result mentioned above tells us that they differ by a scalar.

\begin{rmk} One similarly defines a Whittaker model for every $(V,\pi)$ smooth representation of $\GL_2(E)$. For any Whittaker functional $\lambda$,  one lets $W_{\lambda}: v \mapsto (W_{\lambda,v}:h\mapsto \mu (h\cdot v))$. The image of $W_{\lambda}$ defines a subspace of the space of functions $\Lambda$ on $\GL_2(E)$, satisfying $\Lambda(m_x\cdot h)=\Psi(x)\cdot \Lambda(h)$. The group $\GL_2(E)$ acts naturally by right translation on this space, and the image of $W_{\lambda}$ is invariant for this action. Such image is isomorphic as $\GL_2(E)$-representation to $(V,\pi)$ and indeed provides a ``concrete'' model for $(V,\pi)$.
\end{rmk}

A straightforward computation yields the following property of the function defined above.
\begin{lemma}\label{lemmawhitH}For the Whittaker model of $I(\chi,\psi)$ we have
\begin{align*}
W_{\SmallMatrix{a & b \\ 0 & d}\cdot f}(\SmallMatrix{y & 0 \\ 0 &1}) = \Psi(d^{-1}by)\cdot\chi\psi(d)\cdot W_f(\SmallMatrix{d^{-1}ay & 0 \\ 0 &1}).
\end{align*}
%\begin{proof}
%This is straightforward rewriting
%\begin{displaymath}
%\SmallMatrix{y & 0 \\ 0 &1}\cdot \SmallMatrix{a & b \\ 0 & d}=m_{d^{-1}by}\cdot \SmallMatrix{d&0\\0&d} \cdot \SmallMatrix{d^{-1}ay&0\\0&1}
%\end{displaymath}
%and using $W_f (m_{d^{-1}by}\cdot \SmallMatrix{d&0\\0&d} \cdot \SmallMatrix{d^{-1}ay&0\\0&1})=\Psi(d^{-1}by)\cdot \chi\psi(d)\cdot W_f(\SmallMatrix{d^{-1}ay&0\\0&1})$
%\end{proof}
\end{lemma}

\subsection{Spherical representations}
Let $(V,\pi)$ an irreducible admissible representation of $\GL_2(E)$. One can consider the subspace $V^{\GL_2(\Oo)}$ of vectors fixed by the action of $\GL_2(\Oo)$. This is at most one dimensional (see \cite[Theorem 4.6.2]{bump}).

\begin{defi}
An irreducible admissible representation $(V,\pi)$ of $\GL_2(E)$ is called \textit{spherical} if it contains a $\GL_2(\Oo)-$fixed vector.  
\end{defi}

\begin{rmk}\label{rmkautom} The reason for which we are interested in spherical representations is that automorphic representations of $\GL_2$ decompose into a restricted product of local representations and these are all spherical outside a finite set of places. \end{rmk}

\begin{ex}[Principal series representations]\label{exsphericalvect} The representation $I(\chi,\psi)$ with $\chi,\psi$ unramified and $\chi\psi\neq|\cdot|^{\pm 1}$ is spherical. To see this, one define the \textit{normalised spherical vector} $\varphi_0$ as function on $\GL_2(E)$ by
\begin{displaymath}
\varphi_0 (h)=\varphi_0 (b\cdot k):=|a/d|^{1/2}\chi(a)\psi(d), \ \ \ \text{where }b=\SmallMatrix{a&*\\0&d}, k\in K=\GL_2(\Oo).
\end{displaymath}
To write $h\in \GL_2(E)$ as $b\cdot k$, we use again Iwasawa decomposition. It is clear by the definition that this function is fixed by $K$ and one verifies that it is well defined and that it is an element of $I(\chi,\psi)$. 
%Suppose that $bk=b'k'$ for $b=\SmallMatrix{a&*\\0&d},b'=\SmallMatrix{a'&*\\0&d'} k\in K=\GL_2(\Oo)$. Then we have $b=b'u$ with $u\in K \cap B(E)$ i.e. $u=\SmallMatrix{x&*\\0&y}$ with $x,y\in \Oo^{\times}$. Since $\chi(x)=\psi(y)=1$ and $|x|=|y|=1$, we find $\varphi_0(bk)=\varphi_0 (b'k')$. To check that this defines an element of $I(\chi,\psi)$ we  compute $\varphi_0(b'h)=\varphi(b'bk)$, where as before $h=bk$, $b=\SmallMatrix{a&*\\0&d},b'=\SmallMatrix{a'&*\\0&d'}, k\in K=\GL_2(\Oo)$. Hence
%\begin{displaymath}
%\varphi_0(b'h)=\varphi_0 (\SmallMatrix{aa'&*\\0&dd'}\cdot k)=|aa'/dd'|^{1/2}\chi(aa')\psi(dd')=|a'/d'|^{1/2}\chi(a')\psi(d')\varphi_0(h).
%\end{displaymath}
\end{ex}
It turns out that this example is enough to determine every spherical representation (of dimension greater than 1). 
\begin{thm}\label{sphericalprinc} \cite[Theorem 4.6.4]{bump}. Let $(V,\pi)$ be a spherical representation of $\GL_2(E)$ of dimension greater than $1$. Then $\pi$ is a spherical principal series representation.
\end{thm}

\begin{rmk}\label{sphericalprincrmk}
More precisely, a spherical representation $\pi$ will be isomorphic to $I_H(\chi,\psi)$, where $\chi,\psi$ are the unramified quasicharacters of $E^{\times}$ determined by 
\begin{displaymath}
\chi(\varpi)=\alpha, \ \ \psi(\varpi)=\beta,
\end{displaymath}
and $\alpha,\beta$ are the roots of the polynomial $X^2-q^{-1/2}\lambda X +\mu$, where $\lambda,\mu$ are the eigenvalues of $T(\p),R(\p)$ on the one-dimensional space of spherical vectors of $V$. Indeed the Hecke algebra of locally constant compactly supported complex valued functions on $\GL_2(E)$ acts on $(V,\pi)$ via the formula $\xi\cdot v=\int_{\GL_2(E)}\xi(g)( \pi(g)v) dg$ (see Definition \ref{defiheckealgebra}). The action of the subalgebra of $\GL_2(\Oo)$-biequivariant functions preserves the one-dimensional space $V^{\GL_2(\Oo)}$. In particular we can consider the eigenvalues for the action on spherical vectors of the operators
\begin{displaymath}
T(\p):=\ch(\GL_2(\Oo)\SmallMatrix{\varpi &0 \\0&1}\GL_2(\Oo)), \ \ \ R(\p):=\ch(\GL_2(\Oo)\SmallMatrix{\varpi &0 \\0&\varpi}\GL_2(\Oo)).
\end{displaymath}
\end{rmk}

We now want to characterise the Whittaker model of Definition \ref{defwhitH} for $\varphi_0\in I(\chi,\psi)$ as in the above example. First we let
\begin{displaymath}
\alpha:=\chi(\varpi), \ \ \ \beta:=\psi(\varpi).
\end{displaymath} 
We have the following result, that will be extremely helpful later. We write $W_0:=W_{\varphi_0}$.

\begin{thm}\label{whitvaluesH} \cite[Theorem 4.6.5]{bump}.
Let $\alpha,\beta$ as above. Then for any $y\in E^{\times}$, let $m:=\ord(y)$. We have
\begin{displaymath}
(1-q^{-1}\alpha\beta^{-1})^{-1}W_0\left( \SmallMatrix{y&0\\0&1}\right)=\begin{cases}
0 &\text{if } m<0 \\
q^{-m/2}\cdot \tfrac{\alpha^{m+1}-\beta^{m+1}}{\alpha-\beta} &\text{if } m\geq 0
\end{cases}
\end{displaymath}
\end{thm}

We want to work with a Whittaker model $\W$ such that for $\W_{\varphi_0}\left( \SmallMatrix{y&0\\0&1}\right)=1$ if $y\in \Oo^{\times}$. 
\begin{defi}\label{normwhitH}
The \textit{normalised Whittaker model of $I(\chi,\psi)$} is defined by $(1-q^{-1}\alpha\beta^{-1})^{-1}\cdot W$, for  $\alpha,\beta$ as above.
\end{defi} 

\subsection{Siegel sections}
This section contains exactly the same results and definitions of \cite[§3.2]{GSP4}. We report them for the seek of completeness and refer to \textit{loc. cit.} for the proofs. 

\begin{defi}
Let $\Ss(\Q_{\ell}^2,\C)$ be the space of Schwartz functions on $\Q_{\ell}^2$. For $\phi\in\Ss(\Q_{\ell}^2,\C)$, we write $\hat{\phi}$ for its Fourier transform, i.e.
\begin{displaymath}
\hat{\phi}(x,y)=\int_{\Q_{\ell}}\int_{\Q_{\ell}} e_{\ell}(xv-yu)\phi(u,v)du \ dv,
\end{displaymath}
where $e_{\ell}$ is the standard additive character on $F=\Q_{\ell}$, mapping $\ell^{-n}$ to $\exp(2\pi i/\ell^n)$.
\end{defi}

In the first part of \cite[Proposition 3.2.2]{GSP4}, the authors define a map from $\Ss(\Q_{\ell}^2,\C)$ to $I_H(\chi,\psi)$ for $\chi,\psi$ characters of $\Q_{\ell}^{\times}$ using explicit integrals. With the same notation we write
\begin{align*}
\Ss(\Q_{\ell}^2 &,\C) \rightarrow I_H(\chi,\psi) \\& \phi \mapsto f_{\phi,\chi,\psi}.
\end{align*}
\begin{prop}\label{propdefF}
The above mentioned map satisfies 
\begin{displaymath}
f_{g\cdot\phi,\chi,\psi}(h)=\chi(\det g)^{-1}|\det g|^{-1/2}f_{\phi,\chi,\psi}(hg),
\end{displaymath}
\begin{displaymath}
f_{\widehat{g\cdot\phi},\chi,\psi}(h)=\psi(\det g)^{-1}|\det g|^{-1/2}f_{\hat{\phi},\chi,\psi}(hg).
\end{displaymath}
In particular if $\psi=|\cdot|^{-1/2}$ and $\chi$ is unramified, then the map
\begin{align*}
\Ss(\Q_{\ell}^2 &,\C) \rightarrow I_H(\chi,\psi) \\& \phi \mapsto F_{\phi,\chi,\psi}:=f_{\hat{\phi},\chi,\psi}
\end{align*}
is $H(\Q_{\ell})$-equivariant.
\end{prop}

\begin{prop}\label{propinterchangschw} With notation as above, we have 
\begin{displaymath}
M(f_{\phi,\chi,\psi})=\frac{\varepsilon(\psi/\chi)}{L(\chi\psi^{-1}, 1)}\cdot f_{\hat{\phi},\psi,\chi},
\end{displaymath}
where $\varepsilon(\psi/\chi)$ is the local $\varepsilon$-factor (equal to 1 if $\psi/\chi$ is unramified).
\end{prop}

We now define some special Schwartz function that will be useful later. 

\begin{defi}\label{defischwt}
For integers $t\geq 0$ we define $\phi_t\in \Ss(\Q_{\ell}^2,\C)$ as follows
\begin{itemize}
\item for $t=0$, $\phi_0:=\ch(\Z_{\ell})\ch(\Z_{\ell})$,
\item for $t>0$, $\phi_t:=\ch(\ell^t\Z_{\ell})\ch(\Z_{\ell}^{\times})$.
\end{itemize}
This functions are preserved by the action of 
\begin{displaymath}
K_{H,0}(\ell^t):=\lbrace \mat \in H(\Z_{\ell}): c \equiv 0 \text{ mod } \ell^t \rbrace.
\end{displaymath}
\end{defi}

\begin{lemma}\label{lemmavalscw}
Let $\chi,\psi$ be unramified characters. The function $f_{\phi_t,\chi, \psi}$ is supported on $B(\Q_{\ell})K_{H,0}(\ell^t)$ and
\begin{displaymath}
f_{\phi_t,\chi, \psi}(1)=\begin{cases}
1 &\text{if } t=0 \\
L(\chi\psi^{-1},1)^{-1} &\text{if } t\geq 1.
\end{cases}
\end{displaymath}
\end{lemma}

\begin{defi}\label{schwrt1t}
For integers $t\geq 1$ we define $\phi_{1,t}\in \Ss(\Q_{\ell}^2,\C)$ to be $\ch(\ell^t\Z_{\ell})\ch(1+\ell^t\Z_{\ell})$. This function is preserved by the action of 
\begin{displaymath}
K_{H,1}(\ell^t):=\lbrace \gamma \in H(\Z_{\ell}): \gamma \equiv \SmallMatrix{*&*\\0&1} \text{ mod } \ell^t \rbrace.
\end{displaymath}
\end{defi}

\section{Zeta integrals for \texorpdfstring{$G(\Q_{\ell})=\GL_2(F_{\ell})$}{GL2(F)} (inert prime case)}\label{inertzeta}
Let $E$ be an unramified quadratic extension of $\Q_{\ell}$. We will work with the representation $\sigma=I_G(\chi,\psi)$ of $G=\GL_2(E)$. We denote with $K$ the subgroup $\GL_2(\Oo)$, where $\Oo$ is the ring of integers of $E$.

\subsection{Action of the Hecke algebra on Whittaker model}
First we recall the definition of Hecke algebra acting on $\sigma$.

\begin{defi}\label{defiheckealgebra}We denote with $\hh(G)$ the Hecke algebra of locally constant compactly supported $\C$-valued functions on $G=\GL_2(E)$. It is an algebra under convolution, defined by
\begin{displaymath}
\phi_1\star\phi_2 (g):= \int_G \phi_1(gh^{-1})\phi_2(h)dh,
\end{displaymath}
for $\phi_1,\phi_2\in \hh(G)$. Moreover we regard $\sigma$ as left $\hh(G)$- module via
\begin{displaymath}
\phi\cdot f = \int_G \phi(g)(g\cdot f) dg.
\end{displaymath}
\end{defi}

\begin{lemma}\label{lemmaheckealginert}
We have
\begin{displaymath}
g_1\cdot (\phi\cdot (g_2\cdot f))=\phi(g_1^{-1}(-)g_2^{-1})\cdot f.
\end{displaymath}
\end{lemma}

\begin{ex}[The operator $U(\ell)$]\label{Uellinert} We define $U(\ell)\in \hh(G)$ to be
\begin{displaymath}
U(\ell):=\tfrac{1}{\vol(K')}\ch(K'\cdot \SmallMatrix{\ell & 0 \\ 0 & 1}\cdot K'),
\end{displaymath}
where $K'$ is any subgroup of $K$ contained in $\lbrace \gamma \in K: \gamma \equiv \SmallMatrix{*&*\\0&1}$ mod $ \ell\rbrace$.
We can write $K'\cdot \SmallMatrix{\ell & 0 \\ 0 & 1}\cdot K'$ as union of left cosets
\begin{align*}
(K'\cdot \SmallMatrix{\ell & 0 \\ 0 & 1}\cdot K')&=\bigsqcup_{\gamma\in J }\gamma\cdot \SmallMatrix{\ell & 0 \\ 0 & 1}\cdot K'
\end{align*}
where $J$ is a set of representatives for the left quotient $(\SmallMatrix{\ell & 0 \\ 0 & 1}K'\SmallMatrix{\ell^{-1} & 0 \\ 0 & 1}\cap K')\backslash K'$. We claim that we can take $J=\lbrace \SmallMatrix{1 & u \\ 0 & 1} \rbrace_{u\in \left(\Oo/\ell\Oo\right)}$. Indeed the subgroup for which we are taking the quotient is the subgroup $K''$ of matrices of $\mat\in K'$ such that $c\equiv 0$ mod $\ell$. The matrices considered are clearly in distinct cosets and since for any $\mat\in K'$, $d\not\equiv 0$ mod $\ell$ we can choose $u \in \Oo/\ell\Oo$ such that $b\equiv u d$ mod $\ell$. In other words
\begin{displaymath}
\SmallMatrix{1 & u \\ 0 & 1}^{-1}\mat\in K''.
\end{displaymath}
Hence we can rewrite  
\begin{displaymath}
(K'\cdot \SmallMatrix{\ell & 0 \\ 0 & 1}\cdot K')=\bigsqcup_{u\in \left(\Oo/\ell\Oo\right)}\SmallMatrix{1 & u \\ 0 & 1}\cdot \SmallMatrix{\ell & 0 \\ 0 & 1}\cdot K' \\= \bigsqcup_{u\in \left(\Oo/\ell\Oo\right)}\SmallMatrix{\ell & u \\ 0 & 1}\cdot K'.
\end{displaymath}
\end{ex}

We will need to define an appropriate additive character of $E$ and then work with the normalised Whittaker model for $\sigma$ as in Definition \ref{normwhitH}. Consider $e_{\ell}$ the standard additive character on $\Q_{\ell}$. We define an additive character $\Psi$ on $E$ fixing $\delta\in\Oo_E$ such that $E=\Q_{\ell}\oplus \Q_{\ell}(\delta)$ and letting
\begin{displaymath}
\Psi: x \to e_{\ell}(\tr_{E/\Q_{\ell}}(\delta^{-1}x)).
\end{displaymath}
We can assume $v(\delta)=0$ since $E/\Q_{\ell}$ is unramified. This character has conductor $\Oo_E$ (see for example \cite[Exercise 3(e), Chapter 7]{fourier}).

We describe how the action of the operator $U(\ell)$ of Example \ref{Uellinert} modifies the Whittaker model.

\begin{prop}\label{actionUellWhittinert}
Let $\varphi\in\sigma$ a spherical vector. Then for any $y\in E^{\times}$, we have
\begin{displaymath}
\W_{U(\ell)\cdot \varphi}\left( \SmallMatrix{y&0\\0&1}\right)=\begin{cases}
0 &\text{ if } v(y)< 0 \\
\ell^2\W_{\varphi}\left( \SmallMatrix{\ell y&0\\0&1}\right) &\text{ if } v(y)\geq 0.
\end{cases}
\end{displaymath}
\begin{proof} We prove the result for $W$ as in Definition \ref{defwhitmod}. We can also assume $\phi=\phi_0$ the normalised spherical vector. By definition
\begin{align*}
W_{U(\ell)\cdot \varphi_0}\left( \SmallMatrix{y&0\\0&1}\right)&=\mu\left( \SmallMatrix{y&0\\0&1} \cdot (U(\ell)\cdot\varphi_0)\right)=\sum_{u\in\left(\Oo_E/\ell\Oo_E\right)}\mu\left(\SmallMatrix{y&0\\0&1}\SmallMatrix{\ell & u \\ 0 & 1}\cdot \varphi_0\right),
\end{align*}
where in the second equality we used the decomposition of $U(\ell)$ as in Example \ref{Uellinert} and the fact that $\varphi_0$ is $K$-invariant. Now we write
\begin{displaymath}
\SmallMatrix{y&0\\0&1}\SmallMatrix{\ell & u \\ 0 & 1}=\SmallMatrix{1&yu\\0&1}\SmallMatrix{y\ell & 0 \\ 0 & 1}=m_{yu}\cdot \SmallMatrix{y\ell & 0 \\ 0 & 1}.
\end{displaymath}
So we find 
\begin{displaymath}
W_{U(\ell)\cdot \varphi_0}\left(\SmallMatrix{y&0\\0&1}\right)=\sum_{u\in\left(\Oo_E/\ell\Oo_E\right)}\Psi(yu)W_{\varphi_0} \left( \SmallMatrix{\ell y&0\\0&1}\right).
\end{displaymath}
If $v(\ell y)<0$, i.e. $v(y)<-1$, applying Theorem \ref{whitvaluesH}, we find that the above quantity is zero. If $v(\ell y)=0$, i.e. $v(y)=-1$, the sum is equal to 
\begin{displaymath}
\sum_{u\in\left(\Oo_E/\ell\Oo_E\right)}\Psi(yu)=\sum_{0\leq i,j\leq \ell -1}\Psi(y(i+\delta j))=\sum_{0\leq i,j\leq \ell -1}\Psi(iy)\Psi(\delta jy))=\sum_{0\leq i,j\leq \ell -1}e_{\ell}(\tr(\delta^{-1} y))^ie_{\ell}(\tr(y))^j.
\end{displaymath}
Having assumed that $v(\delta)=0$ and having $v(y)=-1$, we have that at least one of the two terms $e_{\ell}(\tr(\delta^{-1} y)),e_{\ell}(\tr(y))$ is equal to $\zeta_{\ell}=\exp^{{2\pi i}/{\ell}}$. Assume for example $e_{\ell}(\tr(\delta^{-1} y))=\zeta_{\ell}$, we can rewrite the sum as
\begin{displaymath}
\sum_{0\leq j\leq \ell -1}e_{\ell}(\tr(y))^j\cdot \left(\sum_{0\leq i\leq \ell -1}\zeta_{\ell}^i\right)=0.
\end{displaymath}
Finally, if $v(y\ell)>0$, i.e. $y\in \Oo_F$, $\Psi(yu)=1$ and hence
\begin{displaymath}
W_{U(\ell)\cdot \varphi_0}\left( \SmallMatrix{y&0\\0&1}\right)=\ell^2 W_{\varphi_0}\left( \SmallMatrix{\ell y&0\\0&1}\right).
\end{displaymath}
Hence the result.
\end{proof}
\end{prop}

\subsection{Zeta integrals}\label{lfnctinert}

As above fix the irreducible spherical principal series representation $\sigma=I_G(\chi,\psi)$, for $\chi,\psi$ quasicharacters of $E^{\times}$. Let
\begin{displaymath}
\alpha:=\chi(\ell), \ \ \ \beta=\psi(\ell)
\end{displaymath}
and let $\chi_{\sigma}$ the central character of $\sigma$, i.e. $\chi_{\sigma}=\chi\psi$. We define the \textit{local Asai $L$-factor}\footnote{The standard $L$-factor of $\sigma$ is $[(1-\alpha \ell^{-2s})(1-\beta\ell^{-2s})]^{-1}$ and can be obtained by the same integral we consider here, but integrating over $y\in E^{\times}$ with norm and measure on $E$ rather than on $\Q_{\ell}$. It will be clear later the reason of the name Asai $L$-factor.} of $\sigma$ to be
\begin{displaymath}
L(\as(\sigma),s):=[(1-\alpha \ell^{-s})(1-\beta\ell^{-s})(1-\alpha\beta\ell^{-2s})]^{-1},
\end{displaymath}
Moreover if $\eta$ is an unramified character of $\Q_{\ell}^{\times}$, we let
\begin{displaymath}
L(\as(\sigma\otimes \eta),s)=[(1-\alpha\eta(\ell) \ell^{-s})(1-\beta\eta(\ell)\ell^{-s})(1-\alpha\beta\eta(\ell)^2\ell^{-2s})]^{-1}.
\end{displaymath} 
\begin{defi}
Let $\sigma$ as above and $\eta$ an unramified character of $\Q_{\ell}^{\times}$. For every $f\in \sigma$, we define
\begin{displaymath}
Z(\sigma,\eta,f,s):=L(\as(\sigma\otimes \eta),s)^{-1}\int_{\Q_{\ell}^{\times}}|y|^{s-1}\eta(y)\W_f\left( \SmallMatrix{y&0\\0&1}\right)d^{\times}y.
\end{displaymath}
\end{defi}

The following three lemmas will be very useful.
\begin{lemma}[Zeta integral at the spherical vector]\label{sphericalzetainert}
There exist $r(\sigma,\eta)\in\R$ such that for ever $f \in \sigma$ and $s\in \C$ such that $\re(s)>r(\sigma,\eta)$, the above integral is absolutely convergent and, as function of $s$, lies in $\C[\ell^s,\ell^{-s}]$; in particular it has analytic continuation for all $s\in\C$. Moreover, if $\varphi_0$ is the normalised spherical vector as above, we have
\begin{displaymath}
Z(\sigma,\eta,\varphi_0,s)=L(\eta^2\chi_{\sigma},2s)^{-1}.
\end{displaymath}
\begin{proof}
It is enough to check convergence and analytic continuation for $f=g\cdot \varphi_0$, where $g\in G$. The validity of these statements for such $f$ depends only on the class of $g$ in $N\backslash G/G(\Oo_E)$. Since representatives of this quotient are elements of the form $\SmallMatrix{a&0\\0&d}$, Lemma \ref{lemmawhitH} implies that it suffices to look at the integral for $f=\varphi_0$.

Applying Theorem \ref{whitvaluesH} (notice that in our case $q=\ell^2$), we find
\begin{align*}
\int_{\Q_{\ell}^{\times}}|y|^{s-1}\eta(y)\W_{\varphi_0}\left( \SmallMatrix{y&0\\0&1}\right)d^{\times}y&=\sum_{m\geq 0}(\ell^{s-1})^{-m}\ell^{-m}\eta(\ell)^m \cdot \frac{\alpha^{m+1}-\beta^{m+1}}{\alpha-\beta}&=\sum_{m\geq 0} X^m \cdot\frac{\alpha^{m+1}-\beta^{m+1}}{\alpha-\beta},
\end{align*}
where $X=\ell^{-s}\eta(\ell)$. We can manipulate the latter series and obtain
\begin{align*}
\sum_{m\geq 0} X^m \cdot\frac{\alpha^{m+1}-\beta^{m+1}}{\alpha-\beta} & = \frac{1}{\alpha-\beta}\sum_{m\geq 0} \left( \alpha \cdot (X\alpha)^m -\beta (X\beta)^m\right) = \frac{1}{\alpha-\beta} \left( \frac{\alpha}{1-\alpha X} - \frac{\beta}{1-\beta X}\right)\\&=\frac{1}{(1-\alpha X)(1-\beta X)}.
\end{align*}
The series converges for $|\alpha X|_{\C}<1,|\beta X|_{\C}<1$, that is for $\re(s)>r(\sigma,\eta)$, for some real number depending on $\sigma$ and $\eta$. Substituting $X=\ell^{-s}\eta(\ell)$, for $s$ in this region, we find 
\begin{align*}
\int_{\Q_{\ell}^{\times}}|y|^{s-1}\eta(y)\W_{\varphi_0}\left( \SmallMatrix{y&0\\0&1}\right)d^{\times}y&= [(1-\alpha\eta(\ell)\ell^{-s})(1-\beta\eta(\ell)\ell^{-s})]^{-1}=L(\as(\sigma\otimes \eta), s)\cdot (1-\alpha\beta\eta(\ell)^{2}\ell^{-2s})
\end{align*}
and $(1-\alpha\beta\eta(\ell)^{2}\ell^{-2s})=L(\chi_{\sigma}\eta^2,2s)^{-1}$. 
\end{proof}
\end{lemma}

\begin{lemma}[Action of $U(\ell)$ on the zeta integral]\label{zetaU(l)inert}
If $\varphi_0$ is the normalised spherical vector as above, we have
\begin{displaymath}
Z(\sigma,\eta,U(\ell)\varphi_0,s)=\tfrac{\ell^{s+1}}{\eta(\ell)}[Z(\sigma,\eta,\varphi_0,s)-L(\as(\sigma\otimes\eta),s)^{-1}]=\tfrac{\ell^{s+1}}{\eta(\ell)}[L(\eta^2\chi_{\sigma},2s)^{-1}-L(\as(\sigma\otimes\eta),s)^{-1}]
\end{displaymath}
\begin{proof}
First we apply Proposition \ref{actionUellWhittinert} and find
\begin{align*}
Z(\sigma,\eta,U(\ell)\varphi_0,s)&=\ell^{2}L(\as(\sigma\otimes\eta),s)^{-1}\int_{|y|< \ell}|y|^{s-1}\eta(y)\W_0\left(\SmallMatrix{\ell y & 0 \\ 0 &1}\right)d^{\times}y\\&=\ell^{s+1}\eta(\ell)^{-1}L(\as(\sigma\otimes\eta),s)^{-1}\int_{|y|< 1}|y|^{s-1}\eta(y)\W_0\left(\SmallMatrix{y & 0 \\ 0 &1}\right)d^{\times}y,
\end{align*}
where in the second equality we used the change of variables $y \rightsquigarrow \ell y$. We then rewrite the integral in the last term as
\begin{align*}
\int_{\Q_{\ell}^{\times}}|y|^{s-1}\eta(y)\W_0\left(\SmallMatrix{ y & 0 \\ 0 &1},\SmallMatrix{y & 0 \\ 0 &1}\right)d^{\times}y-\int_{|y|\geq 1}|y|^{s-1}\eta(y)\W_0\left(\SmallMatrix{ y & 0 \\ 0 &1},\SmallMatrix{y & 0 \\ 0 &1}\right)d^{\times}y.
\end{align*}
Then we apply Theorem \ref{whitvaluesH} and obtain
\begin{align*}
\int_{|y|\geq 1}|y|^{s-1}\eta(y)\W_0\left(\SmallMatrix{y & 0 \\ 0 &1}\right)d^{\times}y=\int_{\Z_{\ell}^{\times}}\W_0\left(\SmallMatrix{y & 0 \\ 0 &1}\right)d^{\times}y=\int _{\Z_{\ell}^{\times}} d^{\times}y = 1.
\end{align*}
Putting everything together we find 
\begin{align*}
Z(\sigma,\eta,U(\ell)\varphi_0,s)&=\ell^{s+1}\eta(\ell)^{-1}L(\as(\sigma\otimes\eta),s)^{-1}\left( \int_{\Q_{\ell}^{\times}}|y|^{s-1}\eta(y)\W_0\left(\SmallMatrix{ y & 0 \\ 0 &1},\SmallMatrix{y & 0 \\ 0 &1}\right)d^{\times}y- 1\right)\\&=\ell^{s+1}\eta(\ell)^{-1}(Z(\sigma,\eta,\varphi_0,s)-L(\as(\sigma\otimes\eta),s)^{-1}).
\end{align*}
\end{proof}
\end{lemma}

\begin{lemma}[Action of the Borel subgroup of $\GL_2(\Q_{\ell})$]\label{borelHinert}
For any $f\in\sigma$, $a,d\in\Q_{\ell}^{\times}$, we have
\begin{displaymath}
Z(\sigma,\eta,\SmallMatrix{a&*\\0&d}\cdot f,s)=\left|\tfrac{d}{a}\right|^{s-1}\chi_{\sigma}(d)\eta(a^{-1}d) \cdot Z(\sigma,\eta,f,s)
\end{displaymath}
\begin{proof}
We apply Lemma \ref{lemmawhitH} together with the fact that, for our choice of $\Psi$, we have $\Psi(x)=1$ for every $x\in\Q_{\ell}$. We find
\begin{align*}
Z(\sigma,\eta,\SmallMatrix{a&*\\0&d}\cdot f,s)&=\chi_{\sigma}(d)L(\as(\sigma\otimes\eta),s)^{-1}\int_{\Q_{\ell}^{\times}}|y|^{s-1}\eta(y)\W_f\left(\SmallMatrix{d^{-1}ay & 0 \\ 0 &1}\right)d^{\times}y \\&=\chi_{\sigma}(d)|d/a|^{s-1}\eta(a^{-1}d)L(\as(\sigma\otimes\eta),s)^{-1}\int_{\Q_{\ell}^{\times}}|y|^{s-1}\eta(y)\W_f\left(\SmallMatrix{y & 0 \\ 0 &1}\right)d^{\times}y\\&=\chi_{\sigma}(d)|d/a|^{s-1}\eta(a^{-1}d)\cdot Z(\sigma,\eta,f,s)
\end{align*}
where in the second equality we used the change of variable $y \rightsquigarrow d^{-1}ay$.
\end{proof}
\end{lemma}

\section{Zeta integrals for \texorpdfstring{$G(\Q_{\ell})=\GL_2(\Q_{\ell})\times\GL_2(\Q_{\ell})$}{GL2xGL2} (split prime case)}\label{splitzeta}
\subsection{Whittaker models for \texorpdfstring{$G=\GL_2\times\GL_2$}{GL2xGL2}}
Let $\chi_1,\psi_1,\chi_2,\psi_2$ be quasicharacters of $\Q_{\ell}^{\times}$. We now consider a representation of $G=\GL_2(\Q_{\ell})\times\GL_2(\Q_{\ell})$.
\begin{defi} For $\chi_1,\psi_1,\chi_2,\psi_2$ as above, let
\begin{displaymath}
I_G(\underline{\chi},\underline{\psi}):= I_H(\chi_1,\psi_1)\otimes I_H(\chi_2,\psi_2),
\end{displaymath}
i.e. $f\in I_G(\underline{\chi},\underline{\psi})$ is $f:G\rightarrow \C$ such that 
\begin{displaymath}
f((\SmallMatrix{a & * \\ 0 & d},\SmallMatrix{a' & * \\ 0 & d'})\cdot g)=  \left|\tfrac{a}{d}\right|^{1/2}\left|\tfrac{a'}{d'}\right|^{1/2}\chi_1(a)\psi_1(d)\chi_2(a')\psi_2(d')f(g).
\end{displaymath}
\end{defi}

We see $I_G(\underline{\chi},\underline{\psi})$ as a $G$-representation letting $G$ act by right translation.

\begin{defi}\label{sphgl2gl2}
The $G$-representations $I_G(\underline{\chi},\underline{\psi})$, for $\chi_1,\chi_2,\psi_1,\psi_2$ quasicharacters of $\Q_{\ell}^{\times}$ are called \textit{principal series representations for $G$}.
\end{defi}

We need to define what is a Whittaker functional for a representation $V$ of $G$, having fixed an additive character $\Psi$.
\begin{defi}
A \textit{Whittaker functional} on $V$ is a linear functional $\lambda: V \rightarrow \C$ satisfying
\begin{displaymath}
\lambda(m_{x,x'} \cdot v)=\Psi(x-x')\lambda(v),
\end{displaymath}
for every $x,x'\in \Q_{\ell}, v\in V$, where $m_{x,x'}=(m_x,m_{x'})=(\SmallMatrix{1 & x \\ 0 & 1},\SmallMatrix{1 & x' \\ 0 & 1})$.
\end{defi}

We now define a Whittaker model for $\sigma=I_G(\underline{\chi},\underline{\psi})$, for $\underline{\chi}=(\chi_1,\chi_2),\underline{\psi}=(\psi_1,\psi_2)$. We will be using Whittaker models for $I_H(\chi_1,\psi_1)$ and $I_H(\chi_2,\psi_2)$ as constructed above. Recall that everything depends on the choice of the additive character. We will consider the functionals as in Definition \ref{defwhitH}, but with different choices of the additive character. Fix such an additive character $\Psi$ for which we want to obtain a  Whittaker functional for $\sigma$. We then let $\Psi_1=\Psi$ and $\Psi_2=\Psi(-(\cdot))$. And write $\mu_i:I(\chi_i,\psi_i)\rightarrow \C$  where
\begin{displaymath}
\mu_1(f_1)=\int_F f_1(w\cdot m_x) \Psi_1(-x)dx,
\end{displaymath}
\begin{displaymath}
\mu_2(f_2)=\int_F f_2(w\cdot m_x) \Psi_2(-x)dx.
\end{displaymath}
And finally let $\mu: \sigma \rightarrow \C$ to be defined by
\begin{displaymath}
\mu(f_1,f_2)=\mu_1(f_1)\cdot\mu_2(f_2).
\end{displaymath}
It is straightforward to see that it is a Whittaker functional for $\sigma$.

\begin{defi} We let $W$ be the Whittaker model for $\sigma$ defined by
\begin{displaymath}
W: f \mapsto \left(W_f: g=(g_1,g_2) \mapsto \mu(g\cdot f) \right)
\end{displaymath}
\end{defi}

From the definition we have, for $f=f_1\otimes f_2$,
\begin{align*}
W_f(g_1,g_2)&=\mu_1(g_1\cdot f_1)\cdot\mu(g_2\cdot f_2)=\left(\int_F f_1(w\cdot m_x\cdot g_1) \Psi_1(-x)dx \right)\cdot \left(\int_F f_2(w\cdot m_x\cdot g_2) \Psi_2(-x)dx\right)\\&=W_{1,f_1}(g_1)\cdot W_{2,f_2}(g_2),
\end{align*}
where $W_{1,f_1},W_{2,f_2}$ are the Whittaker models for $I_H(\chi_1,\psi_1)$ and $I_H(\chi_2,\psi_2)$ obtained from the functionals $\mu_1,\mu_2$.
\begin{lemma}\label{lemmaactionborelH}For the Whittaker  model of $\sigma$ we have
\begin{align*}
W_{\left(\SmallMatrix{a & b \\ 0 & d},\SmallMatrix{a & b \\ 0 & d}\right)\cdot f}\left(\SmallMatrix{y & 0 \\ 0 &1},\SmallMatrix{y & 0 \\ 0 &1}\right) =\chi_{\sigma}(d)\cdot W_f\left(\SmallMatrix{d^{-1}ay & 0 \\ 0 &1},\SmallMatrix{d^{-1}ay & 0 \\ 0 &1}\right),
\end{align*}
where $\chi_{\sigma}=\chi_1\psi_1\chi_2\psi_2$ will be called the central character of $\sigma$.
\begin{proof}
This is straightforward from Lemma \ref{lemmawhitH}. Indeed, by definition, the left hand side term is equal to 
\begin{align*}
W_{1,\SmallMatrix{a & b \\ 0 & d}\cdot f_1}(\SmallMatrix{y & 0 \\ 0 &1})\cdot W_{2,\SmallMatrix{a & b \\ 0 & d}\cdot f_2}(\SmallMatrix{y & 0 \\ 0 &1})
\end{align*}
and applying the lemma, this is equal to
\begin{displaymath}
\Psi(d^{-1}by)\cdot \Psi(-d^{-1}by)\cdot\chi_1\psi_1(d)\cdot\chi_2\psi_2(d)\cdot W_{1,f_1}(\SmallMatrix{d^{-1}ay & 0 \\ 0 &1})\cdot W_{2,f_2}(\SmallMatrix{d^{-1}ay & 0 \\ 0 &1}).
\end{displaymath}
\end{proof}
\end{lemma}

\begin{defi}
The \textit{normalised Whittaker model for $\sigma$} is defined by
\begin{displaymath}
\W=(1-\ell^{-1}\alpha\beta^{-1})^{-1}(1-\ell^{-1}\gamma\delta^{-1})^{-1}\cdot W
\end{displaymath}
where $\alpha=\chi_1(\ell), \beta=\psi_1(\ell),\gamma=\chi_2(\ell), \delta=\psi_2(\ell)$.
\end{defi}

The definition and properties of spherical representations of $H$ carry over to representations of $G$, using the subgroup $\GL_2(\Z_{\ell})\times \GL_2(\Z_{\ell})$. In particular we define the \textit{normalised spherical vector of $\sigma$} to be 
\begin{displaymath}
\varphi_0= \varphi_{1,0}\otimes\varphi_{2,0},
\end{displaymath}
where $\varphi_{i,0}$ is the normalised spherical vector for $I_H(\chi_i,\psi_i)$ as in Example \ref{exsphericalvect}. Let then 
\begin{displaymath}
\W_0:=\W_{\varphi_0}.
\end{displaymath}
\begin{thm}\label{whitvaluesG}
Let, as above, $\alpha=\chi_1(\ell), \beta=\psi_1(\ell),\gamma=\chi_2(\ell), \delta=\psi_2(\ell)$. Then for any $y\in \Q_{\ell}^{\times}$, let $m:=\ord(y)$. We have
\begin{displaymath}
\W_0\left( \SmallMatrix{y&0\\0&1},\SmallMatrix{y&0\\0&1}\right)=\begin{cases}
0 &\text{if } m<0 \\
\ell^{-m}\cdot \tfrac{\alpha^{m+1}-\beta^{m+1}}{\alpha-\beta}\cdot \tfrac{\gamma^{m+1}-\delta^{m+1}}{\gamma-\delta} &\text{if } m\geq 0
\end{cases}
\end{displaymath}
\begin{proof}
This is a corollary of Theorem \ref{whitvaluesH}
\end{proof}
\end{thm}

\subsection{Action of the Hecke algebra on Whittaker model}
We will now recall the definition of the Hecke algebra acting on $\sigma=I_G(\underline{\chi},\underline{\psi})$.

\begin{defi}We denote with $\hh(G)$ the Hecke algebra of locally constant compactly supported $\C$-valued functions on $G$. It is an algebra under convolution, defined by
\begin{displaymath}
\phi_1\star\phi_2 (g):= \int_G \phi_1(gh^{-1})\phi_2(h)dh,
\end{displaymath}
for $\phi_1,\phi_2\in \hh(G)$. Moreover we regard $\sigma$ as left $\hh(G)$- module via
\begin{displaymath}
\phi\cdot f = \int_G \phi(g)(g\cdot f) dg.
\end{displaymath}
\end{defi}

\begin{lemma}\label{lemmaheckealg}
We have
\begin{displaymath}
g_1\cdot (\phi\cdot (g_2\cdot f))=\phi(g_1^{-1}(-)g_2^{-1})\cdot f.
\end{displaymath}
\end{lemma}

\begin{ex}[The operator $U(\ell)$]\label{Uell} We define $U(\ell)\in \hh(G)$ to be, essentially $(U(\ell),U(\ell))$, i.e. the usual $U(\ell)$ operator on each of the $\GL_2(\Q_{\ell})$. More precisely
\begin{displaymath}
U(\ell):=\tfrac{1}{\vol(K')}\ch(K'\cdot \left( \SmallMatrix{\ell & 0 \\ 0 & 1}, \SmallMatrix{\ell & 0 \\ 0 & 1}\right)\cdot K')
\end{displaymath}
for $K'=K'_1\times K'_2$ subgroup of $\GL_2(\Z_{\ell})\times \GL_2(\Z_{\ell})$, with $K_1',K_2'\subset \lbrace \gamma \in K: \gamma \equiv \SmallMatrix{*&*\\0&1}$ mod $ \ell\rbrace$. Proceeding as in Example \ref{Uellinert}, we can rewrite
\begin{align*}
K'\cdot \left( \SmallMatrix{\ell & 0 \\ 0 & 1}, \SmallMatrix{\ell & 0 \\ 0 & 1}\right)\cdot K'&=\bigsqcup_{0\leq u,v\leq\ell-1}\left(\SmallMatrix{1 & u \\ 0 & 1},\SmallMatrix{1& v \\ 0 & 1}\right)\cdot\left(\SmallMatrix{\ell & 0 \\ 0 & 1},\SmallMatrix{\ell & 0 \\ 0 & 1}\right)\cdot K' \\&= \bigsqcup_{0\leq u,v\leq\ell-1}\left(\SmallMatrix{\ell & u \\ 0 & 1},\SmallMatrix{\ell& v \\ 0 & 1}\right)\cdot K'.
\end{align*}
\end{ex}

From now on we take $\Psi=e_{\ell}$, the standard additive character of $\Q_{\ell}$, i.e. the one mapping $\ell^{-n}$ to $\exp(2\pi i/\ell^n)$.

We describe how the action of the operator $U(\ell)$ of Example \ref{Uell} modifies the Whittaker model.

\begin{prop}\label{actionUellWhitt}
Let $\varphi\in\sigma$ a spherical vector. Then for any $y\in \Q_{\ell}^{\times}$, we have
\begin{displaymath}
\W_{U(\ell)\cdot \varphi}\left( \SmallMatrix{y&0\\0&1},\SmallMatrix{y&0\\0&1}\right)=\begin{cases}
0 &\text{ if } |y|\geq \ell \\
\ell^2\W_{\varphi}\left( \SmallMatrix{\ell y&0\\0&1},\SmallMatrix{\ell y&0\\0&1}\right) &\text{ if } |y|< \ell.
\end{cases}
\end{displaymath}
\begin{proof} We prove the result for the Whittacker model $W$. By definition
\begin{align*}
W_{U(\ell)\cdot \varphi}\left( \SmallMatrix{y&0\\0&1},\SmallMatrix{y&0\\0&1}\right)&=\mu\left( \left( \SmallMatrix{y&0\\0&1},\SmallMatrix{y&0\\0&1}\right)\cdot (U(\ell)\cdot\varphi)\right)\\&=\sum_{0\leq u,v\leq\ell-1}\mu\left(\left( \SmallMatrix{y&0\\0&1}\SmallMatrix{\ell & u \\ 0 & 1},\SmallMatrix{y&0\\0&1}\SmallMatrix{\ell& v \\ 0 & 1}\right)\cdot \varphi\right),
\end{align*}
where in the second equality we used the decomposition of $U(\ell)$ as in Example \ref{Uell} and the fact that $\varphi$ is $K$-invariant. Now we write
\begin{displaymath}
\left( \SmallMatrix{y&0\\0&1}\SmallMatrix{\ell & u \\ 0 & 1},\SmallMatrix{y&0\\0&1}\SmallMatrix{\ell& v \\ 0 & 1}\right)=\left( \SmallMatrix{1&yu\\0&1}\SmallMatrix{y\ell & 0 \\ 0 & 1},\SmallMatrix{1&yv\\0&1}\SmallMatrix{y\ell& 0 \\ 0 & 1}\right)=m_{yu,yv}\cdot\left( \SmallMatrix{y\ell & 0 \\ 0 & 1},\SmallMatrix{y\ell& 0 \\ 0 & 1}\right).
\end{displaymath}
So we find 
\begin{displaymath}
W_{U(\ell)\cdot \varphi}\left( \SmallMatrix{y&0\\0&1},\SmallMatrix{y&0\\0&1}\right)=\sum_{0\leq u,v\leq\ell-1}\Psi(yu)\Psi(-yv)W_{\varphi} \left( \SmallMatrix{\ell y&0\\0&1},\SmallMatrix{\ell y&0\\0&1}\right).
\end{displaymath}
If $|y\ell|>1$, applying Theorem \ref{whitvaluesG}, we find that the above quantity is zero. Similarly if $|y\ell|=1$,i.e. $e_{\ell}(y)=\zeta_{\ell}:=e^{2\pi i/\ell}$, the sum is equal to 
\begin{displaymath}
c\cdot\sum_{0\leq u,v\leq\ell-1}e_{\ell}(yu)e_{\ell}(-yv)=c\cdot \sum_{0\leq u,v\leq\ell-1}\zeta_{\ell}^u\zeta_{\ell}^{-v}=0,
\end{displaymath}
where $\varphi=c\cdot \varphi_0$.
Finally, if $|y\ell|<1$, $e_{\ell}(yu)=e_{\ell}(-yv)=1$ for every $u,v$ and hence
\begin{displaymath}
W_{U(\ell)\cdot \varphi}\left( \SmallMatrix{y&0\\0&1},\SmallMatrix{y&0\\0&1}\right)=\ell^2 W_{\varphi}\left( \SmallMatrix{\ell y&0\\0&1},\SmallMatrix{\ell y&0\\0&1}\right).
\end{displaymath}
Hence the result.
\end{proof}
\end{prop}

\subsection{Zeta integrals}\label{lfnctsplit}

We fix the quasicharacters $\chi_1,\psi_1,\chi_2,\psi_2$ such that $\chi_1\psi_1^{-1}\neq |\cdot |^{\pm 1}, \chi_2\psi_2^{-1}\neq |\cdot |^{\pm 1}$ and then fix the irreducible spherical principal series representation $\sigma=I_G(\underline{\chi},\underline{\psi})$ as above. We define the \textit{local $L$-factor} of $\sigma$ to be
\begin{displaymath}
L(\sigma,s):=[(1-\alpha\gamma\ell^{-s})(1-\alpha\delta\ell^{-s})(1-\beta\gamma\ell^{-s})(1-\beta\delta\ell^{-s})]^{-1},
\end{displaymath}
where, as above, $\alpha=\chi_1(\ell),\beta=\psi_1(\ell),\gamma=\chi_2(\ell),\delta=\psi_2(\ell)$.
Moreover if $\eta$ is an unramified character of $\Q_{\ell}^{\times}$, we define
\begin{displaymath}
L(\sigma\otimes \eta,s)=[(1-\alpha\gamma\eta(\ell)\ell^{-s})(1-\alpha\delta\eta(\ell)\ell^{-s})(1-\beta\gamma\eta(\ell)\ell^{-s})(1-\beta\delta\eta(\ell)\ell^{-s})]^{-1}.
\end{displaymath} 
\begin{defi}
Let $\sigma$ as above and $\eta$ an unramified character of $\Q_{\ell}^{\times}$. For every $f\in \sigma$, we define
\begin{displaymath}
Z(\sigma,\eta,f,s):=L(\sigma\otimes \eta,s)^{-1}\int_{\Q_{\ell}^{\times}}|y|^{s-1}\eta(y)\W_f\left( \SmallMatrix{y&0\\0&1},\SmallMatrix{y&0\\0&1}\right)d^{\times}y.
\end{displaymath}
\end{defi}

The following three useful lemmas are the analogues of Lemmas \ref{sphericalzetainert}, \ref{zetaU(l)inert}, \ref{borelHinert} of the previous section.
\begin{lemma}[Zeta integral at the spherical vector]\label{sphericalzeta}
There exist $r(\sigma,\eta)\in\R$ such that for ever $f \in \sigma$ and $s\in \C$ such that $\re(s)>r(\sigma,\eta)$, the above integral is absolutely convergent and, as function of $s$, lies in $\C[\ell^s,\ell^{-s}]$; in particular it has analytic continuation for all $s\in\C$. Moreover, if $\varphi_0$ is the normalised spherical vector as above, we have
\begin{displaymath}
Z(\sigma,\eta,\varphi_0,s)=L(\eta^2\chi_{\sigma},2s)^{-1}.
\end{displaymath}
\begin{proof}
In order to prove the first statements, we reduce to compute the integral for $f=\varphi_0$, arguing as in the proof of Lemma \ref{sphericalzetainert}. Applying Theorem \ref{whitvaluesG}, we find
\begin{align*}
\int_{\Q_{\ell}^{\times}}|y|^{s-1}\eta(y)\W_{\varphi_0}\left( \SmallMatrix{y&0\\0&1},\SmallMatrix{y&0\\0&1}\right)d^{\times}y&=\sum_{m\geq 0}(\ell^{s-1})^{-m}\ell^{-m}\eta(\ell)^m \cdot \frac{\alpha^{m+1}-\beta^{m+1}}{\alpha-\beta}\cdot \frac{\gamma^{m+1}-\delta^{m+1}}{\gamma-\delta}\\&=\sum_{m\geq 0} X^m \cdot\frac{\alpha^{m+1}-\beta^{m+1}}{\alpha-\beta}\cdot \frac{\gamma^{m+1}-\delta^{m+1}}{\gamma-\delta},
\end{align*}
where $X=\ell^{-s}\eta(\ell)$. We can manipulate the latter series and obtain
\begin{align*}
\sum_{m\geq 0} X^m \cdot\frac{\alpha^{m+1}-\beta^{m+1}}{\alpha-\beta}\cdot & \frac{\gamma^{m+1}-\delta^{m+1}}{\gamma-\delta}=\frac{1}{(\alpha-\beta)(\gamma-\delta)}\left(\frac{\alpha\gamma}{1-\alpha\gamma X}-\frac{\alpha\delta}{1-\alpha\delta X}-\frac{\beta\gamma}{1-\beta\gamma X}+\frac{\beta\delta}{1-\beta\delta X}\right)\\&=\frac{1}{\alpha-\beta}\left( \frac{\alpha+\alpha\beta^2\gamma\delta X^2-\beta-\alpha^2\beta\gamma\delta X^2}{(1-\alpha\gamma X)(1-\alpha\delta X)(1-\beta\gamma X)(1-\beta\delta X)}\right)\\&=\frac{1-\alpha\beta\gamma\delta X^2}{(1-\alpha\gamma X)(1-\alpha\delta X)(1-\beta\gamma X)(1-\beta\delta X)}.
\end{align*}
This is a standard computation, see for example Jacquet's refreshing exercise \cite[Lemma 15.9.4]{Jacquet}. We have conditions on the convergence giving the condition $\re(s)>r(\sigma,\eta)$. 
Substituting $X=\ell^{-s}\eta(\ell)$, we find 
\begin{align*}
\int_{\Q_{\ell}^{\times}}|y|^{s-1}\eta(y)\W_{\varphi_0}\left( \SmallMatrix{y&0\\0&1},\SmallMatrix{y&0\\0&1}\right)d^{\times}y&=(1-\chi_{\sigma}(\ell)\eta^2(\ell)\ell^{-2s})L(\sigma\otimes\eta,s)=L(\eta^2\chi_{\sigma},2s)^{-1}L(\sigma\otimes\eta,s).
\end{align*}
\end{proof}
\end{lemma}

\begin{lemma}[Action of $U(\ell)$ on the zeta integral]\label{zetaU(l)}
If $\varphi_0$ is the normalised spherical vector as above, we have
\begin{displaymath}
Z(\sigma,\eta,U(\ell)\varphi_0,s)=\tfrac{\ell^{s+1}}{\eta(\ell)}[Z(\sigma,\eta,\varphi_0,s)-L(\sigma\otimes\eta,s)^{-1}]=\tfrac{\ell^{s+1}}{\eta(\ell)}[L(\eta^2\chi_{\sigma},2s)^{-1}-L(\sigma\otimes\eta,s)^{-1}]
\end{displaymath}
\begin{proof}
First we apply Proposition \ref{actionUellWhitt} and find
\begin{align*}
Z(\sigma,\eta,U(\ell)\varphi_0,s)&=\ell^{2}L(\sigma\otimes\eta,s)^{-1}\int_{|y|< \ell}|y|^{s-1}\eta(y)\W_0\left(\SmallMatrix{\ell y & 0 \\ 0 &1},\SmallMatrix{\ell y & 0 \\ 0 &1}\right)d^{\times}y\\&=\ell^{s+1}\eta(\ell)^{-1}L(\sigma\otimes\eta,s)^{-1}\int_{|y|< 1}|y|^{s-1}\eta(y)\W_0\left(\SmallMatrix{ y & 0 \\ 0 &1},\SmallMatrix{y & 0 \\ 0 &1}\right)d^{\times}y,
\end{align*}
where in the second equality we used the change of variables $y \rightsquigarrow \ell y$. We then rewrite the integral in the last term as
\begin{align*}
\int_{\Q_{\ell}^{\times}}|y|^{s-1}\eta(y)\W_0\left(\SmallMatrix{ y & 0 \\ 0 &1},\SmallMatrix{y & 0 \\ 0 &1}\right)d^{\times}y-\int_{|y|\geq 1}|y|^{s-1}\eta(y)\W_0\left(\SmallMatrix{ y & 0 \\ 0 &1},\SmallMatrix{y & 0 \\ 0 &1}\right)d^{\times}y.
\end{align*}
Then we apply Theorem \ref{whitvaluesG} and obtain
\begin{align*}
\int_{|y|\geq 1}|y|^{s-1}\eta(y)\W_0\left(\SmallMatrix{ y & 0 \\ 0 &1},\SmallMatrix{y & 0 \\ 0 &1}\right)d^{\times}y=\int_{\Z_{\ell}^{\times}}\W_0\left(\SmallMatrix{ y & 0 \\ 0 &1},\SmallMatrix{y & 0 \\ 0 &1}\right)d^{\times}y=\int _{\Z_{\ell}^{\times}} d^{\times}y = 1.
\end{align*}
Putting everything together we find 
\begin{align*}
Z(\sigma,\eta,U(\ell)\varphi_0,s)&=\ell^{s+1}\eta(\ell)^{-1}L(\sigma\otimes\eta,s)^{-1}\left( \int_{\Q_{\ell}^{\times}}|y|^{s-1}\eta(y)\W_0\left(\SmallMatrix{ y & 0 \\ 0 &1},\SmallMatrix{y & 0 \\ 0 &1}\right)d^{\times}y- 1\right)\\&=\ell^{s+1}\eta(\ell)^{-1}(Z(\sigma,\eta,\varphi_0,s)-L(\sigma\otimes\eta,s)^{-1}).
\end{align*}
\end{proof}
\end{lemma}

\begin{lemma}[Action of the Borel subgroup of $\GL_2(\Q_{\ell})$]\label{borelH}
For any $f\in\sigma$, $a,d\in\Q_{\ell}^{\times}$, we have
\begin{displaymath}
Z(\sigma,\eta,\left(\SmallMatrix{a&*\\0&d},\SmallMatrix{a&*\\0&d}\right)\cdot f,s)=\left|\tfrac{d}{a}\right|^{s-1}\chi_{\sigma}(d)\eta(a^{-1}d) \cdot Z(\sigma,\eta,f,s)
\end{displaymath}
\begin{proof}
We apply Lemma \ref{lemmaactionborelH} and find
\begin{align*}
Z(\sigma,\eta,\left(\SmallMatrix{a&*\\0&d},\SmallMatrix{a&*\\0&d}\right)\cdot f,s)&=\chi_{\sigma}(d)L(\sigma\otimes\eta,s)^{-1}\int_{\Q_{\ell}^{\times}}|y|^{s-1}\eta(y)\W_f\left(\SmallMatrix{d^{-1}ay & 0 \\ 0 &1},\SmallMatrix{d^{-1}ay & 0 \\ 0 &1}\right)d^{\times}y \\&=\chi_{\sigma}(d)|d/a|^{s-1}\eta(a^{-1}d)L(\sigma\otimes\eta,s)^{-1}\int_{\Q_{\ell}^{\times}}|y|^{s-1}\eta(y)\W_f\left(\SmallMatrix{y & 0 \\ 0 &1},\SmallMatrix{y & 0 \\ 0 &1}\right)d^{\times}y\\&=\chi_{\sigma}(d)|d/a|^{s-1}\eta(a^{-1}d)\cdot Z(\sigma,\eta,f,s)
\end{align*}
where in the second equality we used the change of variable $y \rightsquigarrow d^{-1}ay$.
\end{proof}
\end{lemma}

\section{Towards norm relations}\label{towards}
Let $G$ be the algebraic group over $\Q$ defined in the introduction, i.e. $G=\res^{\Q}_F\GL_2$,  for $F$ real quadratic field. We will now prove some results using the zeta integrals of the two previous sections. We will denote with $\sigma$ an unramified irreducible principal series representation of $G(\Q_{\ell})$, i.e. $\sigma= I_{\GL_2(\Q_{\ell})\times\GL_2(\Q_{\ell})}(\underline{\chi},\underline{\psi})$ if $\ell$ splits and $\sigma=I_{\GL_2(F_{\ell})}(\tilde{\chi},\tilde{\psi})$ for $F_{\ell}$ the unramified quadratic extension of $\Q_{\ell}$ if $\ell$ is inert. We will denote with $\chi_{\sigma}$ the central character of $\sigma$, namely $\chi_{\sigma}=\chi_1\psi_1\chi_2\psi_2$ in the first case and $\chi_{\sigma}=\tilde{\chi}\tilde{\psi}$ in the second one. By abuse of notation, we will often write $H$ for $H(\Q_{\ell})=\GL_2(\Q_{\ell})$ and denote with $L(\as(\sigma),s)$ both the local $L$-factor we considered in $\S$\ref{lfnctinert} and $\S$\ref{lfnctsplit}, i.e.
\begin{displaymath}
L(\as(\sigma),s)=\begin{cases}
L(\sigma,s) &\text{if $\ell$ splits and $\sigma=I_{\GL_2(\Q_{\ell})\times\GL_2(\Q_{\ell})}(\underline{\chi},\underline{\psi})$}\\
L(\as(\sigma),s)  &\text{if $\ell$ is inert and $\sigma=I_{\GL_2(F)}(\tilde{\chi},\tilde{\psi})$}.
\end{cases}
\end{displaymath}
We also let 
\begin{displaymath}
\alpha_i=\chi_i(\ell), \beta_i=\psi_i(\ell) \ \ \ \text{if $\ell$ splits,}
\end{displaymath}
\begin{displaymath}
\alpha=\tilde{\chi}(\ell), \beta=\tilde{\psi}(\ell) \ \ \ \text{if $\ell$ splits,}
\end{displaymath}
\subsection{Multiplicity one} We will fix $\sigma$ as above and another pair of unramified characters $\chi,\psi$ satisfying
\begin{displaymath}
\chi\psi\cdot\chi_{\sigma}=1.
\end{displaymath}
We will moreover assume that $I_H(\chi,\psi)$ is either irreducible or it has an infinite dimensional subrepresentation, i.e. $\chi\psi^{-1}\neq |\cdot|^{- 1}$. 
We will be considering the embedding
\begin{displaymath}
\iota: H(\Q_{\ell}) \hookrightarrow G(\Q_{\ell}).
\end{displaymath}
In the split case $\iota(h):=(h,h)\in \GL_2(\Q_{\ell})\times \GL_2(\Q_{\ell})$, while in the inert case $\iota(h):=h \in \GL_2(\Q_{\ell})\subset \GL_2(F)$.

\begin{thm}[Multiplicity one]\label{thmmult1}
Let $\sigma,\chi,\psi$ as above. We assume that
\begin{equation}\label{inertmult1}
|\alpha_i|_{\C}=|\beta_i|_{\C} \ \ \text{if $\ell$ splits, }\ \ \ |\alpha|_{\C}=|\beta|_{\C} \ \ \text{if $\ell$ is inert}\tag{$\star$}
\end{equation}
and that $|\chi(\ell)|_{\C}\neq|\psi(\ell)|_{\C}$. We then have 
\begin{displaymath}
\dim\left( \Hom_H(I_H(\chi,\psi)\otimes\sigma, \C) \right)\leq 1.
\end{displaymath}
\begin{proof}
We treat separately the split and inert case. The main ingredient is the exact sequence $(\ref{eq:exseq})$, which will allow us to describe $\sigma_{|H}$ in terms of simpler representations of $H$.

If $\ell$ splits and $I_H(\chi,\psi)$ is irreducible, i.e. $\chi\psi^{-1}\neq |\cdot|$,  this is Theorem 1.1 of \cite{prasad}. We apply it for $V_1=I_H(\chi,\psi),V_2=I_H(\chi_1,\psi_1),V_3=I_H(\chi_2,\psi_2)$ (i.e. $V_2\otimes V_3 = \sigma$). To deal with the case when $\chi\psi^{-1}= |\cdot|$ we will make use of the exact sequence $(\ref{eq:exseq})$, following the strategy of \cite[Proof of Theorem 1.2 Case 2]{prasad}. Let $G=\GL_2(\Q_{\ell})\times \GL_2(\Q_{\ell})$, $J=B(\Q_{\ell})\times B(\Q_{\ell})$, $H=\GL_2(\Q_{\ell})$ and $\tau$ is given by
\begin{displaymath}
\tau \left( \SmallMatrix{a&b\\0&d}, \SmallMatrix{a'&b'\\0&d'} \right)=\chi_1(a)\psi_1(d)\chi_2(a')\psi_2(d').
\end{displaymath}
Recall that we denoted with $H_1=\stab_H(1\cdot J)$ and $H_2=\stab_H(\epsilon\cdot J)$, where $\epsilon\in G$ such that $\epsilon\cdot J$ is in the open orbit. In this case $H_1=B(\Q_{\ell})$ and $H_2=T=\{\SmallMatrix{\lambda_1&0\\0&\lambda_2}, \lambda_i\in \Q_{\ell}^{\times}\}$ (the maximal split torus). To prove the last equality, one can take $\epsilon= (\id, \SmallMatrix{0&1\\-1&0})$. Using the fact that $T$ is unimodular and 
\begin{displaymath}
\delta_{B\times B}\left( \SmallMatrix{a&b\\0&d}, \SmallMatrix{a'&b'\\0&d'} \right)=|\tfrac{a}{d}||\tfrac{a'}{d'}|,
\end{displaymath}
\begin{displaymath}
\delta_{B}\left( \SmallMatrix{a&b\\0&d}\right)=|\tfrac{a}{d}|,
\end{displaymath}
\begin{displaymath}
\delta_{B\times B}\left( \epsilon (\SmallMatrix{\lambda_1&0\\0&\lambda_2},\SmallMatrix{\lambda_1&0\\0&\lambda_2}) \epsilon ^{-1}\right)=1, 
\end{displaymath}
we find the exact sequence of $\GL_2(\Q_{\ell})$-modules 
\begin{displaymath}
0\to \text{c-}\ind_{H_2}^{\GL_2(\Q_{\ell})}\tilde{\tau} \to \sigma_{|H} \to I_H(\chi_1\chi_2 |\cdot|^{1/2},\psi_1\psi_2|\cdot|^{-1/2}) \to 0,
\end{displaymath}
where $\tilde{\tau}(\SmallMatrix{\lambda_1&0\\0&\lambda_2})=\chi_1\psi_2(\lambda_1)\psi_1\chi_2(\lambda_2)$. Let now $V=I_H(\chi,\psi)$. Applying $\Hom_H(-,V^{\vee})$ to the above exact sequence we find
\begin{align*}
0\to &\Hom_H(I_H(\chi_1\chi_2 |\cdot|^{1/2},\psi_1\psi_2|\cdot|^{-1/2}),V^{\vee})\to \Hom_H(\sigma_{|H},V^{\vee}) \\&\to \Hom_H(\text{c-}\ind_{H_2}^{\GL_2(\Q_{\ell})}\tilde{\tau},V^{\vee}) \to \Ext^1_H(I_H(\chi_1\chi_2 |\cdot|^{1/2},\psi_1\psi_2|\cdot|^{-1/2}),V^{\vee}))\to \dots
\end{align*}
Since the smooth dual of $V$ is $V^{\vee}=I_H(\chi^{-1},\psi^{-1})$, we have, from the second part of Theorem \ref{thmmult}, that 
\begin{displaymath}
\Hom_H(I_H(\chi_1\chi_2 |\cdot|^{1/2},\psi_1\psi_2|\cdot|^{-1/2}),V^{\vee})\neq 0,
\end{displaymath}
if and only if $\chi_1\chi_2 |\cdot|^{1/2}=\chi^{-1},\psi_1\psi_2|\cdot|^{-1/2}=\psi^{-1}$ or
$\chi_1\chi_2 |\cdot|^{1/2}=\psi^{-1},\psi_1\psi_2|\cdot|^{-1/2}$. Since $\chi\psi^{-1}=|\cdot|$ we can write $\chi=|\cdot|^{1/2}\gamma^{-1},\psi=|\cdot|^{-1/2}\gamma^{-1}$, for some quasicharacter $\gamma$. Hence the above space is non-zero if and only if the pair $(\chi_1\chi_2,\psi_1\psi_2)$ is equal to
\begin{displaymath}
(|\cdot|^{-1}\gamma,|\cdot|\gamma) \ \ \ \text{or} \ \ \ (\gamma,\gamma) .
\end{displaymath}
The first case can not happen because of the assumption (\ref{inertmult1}). In the second case we are reduced to compute the dimension of $\Hom_H(\sigma \otimes \tau, \C)$, where 
\begin{displaymath}
\sigma = I_H(\chi,\chi)\otimes I_H(\gamma\chi^{-1},\gamma\chi^{-1}), \ \ \ \tau = I_H(|\cdot|^{1/2}\gamma^{-1},|\cdot|^{-1/2}\gamma^{-1}).
\end{displaymath}
The argument of \cite[Lemma 2.8]{LLind} works for showing that the number of irreducible components of the restriction to $G^*(\Q_{\ell})$ of the $G(\Q_{\ell})$-representation $\sigma$ is equal to the number of characters $\omega$ of $G(\Q_{\ell})/G^*(\Q_{\ell})\simeq \Q_{\ell}^{\times}/(\Q_{\ell}^{\times})^2$ such that $\sigma\otimes \omega \simeq \sigma$. In this case we find, again applying Theorem \ref{thmmult}, that the only such character is $\omega=\mathbf{1}$ and hence $\sigma$ is irreducible as $G^*(\Q_{\ell})$ representation. Without loss of generality we can assume that $\gamma$ is trivial and hence that both $\sigma$ and $\tau$ are trivial on the centre. Hence they factor through $G^*(\Q_{\ell})/\Q_{\ell}^{\times}\simeq$SO$(4)$ and PGL$(2)\simeq $SO$(3)$ respectively. The main result of \cite{shintani} shows that for any representations $\sigma$ of SO$(n+1)$ and $\tau$ of SO$(n)$ generated by a spherical vector, the dimension of $\Hom_{SO(n)}(\sigma\otimes \tau,\C)$ is at most 1. 
Finally, if $(\chi_1\chi_2,\psi_1\psi_2)\neq (\gamma,\gamma)$ the first space in the sequence is zero which implies, by \cite[Corollary 5.9]{prasad}, that also the $\Ext^1$ is zero. We hence find
\begin{displaymath}
\Hom_H(\sigma_{|H},V^{\vee})\simeq \Hom_H(\text{c-}\ind_{H_2}^{GL_2(\Q_{\ell})}\tilde{\tau},V^{\vee})\simeq \Hom_{H_2}(\tilde{\tau},(V_{|H_2})^{\vee}).
\end{displaymath}
By assumption we have $\chi^{-1}\psi^{-1}=\chi_1\psi_1\chi_2\psi_2$, hence the central characters of $\tilde{\tau}$ and $V^{\vee}$ agree. We can apply \cite[Lemma 9]{wald}, saying that this space is at most one dimensional.

We now prove the inert case, applying again the exact sequence $(\ref{eq:exseq})$. Let $G=GL_2(F_{\ell})$, where $F_{\ell}$ is the quadratic unramified extension of $\Q_{\ell}$ and choose a $\Q_{\ell}$ basis $\{1,\alpha\}$. Then let $J=B(F_{\ell})$, $H=\GL_2(\Q_{\ell})$ and $\tau$ be the smooth representation of $J$ given by
\begin{displaymath}
\tau \left( \SmallMatrix{a&b\\0&d}\right)=\tilde{\chi}(a)\tilde{\psi}(d).
\end{displaymath}
The two orbits of the action of $\GL_2(\Q_{\ell})$ on $\GL_2(F_{\ell})/B(F_{\ell})\simeq \mathbb{P}^1_{F_{\ell}}$ are the $\GL_2(\Q_{\ell})$ orbit of $(1:0)$ (essentially $\mathbb{P}^1_{\Q_{\ell}}$) and the $\GL_2(\Q_{\ell})$ orbit of $(1:\alpha)$, which is given by $(a:b)\in \mathbb{P}^1_{F_{\ell}}$ such that $a/b\not\in \Q_{\ell}, ab\neq 0$. Writing $a=x_0+y_0\alpha,b=x_1+y_1\alpha$, we have
\begin{displaymath}
\SmallMatrix{y_1&-y_0\\-x_1&x_0}\cdot (a:b)=(x_0y_1-x_1y_0:\alpha(x_0y_1-x_1y_0))=(1:\alpha)
\end{displaymath}
and $x_0y_1-x_1y_0\neq 0$, otherwise if $y_1\neq 0$, $a/b=y_0/y_1\in \Q_{\ell}$. If $y_1=0$, then $x_1\neq 0$ and $y_0\neq 0$. We find that the stabiliser of the closed orbit is $H_1=B(\Q_{\ell})$ and, taking $\epsilon=\SmallMatrix{1&0\\ \alpha&1}$, the stabiliser of the open one is $H_2=\{\SmallMatrix{a&b\\ \alpha^2b&a}, a,b \in\Q_{\ell}^2-(0,0)\}$. To see that $H_2$ can be written of this form, we compute the conjugate of $B(F_{\ell})$ by $\epsilon$.
\begin{displaymath}
\SmallMatrix{1&0\\ \alpha&1}\SmallMatrix{a&b\\0&d}\SmallMatrix{1&0\\ -\alpha&1}= \SmallMatrix{a-b\alpha &b\\a\alpha-b\alpha^2-d\alpha &b\alpha-d}.
\end{displaymath}
Requiring that such matrices lie in $\GL_2(\Q_{\ell})$ implies that $b\in\Q_{\ell}$ and $a=a_1+bt,d=a_1-bt$, for $a_1\in\Q_{\ell}$. The group $H_2$ is a (non-split) maximal torus in $\GL_2(\Q_{\ell})$. This is again unimodular and
\begin{displaymath}
\delta_{B(F_{\ell})}(\epsilon\SmallMatrix{a&b\\ \alpha^2b&a}\epsilon^{-1})=\delta_{B(F_{\ell})}(\SmallMatrix{a+b\alpha &b\\0&a-b\alpha})=|a+b\alpha||a-b\alpha|^{-1}=1.
\end{displaymath} We find the exact sequence of $\GL_2(\Q_{\ell})$-modules 
\begin{displaymath}
0\to \text{c-}\ind_{H_2}^{GL_2(\Q_{\ell})}\tilde{\tau} \to \sigma_{|H} \to I_H(\tilde{\chi}_{|\Q_{\ell}^{\times}},\tilde{\psi}_{|\Q_{\ell}^{\times}}) \to 0,
\end{displaymath}
where $\tilde{\tau}(\SmallMatrix{a&b\\ \alpha^2b&a})=\tilde{\chi}(a+bt)\cdot \tilde{\psi}(a-bt)$.

Applying $\Hom_H(-,V^{\vee})$ as before to the above exact sequence we find
\begin{align*}
0\to &\Hom_H(I_H(\tilde{\chi}_{|\Q_{\ell}^{\times}},\tilde{\psi}_{|\Q_{\ell}^{\times}}),V^{\vee})\to \Hom_H(\sigma_{|H},V^{\vee}) \\&\to \Hom_H(\text{c-}\ind_{H_2}^{GL_2(\Q_{\ell})}\tilde{\tau},V^{\vee}) \to \Ext^1_H(I_H(\tilde{\chi}_{|\Q_{\ell}^{\times}},\tilde{\psi}_{|\Q_{\ell}^{\times}}),V^{\vee}))\to \dots
\end{align*}
Since the smooth dual of $V$ is $V^{\vee}=I_H(\chi^{-1},\psi^{-1})$, we have, arguing as above, that 
\begin{displaymath}
\Hom_H(I_H(\tilde{\chi}_{|\Q_{\ell}^{\times}},\tilde{\psi}_{|\Q_{\ell}^{\times}}),V^{\vee})\neq 0,
\end{displaymath}
if and only if $\tilde{\chi}_{|\Q_{\ell}^{\times}}=\chi^{-1},\tilde{\psi}_{|\Q_{\ell}^{\times}}=\psi^{-1}$ or
$\tilde{\chi}_{|\Q_{\ell}^{\times}}=\psi^{-1},\tilde{\psi}_{|\Q_{\ell}^{\times}}=\chi^{-1}$. Assumption (\ref{inertmult1}) implies that this is not the case, hence the above space is zero and so is $\Ext^1$. We hence find
\begin{displaymath}
\Hom_H(\sigma_{|H},V^{\vee})\simeq \Hom_H(\text{c-}\ind_{H_2}^{GL_2(\Q_{\ell})}\tilde{\tau},V^{\vee})\simeq \Hom_{H_2}(\tilde{\tau},(V_{|H_2})^{\vee}).
\end{displaymath}
By assumption we have $\chi^{-1}\psi^{-1}=(\tilde{\chi}\cdot \tilde{\psi})_{|\Q_{\ell}^{\times}}$, and we can again apply \cite[Lemma 9]{wald}, saying that this space is at most one dimensional. 
\end{proof}
\end{thm}

\subsection{A basis for \texorpdfstring{$\Hom_H(I_H(\chi,\psi)\otimes\sigma, \C)$}{Hom_GL2}} Using the zeta integral defined above, we now want to construct an explicit nonzero element of $\Hom_H(I_H(\chi,\psi)\otimes\sigma, \C)$, which by the above theorem will be a basis. 

\begin{defi}\label{defizetino}
Let $\eta=\psi$, for $\psi$ as above. For any $\varphi\in\sigma,s\in\C$, we define a function $z_{s,\varphi}$ on $H(\Q_{\ell})$ by
\begin{displaymath}
z_{s,\varphi}(h):=Z(\sigma, \eta, \iota(h)\cdot \varphi,s+\tfrac{1}{2}),
\end{displaymath}
for any $h\in H(\Q_{\ell})$.
\end{defi}

We now let, for $s\in \C$, $\psi_s:=\psi|\cdot|^{-s}, \chi_s:=\chi|\cdot|^{s}$.

\begin{prop}\label{zetino}
The above function defines an element $z_{\varphi}\in I_H(\psi_s^{-1}, \chi_s^{-1})$ for every $\varphi\in\sigma$. Moreover
\begin{displaymath}
z_{s,\varphi_0}(1)=L(\tfrac{\psi}{\chi}, 2s+1)^{-1},
\end{displaymath}
\begin{displaymath}
z_{s,U(\ell)\varphi_0}(1)=\tfrac{\ell^{s+3/2}}{\eta(\ell)}[L(\tfrac{\psi}{\chi}, 2s+1)^{-1}-L(\as(\sigma\otimes\eta),s+\tfrac{1}{2})^{-1}]
\end{displaymath}
\begin{proof}
The first assertion is a straightforward corollary of Lemma \ref{borelHinert} and Lemma \ref{borelH}. Indeed 
\begin{align*}
z_{\varphi}(\SmallMatrix{a&*\\0&d}\cdot h)&=Z(\sigma,\eta,\iota(\SmallMatrix{a&*\\0&d})\iota(h)\cdot {\varphi},\tfrac{1}{2})=\left|\tfrac{a}{d}\right|^{s+1/2}\chi_{\sigma}(d)\eta(a^{-1}d) \cdot Z(\sigma,\eta,\iota(h)\cdot {\varphi},\tfrac{1}{2})\\&=\left| \tfrac{a}{d}\right|^{1/2} \psi^{-1}(a)|a|^s\chi^{-1}(d)|d|^{-s}\cdot z_{\varphi}(h),
\end{align*}
using $\chi_{\sigma}=\psi^{-1}\chi^{-1}$ and $\eta=\psi$. The formula for the value at $\varphi_0$ ($U(\ell)\varphi_0$ respectively) follows from Lemma \ref{sphericalzetainert} and Lemma \ref{sphericalzeta} (Lemma \ref{zetaU(l)inert} and Lemma \ref{zetaU(l)} respectively).
\end{proof}
\end{prop}

By definition, the map 
\begin{align*}
z_s: \ &\sigma \rightarrow I_H(\psi_s^{-1}, \chi_s^{-1}) \\& {\varphi}\mapsto z_{s,\varphi}
\end{align*}
is $H$-equivariant. Moreover it follows from the proposition that $z_0$ is different from zero if $L(\tfrac{\psi}{\chi}, 2s+1)^{-1}$ and $L(\as(\sigma\otimes\eta),s+\tfrac{1}{2})^{-1}$ do not both vanish at $s=0$. Notice that if $\ell$ is inert then $L(\tfrac{\psi}{\chi}, 2s+1)^{-1}$ divides $L(\as(\sigma\otimes\eta),s+\tfrac{1}{2})^{-1}$.

\begin{lemma}\label{vanish1}
If $\ell$ splits, assume that $L(\tfrac{\psi}{\chi}, 2s+1)^{-1}$ and $L(\as(\sigma\otimes\eta),s+\tfrac{1}{2})^{-1}$ do not both vanish at $s=0$. If $\ell$ is inert, assume that $L(\tfrac{\psi}{\chi}, 2s+1)^{-1}$ and $L(\tfrac{\psi}{\chi}, 2s+1)L(\as(\sigma\otimes\eta),s+\tfrac{1}{2})^{-1}$ do not both vanish at $s=0$. Then the image of the homomorphism $z_0$ is contained in the unique irreducible subrepresentation of $I_H(\psi^{-1},\chi^{-1})$.
\begin{proof}
If $L(\tfrac{\psi}{\chi}, 2s+1)^{-1}$ does not vanish, then $I_H(\psi^{-1},\chi^{-1})$ is irreducible and there is nothing to prove. Otherwise, $\chi\psi^{-1}=|\cdot|$ and $I_H(\psi^{-1},\chi^{-1})$ has a unique infinite dimensional irreducible subrespresentation St$(\gamma)$ and one dimensional quotient with action given by $\gamma(\det)$. We claim that if $L(\as(\sigma\otimes\eta),s+\tfrac{1}{2})^{-1}$ does not vanish at $s=0$, the space $\Hom_H(\sigma_{|H}, \gamma(\det))$ is zero, and, consequently, the image of $z_0$ is contained in St$(\gamma)$. The proof uses the same methods as the one of Theorem \ref{thmmult1}. With the same notation, in the split case one finds an exact sequence 
\begin{align*}
0\to &\Hom_H(I_H(\chi_1\chi_2 |\cdot|^{1/2},\psi_1\psi_2|\cdot|^{-1/2}),\gamma(\det))\to \Hom_H(\sigma_{|H},\gamma(\det)) \\&\to \Hom_H(\text{c-}\ind_{H_2}^{\GL_2(\Q_{\ell})}\tilde{\tau},\gamma(\det)) \to \Ext^1_H(I_H(\chi_1\chi_2 |\cdot|^{1/2},\psi_1\psi_2|\cdot|^{-1/2}),\gamma(\det))\to \dots
\end{align*}
Since $L(\as(\sigma\otimes\psi),\tfrac{1}{2})^{-1}$ is not zero, then none of the characters $\chi_1\psi_2,\chi_1\chi_2,\psi_1\chi_2, \psi_1\psi_2$ is equal to $\psi^{-1}|\cdot|^{-1/2}=\gamma$. On the other hand, applying Frobenius reciprocity one finds that the first space in the sequence is non zero if and only if $\chi_1\chi_2=\gamma$ and $\psi_1\psi_2=\gamma$, while the third one is non zero if and only if $\chi_1\psi_2=\gamma$ and $\psi_1\chi_2=\gamma$. This proves the claim.  

Similarly in the inert case, we have
\begin{align*}
0\to &\Hom_H(I_H(\tilde{\chi}_{|\Q_{\ell}^{\times}},\tilde{\psi}_{|\Q_{\ell}^{\times}}),\gamma(\det))\to \Hom_H(\sigma_{|H},\gamma(\det)) \\&\to \Hom_H(\text{c-}\ind_{H_2}^{GL_2(\Q_{\ell})}\tilde{\tau},\gamma(\det)) \to \Ext^1_H(I_H(\tilde{\chi}_{|\Q_{\ell}^{\times}},\tilde{\psi}_{|\Q_{\ell}^{\times}}),\gamma(\det))\to \dots
\end{align*} 
The assumption $L(\tfrac{\psi}{\chi}, 2s+1)L(\as(\sigma\otimes\eta),s+\tfrac{1}{2})^{-1}=(1-\tilde{\chi}\psi(\ell)\ell^{-s-1/2})(1-\tilde{\psi}\psi(\ell)\ell^{-s-1/2})\neq 0$ at $s=0$ and Frobenius reciprocity again imply that the first and third space in the sequence are zero. 
\end{proof}
\end{lemma}

Unlike in the GSp$_4$ case where temperedness considerations allow to assume the non-vanishing of both the abelian L-factor and the one where the principal series appears, in our setting it can actually happen that they both are zero, e.g. it is possible that at some split primes $\sigma= I_H(\chi_1,\psi_1)\otimes I_H(\gamma\chi_1^{-1},\gamma \psi_1^{-1})$. The following lemma shows that in this case $z_0$ is identically zero.

\begin{lemma}\label{vanish2}
If the assumptions of the previous lemma are not satisfied, then $z_0$ is identically zero.
\begin{proof}
This follows from the explicit description of functions in the Krillov model of principal series representations of $\GL_2(\Q_{\ell})$ as recalled for example in \cite[Lemma 14.3]{Jacquet}. In the split case, the functions $y\mapsto W_{i, \varphi}(\SmallMatrix{y&0\\0&1})$ are indeed in the Krillov model of $I_H(\chi_i,\psi_i)$. By definition 
\begin{displaymath}
L(\as(\sigma\otimes\psi),s+1/2)=L(\chi_1\chi_2\gamma^{-1},s)L(\chi_1\psi_2\gamma^{-1},s)L(\psi_1\chi_2\gamma^{-1},s)L(\psi_1\psi_2\gamma^{-1},s)
\end{displaymath}
Firstly we assume that the order of vanishing of $L(\as(\sigma\otimes\psi),s+1/2)^{-1}$ is 2 and, without loss of generality, we can assume $\chi_1\chi_2=\psi_1\psi_2=\gamma$, where $\psi^{-1}=|\cdot|^{1/2}\gamma$. Since the order is 2, we have $\chi_i\neq\psi_i$, and $W_{i, \varphi}(\SmallMatrix{y&0\\0&1})$ can be written 
\begin{displaymath}
f_i(y)\chi_i(y)|y|^{1/2}+g_i(y)\psi_i(y)|y|^{1/2}
\end{displaymath}
for some $f_i,g_i\in \Ss(\Q_{\ell})$. Hence the function $Z(\sigma,\psi, \varphi_1\otimes\varphi_2, s+1/2)$ is equal to $L(\as(\sigma\otimes\psi),s+1/2)^{-1}$ multiplied by the integral
\begin{align*}
&\int_{\Q_{\ell}^{\times}}|y|^{s-1}\gamma^{-1}(y)(f_1(y)\chi_1(y)|y|^{1/2}+g_1(y)\psi_1(y)|y|^{1/2})(f_2(y)\chi_2(y)|y|^{1/2}+g_2(y)\psi_2(y)|y|^{1/2})d^{\times}y \\&=P_1(s)L(\mathbf{1},s)+P_2(s)L(\chi_1\psi_2\gamma^{-1},s)+P_3(s)L(\psi_1\chi_2\gamma^{-1},s),
\end{align*}
where $P_i(s)$ are polynomials in $\ell^{-s},\ell^s$. The equality follows from the description of the $L$-factor $L(\mu,s)$ for any quasicharacter of $\Q_{\ell}^{\times}$ (see \cite{Jacquet}, the discussion after Lemma 14.3). Since in our situation we have 
\begin{displaymath}
L(\as(\sigma\otimes\psi),s+1/2)^{-1}=L(\mathbf{1},s)^{-2}L(\chi_1\psi_2\gamma^{-1},s)^{-1}L(\psi_1\chi_2\gamma^{-1},s)^{-1},
\end{displaymath}
the result follows. If the order of vanishing of $L(\as(\sigma\otimes\psi),s+1/2)^{-1}$ is 4, we have $\chi_1=\psi_1$ and $\chi_2=\psi_2=\gamma\chi_1^{-1}$. In this case 
\begin{align*}
&\int_{\Q_{\ell}^{\times}}|y|^{s-1}\gamma^{-1}(y)(f_1(y)\chi_1(y)|y|^{1/2}+g_1(y)\chi_1(y)v(y)|y|^{1/2})(f_2(y)\gamma\chi_1^{-1}(y)|y|^{1/2}+g_2(y)\gamma\chi_1^{-1}(y)v(y)|y|^{1/2})d^{\times}y\\&=P_1(s)L(\mathbf{1},s)+P_2(s)L(\mathbf{1},s)^2+P_3(s)L(\mathbf{1},s)^3.
\end{align*}
Here $v$ is the valuation of $\Q_{\ell}$ and the equality follows from what said above and \cite[(14.2.1)]{Jacquet} in the case where $v$ or $v^2$ appears in the integral. Being $L(\as(\sigma\otimes\psi),s+1/2)^{-1}$ equal to $L(\mathbf{1},s)^{-4}$, again the result follows.

Similarly, in the inert case we have that the Krillov function $y\mapsto W_{\varphi}(\SmallMatrix{y&0\\0&1})$ can be written as
\begin{displaymath}
|y|_{F_{\ell}}^{1/2}(f(y)\gamma(y)+g(y)\gamma(y)v(y)).
\end{displaymath} 
for $f,g\in\Ss(F_{\ell})$ and $|\cdot|_{F_{\ell}}$ is equal to $|\cdot|^2$ when restricted to $\Q_{\ell}$. Indeed the fact that both the $L$ factor vanish implies that $\tilde{\chi}=\tilde{\psi}=\gamma$, where $\psi^{-1}=|\cdot|^{1/2}\gamma$. Hence the integral in the definition of $z_{s,\varphi}$ is
\begin{align*}
&\int_{\Q_{\ell}^{\times}}|y|^{s}\gamma^{-1}(y)(f(y)\gamma(y)+g(y)\gamma(y)v(y))d^{\times}y=P_1(s)L(\mathbf{1},s)+P_2(s)L(\mathbf{1},s)^2.
\end{align*}
In this case $L(\as(\sigma\otimes \psi),s+1/2)=L(\mathbf{1},s)^{-4}$ and the result follows.
\end{proof}
\end{lemma}

We recall then the intertwining operator defined thanks to Proposition \ref{interwop} and the pairing of definition \ref{pairing}
\begin{displaymath}
M: I_H(\chi_s,\psi_s)\rightarrow I_H(\psi_s,\chi_s) \ \ \text{ and }\ \ \langle-,-\rangle: I_H(\psi_s,\chi_s)\times I_H(\psi_s^{-1},\chi_s^{-1})\rightarrow \C.
\end{displaymath}
\begin{defiprop}\label{defibasis} For every $f\in I_H(\chi,\psi), {\varphi}\in\sigma$, we let 
\begin{displaymath}
\zita_{\chi,\psi}(f\otimes {\varphi}):=\lim_{s\to 0}L(\tfrac{\psi}{\chi},2s+1)\langle Mf_s, z_{s,{\varphi}} \rangle,
\end{displaymath}
where $f_s\in I_H(\chi_s,\psi_s)$ is any polynomial section passing through $f$. 
This gives a well defined element $\zita_{\chi,\psi}\in \Hom_H(I_H(\chi,\psi)\otimes\sigma, \C)$ which is not zero. 
\begin{proof}
First of all one notices that $\langle Mf, z_{s,{\varphi}} \rangle = \langle f_s, Mz_{s,{\varphi}}\rangle$. If $\chi,\psi,\sigma$ are as in Lemma \ref{vanish2}, then $z_s$ vanishes for $s\to 0$. If $\chi,\psi,\sigma$ are as in Lemma \ref{vanish1}, then $z_s$ is an element of the non-generic irreducible subrepresentation of $I_H(\psi^{-1},\chi^{-1})$, which is the kernel of the $M$ operator if $L(\tfrac{\psi}{\chi},2s+1)$ has a pole at $s=0$. Hence in both cases the limit is well defined and depends only on $f$. Moreover, the first formula of Proposition \ref{zetino} implies that $\zita_{\chi,\psi}(f\otimes {\varphi_0})\neq 0$, for some nice choice of $f$, e.g. for $f=F_{\phi_0}$ one can see this from the computation in the proof of Theorem \ref{thmzita}, where we show
\begin{displaymath}
\zita_{\chi,\psi}(F_{\phi_0}\otimes {\varphi_0})=L(\chi\psi^{-1},1)^{-1}\vol(H(\Z_{\ell})),
\end{displaymath}
which is different from zero since we assumed $\chi\psi^{-1}\neq |\cdot|^{-1}$. 
\end{proof}
\end{defiprop}
The following corollary is then a straightforward consequence of the above proposition and of the multiplicity one theorem.
\begin{cor}
Take $\sigma, \psi,\chi$ as above and assume that condition $($\ref{inertmult1}$)$ holds and that $|\chi(\ell)|_{\C}\neq|\psi(\ell)|_{\C}$. Then $\zita_{\chi,\psi}$ is a basis for $\Hom_H(I_H(\chi,\psi)\otimes\sigma, \C)$.
\end{cor}

Using this specific element of $\Hom_H(I_H(\chi,\psi)\otimes\sigma, \C)$, we now prove a theorem that will play a key role in the proof of norm relations. Having fixed $\chi,\psi$, for every $\phi\in \Ss(\Q_{\ell}^2,\C)$ we let, as in Proposition \ref{propdefF},
\begin{displaymath}
F_{\phi}:=F_{\phi,\chi,\psi}=f_{\hat{\phi},\chi,\psi}\in I_H(\chi,\psi). 
\end{displaymath}
We also recall the special elements $\phi_0,\phi_1\in \Ss(\Q_{\ell}^2,\C)$ as in Definition \ref{defischwt}.

\begin{thm}\label{thmzita}
With notation as above, we assume that $($\ref{inertmult1}$)$ holds and that the characters $\chi,\psi$ are as follows
\begin{itemize}
\item $\chi=|\cdot |^{1/2+k+h}\cdot \tau$, for $\tau$ a finite order unramified character and $k,h\geq 0$ integers;
\item $\psi=|\cdot |^{-1/2+h}$.
\end{itemize}
Then \underline{for any $\zita$} $\in \Hom_H(I_H(\chi,\psi)\otimes\sigma, \C)$ we have
\begin{itemize}
\item[(i)]$\zita(F_{\phi_1}\otimes \varphi_0)=\frac{1}{(\ell +1)}\cdot \left( 1 - \tfrac{\ell^k}{\tau(\ell)}\right)\cdot \zita(F_{\phi_0}\otimes \varphi_0)$;
\item[(ii)]$\zita(F_{\phi_1}\otimes U(\ell) \cdot \varphi_0)=\frac{\ell}{(\ell +1)}\cdot \left[ \left( 1 - \tfrac{\ell^k}{\tau(\ell)}\right) - L(\as(\sigma), h)^{-1}\right] \cdot \zita(F_{\phi_0}\otimes \varphi_0)$,
\end{itemize}
where in (ii) the Hecke operator $U(\ell)$ is the one of Examples \ref{Uellinert} and \ref{Uell}, in the inert and split prime case respectively.
\begin{proof}
We will prove both the statements for the specific function $\zita_{\chi,\psi}$, which is a basis of $\Hom_H(I_H(\chi,\psi)\otimes\sigma, \C)$.

First we notice that $F_{\phi_t}$ is the value at $s=0$ of the Siegel section $f_{\hat{\phi_t},\chi_s,\psi_s}$. We apply Proposition \ref{propinterchangschw} and find
\begin{displaymath}
M(f_{\hat{\phi_t},\chi_s,\psi_s})=L(\chi\psi^{-1},1-2s)^{-1} f_{\phi_t, \psi_s,\chi_s}.
\end{displaymath}
We then apply the definition of the pairing $\langle -,-\rangle$ and get
\begin{align*}
\langle M(f_{\hat{\phi_t},\chi_s,\psi_s}),z_{\varphi_0} \rangle &= L(\chi\psi^{-1},1-2s)^{-1}\cdot \int_{H(\Z_{\ell})}f_{\phi_t, \psi_s,\chi_s}(h)z_{\varphi_0}(h)dh\\&=L(\chi\psi^{-1},1-2s)^{-1}f_{\phi_t, \psi_s,\chi_s}(1)\cdot \int_{K_0(\ell^t)}z_{s,\varphi_0}(h)dh.
\end{align*}
For the last equality we used that $f_{\phi_t, \psi_s,\chi_s}$ restricted to $H(\Z_{\ell})$ is a scalar multiple of $\ch(K_0(\ell^t))$. This follows from the fact that, by Lemma \ref{lemmavalscw}, $f_{\phi_t, \psi_s,\chi_s}$ restricted to $H(\Z_{\ell})$ is supported on $K_0(\ell^t)$ and $\phi_t$ is invariant by the action of $K_0(\ell^t)$. Recall that $g\cdot \varphi_0=\varphi_0$ for any $g\in G(\Z_{\ell})$. Hence $z_{s,\varphi_0}(h)=z_{s,\varphi}(1)$ for any $h\in H(\Z_{\ell})$ and we can continue the chain of equality writing
\begin{align*}
&\zita(F_{\phi_t}\otimes \varphi_0)=L(\chi\psi^{-1},1)^{-1}f_{\phi_t, \psi,\chi}(1)\vol(K_0(\ell^t))\cdot\lim_{s\to 0}  \left(L(\tfrac{\psi}{\chi}, 1+2s)z_{s,\varphi_0}(1)\right)\\&=\begin{cases}
L(\chi\psi^{-1},1)^{-1}L(\psi\chi^{-1},1)^{-1}\vol(K_0(\ell)) &\text{ if } t=1 \\
L(\chi\psi^{-1},1)^{-1}\vol(H(\Z_{\ell})) &\text{ if } t=0, \\
\end{cases} 
\end{align*}
where we applied Lemma \ref{lemmavalscw} for the value $f_{\phi_t, \psi,\chi}(1)$ and the first formula of Proposition \ref{zetino} to show that the limiting value is exactly equal to 1. Since $\frac{\vol(K_0(\ell))}{H(\Z_{\ell})}=[H(\Z_{\ell}):K_0(\ell)]^{-1}=\tfrac{1}{\ell + 1}$ and 
\begin{displaymath}
L(\psi\chi^{-1},1)=L(|\cdot|^{-1-k}\tau^{-1},1)=( 1 - \tfrac{\ell^k}{\tau(\ell)})^{-1}
\end{displaymath} we obtain (i). We proceed similarly to get (ii), using in addition the second formula of Proposition \ref{zetino}. We find
\begin{align*}
\zita_{\chi,\psi}(F_{\phi_1},U(\ell)\varphi_0)&=\vol(K_0(\ell))L(\chi\psi^{-1},1)^{-1}L(\psi\chi^{-1},1)^{-1}\cdot \lim_{s\to 0}L(\tfrac{\psi}{\chi}, 1+2s)z_{U(\ell)\varphi_0}(1)\\&=\tfrac{\ell^{3/2}}{\eta(\ell)}\vol(K_0(\ell))L(\chi\psi^{-1},1)^{-1}\lim_{s\to 0}L(\tfrac{\psi}{\chi},1+2s)^{-1}L(\tfrac{\psi}{\chi}, 1+2s)(z_{\varphi_0}(1)-L(\as(\sigma\otimes \eta),s+\tfrac{1}{2})^{-1})\\&=\ell \vol(K_0(\ell))L(\chi\psi^{-1},1)^{-1}\left[L(\tfrac{\psi}{\chi},1)^{-1}-L(\as(\sigma\otimes \eta),\tfrac{1}{2})^{-1}\right].
\end{align*} 
Using the formula proved above for the value $\zita(F_{\phi_0}\otimes \varphi_0)$ and noticing that $L(\as(\sigma\otimes \eta),\tfrac{1}{2})=L(\as(\sigma\otimes |\cdot|^{-1/2+h}),\tfrac{1}{2})=L(\as(\sigma),h)$, we obtain (ii).
\end{proof}
\end{thm}

\begin{rmk}
We emphasise that, in order to prove this theorem for any $\zita \in \Hom_H(I_H(\chi,\psi)\otimes\sigma, \C)$, we used
\begin{itemize}
\item $\sigma$ is a principal series representation for $G$ with central character such that $\chi\psi\cdot\chi_{\sigma}=1$ and $($\ref{inertmult1}$)$ holds;
\item $\chi,\psi$ are in the form $\chi=|\cdot |^{1/2+k+h}\cdot \tau$, $\psi=|\cdot |^{-1/2+h}$;
\item $\dim\left( \Hom_H(I_H(\chi,\psi)\otimes\sigma, \C) \right)=1$.
\end{itemize}
\end{rmk}

\subsection{From \texorpdfstring{$\Hom_H(\tau\otimes\sigma, \C)$ to $\X(\tau,\sigma^{\vee})$}{Hom_GL2 to X}} Let $\tau,\sigma$ be smooth representations of $H(\Q_{\ell})$ and $G(\Q_{\ell})$ respectively. We will now establish a bijection from the space $\Hom_H(\tau\otimes\sigma, \C)$ and the space $\X(\tau,\sigma^{\vee})$ of linear maps $\Zita: \tau \otimes_{\C} \hh(G) \rightarrow \sigma^{\vee}$ satisfying certain properties. For the specific choice $\tau=\Ss(\Q_{\ell}^2,\C)$, we will prove results that are essential in the proof of the norm relations (in motivic cohomology). In particular for $\sigma$ an unramified principal series representation as above, we will use the above mentioned bijection and will be able to combine these results with Theorem \ref{thmzita}, obtaining a result that is a key point in the proof of tame norm relations (in Galois cohomology).

\begin{defi}\label{defZita}
Let $\tau,\sigma$ be smooth representations of $H(\Q_{\ell})$ and $G(\Q_{\ell})$ respectively. We define $\X(\tau,\sigma^{\vee})$ to be the space of linear maps $\Zita: \tau \otimes_{\C} \hh(G) \rightarrow \sigma^{\vee}$, which are $H(\Q_{\ell})\times G(\Q_{\ell})$-equivariant, with the actions defined as follows:
\begin{itemize}
\item $H(\Q_{\ell})$ acts trivially on $\sigma^{\vee}$ and on $\tau \otimes \hh(G)$ via
\begin{displaymath}
h\cdot (v\otimes \xi) = (h\cdot v) \otimes \xi(h^{-1}(-)).
\end{displaymath}
\item $G(\Q_{\ell})$ acts naturally on $\sigma^{\vee}$ (which is a $G(\Q_{\ell})$-representation) and on $\tau \otimes \hh(G)$ via 
\begin{displaymath}
g\cdot (v\otimes \xi) = v \otimes \xi((-)g).
\end{displaymath}
\end{itemize}
\end{defi}  

We now state explicitly the bijection we were mentioning above.

\begin{prop}\label{bij}
There is a canonical bijection between $\Hom_H(\tau\otimes\sigma, \C)$ and $\X(\tau,\sigma^{\vee})$ characterised as follows
\begin{align*}
\Hom_H(\tau\otimes & \sigma, \C) \longrightarrow \X(\tau,\sigma^{\vee}) \\&\zita \longmapsto \Zita,
\end{align*}
where $\Zita(f\otimes \xi)(F)=\zita(f\otimes (\xi\cdot F))$, for every $f\in\tau,\xi\in \hh(G)$ and $F\in \sigma$.
\begin{proof}
We start by rewriting Lemma \ref{lemmaheckealg} as
\begin{itemize}
\item[(1)] $g\cdot (\xi\cdot F)=\xi (g^{-1}(-))\cdot F$;
\item[(2)] $\xi \cdot (g\cdot F)= \xi ((-)g^{-1}) \cdot F$;
\end{itemize} 
for every $\xi\in\hh(G),F\in\sigma ,g\in G(\Q_{\ell})$.

Firstly we check that $\Zita$ is $G(\Q_{\ell})$-equivariant. By definition of the action on the smooth dual of $\sigma$, for every $g\in G(\Q_{\ell})$ and $\Phi\in \sigma^{\vee}$, $g\cdot \Phi (-)=\Phi(g^{-1}\cdot(-))$. We have
\begin{align*}
[g\cdot\Zita(f\otimes \xi)](F)&=\Zita(f\otimes \xi)(g^{-1}\cdot F)=\zita(f\otimes (\xi\cdot(g^{-1}\cdot F)))\stackrel{(2)}{=} \zita(f \otimes (\xi ((-)g) \cdot F))\\&= \Zita(f\otimes (\xi ((-)g) )(F)=\Zita(g \cdot(f\otimes\xi))(F).
\end{align*}

Then we check that $\Zita$ is $H(\Q_{\ell})$-equivariant, recalling that $H(\Q_{\ell})$ acts trivially on $\sigma^{\vee}$. For $h\in H(\Q_{\ell})$ we have
\begin{align*}
\Zita(h\cdot (f\otimes \xi))(F)&=\Zita((h\cdot f) \otimes \xi(h^{-1}(-)))(F)=\zita((h\cdot f)\otimes (\xi(h^{-1}(-))\cdot F))\\& \stackrel{(1)}{=} \zita ((h\cdot f)\otimes (h\cdot (\xi\cdot F)))\stackrel{\zita\in\Hom_H}{=}\zita(f\otimes \xi\cdot F)=\Zita(f\otimes \xi)(F).
\end{align*}

Hence $\Zita\in\X(\tau,\sigma^{\vee})$. 

The fact that this defines a bijection follows from the isomorphism
\begin{displaymath}
\Hom_G(\text{c-}\ind_H^G(\tau),\sigma^{\vee})\simeq \Hom_H(\tau\otimes\sigma, \C),
\end{displaymath}
which is essentially given by Frobenius reciprocity (see \cite[Proposition 3.8.1]{GSP4}). Here we denoted with $\text{c-}\ind_H^G(\tau)$ the compact induction. Using $\tau \otimes \hh(G)= \text{c-}\ind_H^G\tau$, one finds
\begin{displaymath}
\Hom_G(\text{c-}\ind_H^G(\tau),\sigma^{\vee})\simeq \X(\tau,\sigma^{\vee}).
\end{displaymath}
\end{proof}
\end{prop}

\begin{defi}\label{defiR'}
Let $R\in\hh(G)$. we define $R'\in\hh(G)$ by $R'(g):=R(g^{-1})$.
\end{defi}
\begin{rmk}\label{rmkR'}
It is an easy computation to check that for every $\Phi\in\sigma^{\vee}$,$F\in\sigma$, we have
\begin{displaymath}
\Phi(R\cdot F)= R'\cdot \Phi(\F).
\end{displaymath}
Indeed one one side we have,
$\Phi(R\cdot F)=\Phi\left(\int_G R(g)g\cdot F dg \right)= \int_G R(g)\Phi(g\cdot F)dg$,
using linearity of $\Phi$. On the other we find
$R'\cdot \Phi(\F)=\int_G R(g^{-1})g\Phi(F)dg=\int_G R(g^{-1})\Phi(g^{-1}\cdot F)dg$.
This integrals are equal since $G$ is unimodular. 
\end{rmk}

\begin{cor} Let $\zita \leftrightarrow \Zita$ as in the above Proposition. Let $U_1\leq U_0$ be subgroups of $G$, $f_0,f_1\in \tau$ and $g_0,g_1\in G$ such that
\begin{displaymath}
\zita(f_1,g_1\cdot F)=\zita(f_0,g_0\cdot(R\cdot F))
\end{displaymath}
for some $R\in \hh(U_0\backslash G/U_0)$ and for every $F\in\sigma^{U_0}$. Then the elements $\Zita_i:=\Zita(f_i\otimes \ch(g_iU_i))\in (\sigma^{\vee})^{U_i}$ satisfy
\begin{displaymath}
\sum_{u\in U_0/U_1}u\cdot \Zita_1 = R'\cdot \Zita_0 \in (\sigma^{\vee})^{U_0}.
\end{displaymath}
\begin{proof}
It is clear by the definition of the action of $G$ that $\Zita_i\in (\sigma^{\vee})^{U_i}$, moreover summing over quotient representatives gives also $\sum_{u\in U_0/U_1}u\cdot \Zita_1\in (\sigma^{\vee})^{U_0}.$ Writing $(\sigma^{\vee})^{U_0}=(\sigma^{U_0})^{\vee}$, we are then left to check that both the L.H.S. and the R.H.S. take the same value at every $F\in \sigma^{U_0}$. Applying the Lemma above, we find
\begin{displaymath}
R'\cdot\Zita_0(F)=\Zita_0(R\cdot F)=\Zita(f_0\otimes\ch(g_0U_0))(R\cdot F)=\zita(f_0\otimes (\ch(g_0U_0)R)\cdot F)=\vol(U_0) \zita(f_1,g_1\cdot F).
\end{displaymath}
In the last equality we used the assumption $\zita(f_1,g_1\cdot F)=\zita(f_0,g_0\cdot(R\cdot F))$ together with the fact that $g_0\cdot (R\cdot F)=R(g_0^{-1}(-))\cdot F$ and
\begin{displaymath}
\ch(g_0U_0)\star R(g)=\int_G\ch(g_0U_0)(gh)R(h^{-1})dh=\int_{g^{-1}g_0U_0}R(h^{-1})dh=\vol(U_0)R(g_0^{-1}g),
\end{displaymath}
where we obtained the last equality from the fact that $R$ is in $\hh(U_0\backslash G/U_0)$. Moreover for every $u\in U_0/U_1$ we have
\begin{displaymath}
u\cdot \Zita_1(F)=u\cdot \Zita(f_1\otimes\ch(g_1U_1))(F)=\Zita(f_1\otimes\ch(g_1U_1))(u^{-1}F)=\zita(f_1\otimes \ch(g_1U_1)((-)u)\cdot F).
\end{displaymath}
We also find that 
\begin{align*}
\sum_{u\in U_0/U_1}\ch(g_1U_1)((-)u)\cdot F&=\sum_u \int_G \ch(g_1U_1)(gu) \ g\cdot F dg=\sum_u\int_{g_1U_1}gu^{-1}\cdot F dg\\&=\sum_u \int_{U_1}g_1\cdot F dg=\vol(U_1)[U_0:U_1]\ g_1\cdot F=\vol(U_0)\ g_1\cdot F,
\end{align*}
where we used the fact that $F$ is invariant by $U_0\geq U_1$. The result follows using linearity and the above expression for $u\cdot \Zita_1(F)$.
\end{proof}
\end{cor}

We now work in the setting where 
\begin{itemize}
\item we take $\tau=\Ss(\Q_{\ell}^2,\C)$
\item we replace $\sigma^{\vee}$ by an arbitrary smooth complex representation $W$ of $G(\Q_{\ell})$.
\end{itemize}
We consider $\X(W)$ to be, similarly as above, the space of functions
\begin{displaymath}
\Zita: \Ss(\Q_{\ell}^2,\C)\otimes \hh(G)\rightarrow W
\end{displaymath}
satisfying the $H(\Q_{\ell})\times G(\Q_{\ell})$ equivariance property with actions defined as above. 

\begin{lemma}
Let $\xi\in \hh(G)$ be invariant by left translation of the principal congruence subgroup of level $\ell^T$ in $H(\Z_{\ell})$ for some $T\geq 0$. Then \underline{for any $\Zita$}$\in\X(W)$ the expression
\begin{displaymath}
\frac{1}{\vol(K_{H,1}(\ell^t))}\Zita(\phi_{1,t}\otimes \xi)
\end{displaymath}
is independent of $t\geq T$, where $K_{H,1}(\ell^t),\phi_{1,t}$ are as in Definition \ref{schwrt1t}.
\begin{proof}
This is the analogous of \cite[Lemma 3.9.2]{GSP4}. The proof carries over, we sketch it for the seek of completeness. For any $t\geq T$ we fix $J$ a set of coset representatives for the quotient $K_{H,1}(\ell^T)/K_{H,1}(\ell^t)$ such that $J$ is contained in the principal congruence subgroup of level $\ell^T$. We can write $\phi_{1,T}=\sum_{k\in J}k\cdot \phi_{1,t}$. From that, using $H(\Q_{\ell})$-equivariance of $\Zita$ and the fact that $\xi$ is invariant by the action of the principal congruence subgroup of $H$ of level $\ell^T$, we obtain
\begin{displaymath}
\Zita(\phi_{1,T}\otimes \xi)=\sum_{k\in J}\Zita(k \cdot (\phi_{1,t}\otimes (k^{-1}\cdot \xi))=\tfrac{\vol(K_{H,1}(\ell^T))}{\vol(K_{H,1}(\ell^t))}\Zita(\phi_{1,t}\otimes \xi).
\end{displaymath}
\end{proof}
\end{lemma}
\begin{defi}\label{defischinfinity}
We define $\Zita(\phi_{1,\infty}\otimes \xi)$ to be the limiting value defined by the above lemma.
\end{defi}
We now define a precise choice for $\xi$, that will be used for the definition of the Euler system classes.
\begin{defi}\label{etadef}
Let $m\geq 0$ integer and $a\in\Z_{\ell}^{\times}$, we define $\eta^{(a)}_m\in G(\Q_{\ell})$ by 
\begin{displaymath}
\eta^{(a)}_m:=\begin{cases}
\left( \SmallMatrix{1&0\\0&1},\SmallMatrix{1&\tfrac{a}{\ell^m}\\0&1}\right)\in \GL_2(\Q_{\ell})\times \GL_2(\Q_{\ell}) &\text{if $\ell$ splits} \\
\SmallMatrix{1&\delta\cdot\tfrac{a}{\ell^m}\\0&1}\in \GL_2(F_{\ell}) &\text{if $\ell$ is inert. }
\end{cases}
\end{displaymath}
In the second case we fix $\delta\in \Oo_{F_{\ell}}$ such that $F_{\ell}=\Q_{\ell}\oplus \Q_{\ell}(\delta)$ as in $\S$\ref{inertzeta}.
We will write $\eta_m=\eta^{(1)}_m$. We also let for $n\geq \max (m,1)$.
\begin{displaymath}
K_{m,n}^{(a)}:=\begin{cases}
\lbrace (g_1,g_2)\in \GL_2(\Z_{\ell})\times \GL_2(\Z_{\ell}): g_1,g_2\equiv \SmallMatrix{*&*\\0&1}) \text{ mod } \ell^n, \det g_1,\det g_2 \equiv a \text{ mod } \ell^m \rbrace &\text{if $\ell$ splits}\\
\lbrace g\in G(\Oo_{F_{\ell}}): g\equiv  \SmallMatrix{*&*\\0&1} \text{ mod } \ell^n, \det g \equiv a \text{ mod } \ell^m \rbrace &\text{if $\ell$ is inert}
\end{cases}
\end{displaymath}
\end{defi}

\begin{rmk}[On the choice of $\eta_m$] The choice of these elements in $G(\Q_{\ell})$ corresponds to the choice of the ``embedding twist'' in the original definition of the Asai-Flach classes of \cite{HMS} (and of Beilinson-Flach classes in \cite{RSCMF}). The choice of the matrices is given by something of the form $\iota(\SmallMatrix{1 & \ell^{-m}\\0&1})$ ``twisted'' by some upper triangular matrix of $G(\Z_{\ell})$ not coming from $H(\Z_{\ell})$, i.e. something of the form
\begin{displaymath}
\left(\SmallMatrix{1&a_1\\0&1},\SmallMatrix{1&a_2\\0&1}\right) \ \ \ \text{and  }\ \ \ \SmallMatrix{1&a\\0&1} \ \ \text{respectively}
\end{displaymath}
for $a_1,a_2\in \Z_{\ell}/ \ell \Z_{\ell}, a_1\neq a_2$ and $a\in\Oo_{F_{\ell}}/ (\ell \Oo_{F_{\ell}}+\Z_{\ell})$ respectively.
\end{rmk}

Recall then the Hecke operator $R=U(\ell)$ in Example \ref{Uellinert} and \ref{Uell}. Taking $K'=K_{m,n}$ we have a decomposition as left cosets as in the mentioned examples. We now denote with $U'(\ell)$ the element $R'$ (see definition \ref{defiR'}) of the Hecke algebra invariant (on the left and on the right) by $K_{m,n}$, explicitly it is
\begin{displaymath}
U'(\ell)=\frac{1}{\vol(K_{m,n})}\ch( K_{m,n}\iota(\SmallMatrix{\ell^{-1}&0\\0&1}) K_{m,n})\in \hh(K_{m,n}\backslash G / K_{m,n}).
\end{displaymath}

\begin{prop}\label{THEprop}
\underline{For any $\Zita$}$\in \X(W)$, we have
\begin{displaymath}
\Zita(\phi_{1,\infty}\otimes \ch(\eta_{m+1}K_{m,n}))=\begin{cases}
\tfrac{1}{\ell}U'(\ell)  \\
\tfrac{1}{\ell -1}(U'(\ell)-1) 
\end{cases} \cdot \Zita(\phi_{1,\infty}\otimes \ch(\eta_{m}K_{m,n})) \ \  \begin{matrix}
\text{if } m\geq 1\\
\text{if } m=0.
\end{matrix}
\end{displaymath}
\begin{proof}
First of all we notice that (similarly as in Remark \ref{rmkR'}) we have $U'(\ell)\cdot \Zita(\phi_{1,\infty}\otimes \ch(\eta_{m}K_{m,n}))=\Zita(\phi_{1,\infty}\otimes (\ch(\eta_{m}K_{m,n})\star U(\ell)))$. We moreover apply the decomposition of Examples \ref{Uellinert} and \ref{Uell} to find
\begin{displaymath}
U'(\ell)\cdot \Zita(\phi_{1,\infty}\otimes \ch(\eta_{m}K_{m,n})) = \begin{matrix}
&\text{(S) }\sum_{0\leq u,v\leq\ell-1} \Zita\left(\phi_{1,\infty}\otimes \ch\left(\eta_{m}\left(\SmallMatrix{\ell & u \\ 0 & 1},\SmallMatrix{\ell& v \\ 0 & 1}\right) K_{m,n}\right)\right)   
\\&\text{(I) }\sum_{0\leq i,j\leq\ell-1} \Zita\left(\phi_{1,\infty}\otimes \ch\left(\eta_{m}\SmallMatrix{\ell & i + \alpha j \\ 0 & 1}K_{m,n}\right)\right),
\end{matrix}
\end{displaymath}
where (S) denotes the split case and (I) the inert one.
In both cases we are going to rewrite the Hecke algebra element using the invariance of $K_{m,n}$ by $\Z_{\ell}$ translation. 

(S) In the first case
\begin{displaymath}
\eta_m \cdot \left(\SmallMatrix{\ell & u \\ 0 & 1},\SmallMatrix{\ell& v \\ 0 & 1}\right) = \left(\SmallMatrix{1 & u \\ 0 & 1},\id \right)\cdot \left(\SmallMatrix{\ell & 0 \\ 0 & 1},\SmallMatrix{\ell & v+\ell^{-m} \\ 0 & 1} \right)
\end{displaymath} 
and we have 
$g\in  \left(\SmallMatrix{1 & u \\ 0 & 1},\id \right)\cdot\left(\SmallMatrix{\ell & 0 \\ 0 & 1},\SmallMatrix{\ell & v+\ell^{-m} \\ 0 & 1} \right)K_{m,n}$ if and only if $g\in \left(\SmallMatrix{\ell & 0 \\ 0 & 1},\SmallMatrix{\ell & v+\ell^{-m} \\ 0 & 1} \right)K_{m,n}$. Hence writing $\left(\SmallMatrix{\ell & 0 \\ 0 & 1},\SmallMatrix{\ell & v+\ell^{-m} \\ 0 & 1} \right)=\iota(\SmallMatrix{\ell & 0 \\ 0 & 1})\cdot \left(\id, \SmallMatrix{1 & \tfrac{1+\ell^{m}v}{\ell^{m+1}}\\ 0 & 1}\right)$, we get
\begin{displaymath}
\ch\left(\eta_{m}\left(\SmallMatrix{\ell & u \\ 0 & 1},\SmallMatrix{\ell& v \\ 0 & 1}\right) K_{m,n}\right)(g)=\ch(\eta_{m+1}^{(1+\ell^m v)}K_{m,n})\left(\iota (\SmallMatrix{\ell & 0 \\ 0 & 1}^{-1})g\right)=\SmallMatrix{\ell & 0 \\ 0 & 1}\cdot \ch(\eta_{m+1}^{(1+\ell^m v)}K_{m,n})(g).
\end{displaymath}

(I) Similarly in the second case, writing
\begin{displaymath}
\eta_{m}\SmallMatrix{\ell & i + \delta j \\ 0 & 1}=\SmallMatrix{1&i\\0&1}\SmallMatrix{\ell & \delta(j+\ell^{-m}) \\ 0 & 1}
\end{displaymath}
and using $g\in \SmallMatrix{1&i\\0&1}\SmallMatrix{\ell & \delta(j+\ell^{-m}) \\ 0 & 1}K_{m,n}$ if and only if $g\in \SmallMatrix{\ell & \delta(j+\ell^{-m}) \\ 0 & 1}K_{m,n}$, we find
\begin{displaymath}
\ch\left(\eta_{m}\SmallMatrix{\ell & i + \delta j \\ 0 & 1}K_{m,n}\right)=\SmallMatrix{\ell & 0 \\ 0 & 1}\cdot \ch(\eta_{m+1}^{(1+\ell^m j)}K_{m,n}).
\end{displaymath}

Now we write the above sum (both in the (S) and (I) case) as
\begin{align*}
\sum_{0\leq u,v\leq\ell-1} &\Zita\left(\phi_{1,\infty}\otimes \SmallMatrix{\ell & 0 \\ 0 & 1} \cdot \ch(\eta_{m+1}^{(1+\ell^m v)}K_{m,n})\right) = \tfrac{1}{\vol(K_{H,1}(\ell^n))}\sum_{0\leq u,v\leq\ell-1} \Zita\left(\phi_{1,n}\otimes\SmallMatrix{\ell & 0 \\ 0 & 1}\cdot \ch(\eta_{m+1}^{(1+\ell^m v)}K_{m,n})\right)  \\&= \tfrac{\ell}{\vol(K_{H,1}(\ell^n))}\sum_{0\leq v\leq\ell-1} \Zita\left(\SmallMatrix{\ell & 0 \\ 0 & 1}\cdot\left(\SmallMatrix{\ell & 0 \\ 0 & 1}^{-1}\cdot\phi_{1,n}\otimes\ch(\eta_{m+1}^{(1+\ell^m v)}K_{m,n})\right)\right)  \\&= \tfrac{\ell}{\vol(K_{H,1}(\ell^n))}\sum_{0\leq v\leq\ell-1} \Zita\left(\SmallMatrix{\ell & 0 \\ 0 & 1}^{-1}\cdot\phi_{1,n}\otimes\ch(\eta_{m+1}^{(1+\ell^m v)}K_{m,n})\right)\\&= \tfrac{\ell}{\vol(K_{H,1}(\ell^n))}\sum_{0\leq v\leq\ell-1} \Zita\left(\ch(\ell^{n+1}\Z_{\ell}\times(1+\ell^n\Z_{\ell}))\otimes\ch(\eta_{m+1}^{(1+\ell^m v)}K_{m,n})\right)\\&=\sum_{0\leq v\leq\ell-1} \Zita\left(\phi_{1,\infty}\otimes\ch(\eta_{m+1}^{(1+\ell^m v)}K_{m,n})\right)
\end{align*}
In the third equality we used the fact that $\Zita$ is $H(\Q_{\ell})$ equivariant and the action on the target is trivial. The fourth one is a consequence of the definition of $\phi_{1,n}$ and the action of $H(\Q_{\ell})$ on Schwartz functions. For the last one one reasons as follows. Write $S'= \ell^{n+1}\Z_{\ell}\times(1+\ell^n\Z_{\ell})$ ans $S=\ell^{n+1}\Z_{\ell}\times(1+\ell^{n+1}\Z_{\ell})$; in particular $\ch(S)=\phi_{1,n+1}$. Now we write
\begin{displaymath}
\stab(S)=K_{H,1}(\ell^{n+1})\subset \stab(S')=\{ \SmallMatrix{a&b\\c&d}\in H(\Z_{\ell}): c\equiv 0 (\ell^{n+1}),d\equiv 1 (\ell^{n}) \}.
\end{displaymath}
We also write $\Sigma=\{ \SmallMatrix{1&0\\0&1+\ell^nk}\}_{0\leq k\leq \ell -1}$, which is a set of representatives for the quotient $\stab(S')/\stab(S)$. We can then write
\begin{displaymath}
\ch(S')=\sum_{\sigma\in\Sigma}\sigma \cdot \ch(S)=\sum_{\sigma\in\Sigma}\sigma \cdot \phi_{1,n+1}.
\end{displaymath}
It is easy to check that for every $\sigma\in\Sigma$, letting $\xi_v:=\ch(\eta_{m+1}^{(1+\ell^m v)}K_{m,n})$ we have $\sigma\cdot \xi_v=\xi_v$. Hence we find
\begin{displaymath}
\Zita(\ch(S')\otimes \xi_v)=\sum_{\sigma}\Zita(\sigma \cdot \phi_{1,n+1}\otimes \xi_v)=\sum_{\sigma}\Zita(\sigma \cdot( \phi_{1,n+1}\otimes \sigma^{-1}\cdot\xi_v))=\sum_{\sigma}\Zita(\phi_{1,n+1}\otimes \xi_v)=\ell\cdot \Zita(\phi_{1,n+1}\otimes \xi_v).
\end{displaymath}
Hence we can write 
\begin{align*}
\tfrac{\ell}{\vol(K_{H,1}(\ell^n))}\sum_{0\leq v\leq\ell-1} \Zita\left(\ch(S')\otimes \xi_v\right)&=\tfrac{\ell^2}{\vol(K_{H,1}(\ell^n))}\sum_{0\leq v\leq\ell-1} \Zita\left(\phi_{1,n+1}\otimes \xi_v\right)\\&=\tfrac{1}{\vol(K_{H,1}(\ell^{n+1}))}\sum_{0\leq v\leq\ell-1} \Zita\left(\phi_{1,n+1}\otimes \xi_v\right)\\&=\sum_{0\leq v\leq\ell-1} \Zita\left(\phi_{1,\infty}\otimes\ch(\eta_{m+1}^{(1+\ell^m v)}K_{m,n})\right),
\end{align*}
where for the second equality we used $[K_{H,1}(\ell^{n}):K_{H,1}(\ell^{n+1})]=\ell^2$.

Now we notice that 
\begin{displaymath}
\SmallMatrix{a&0\\0&1}\eta_{m+1}\SmallMatrix{a^{-1}&0\\0&1}=\eta_{m+1}^{(a)}.
\end{displaymath}
Moreover for $a\equiv 1$ modulo $\ell^m\Z_{\ell}$, $\SmallMatrix{a^{-1}&0\\0&1}\in K_{m,n}$ and hence
\begin{displaymath}
\ch(\eta_{m+1}K_{m,n})\left( \SmallMatrix{a^{-1}&0\\0&1}(-)\right)=\ch(\SmallMatrix{a&0\\0&1}\eta_{m+1}\SmallMatrix{a^{-1}&0\\0&1}K_{m,n})=\ch(\eta_{m+1}^{(a)}K_{m,n}).
\end{displaymath}
Moreover, for such $a$'s, we have $\SmallMatrix{a&0\\0&1}\cdot \phi_{1,t}=\phi_{1,t}$ hence we can write
\begin{displaymath}
\Zita(\phi_{1,\infty}\otimes \ch(\eta_{m+1}^{(a)}K_{m,n}) )=\Zita\left(\phi_{1,\infty}\otimes \ch(\eta_{m+1}K_{m,n})\left( \SmallMatrix{a^{-1}&0\\0&1}(-)\right)\right)=\Zita(\phi_{1,\infty}\otimes\ch(\eta_{m+1}K_{m,n}))
\end{displaymath}
\fbox{$m\geq1$} Applying what we wrote above, we get that all the terms in the sum are equal to $\Zita(\phi_{1,\infty}\otimes\ch(\eta_{m+1}K_{m,n}))$ and hence 
\begin{displaymath}
U'(\ell)\cdot \Zita(\phi_{1,\infty}\otimes \ch(\eta_{m}K_{m,n}))=\ell \cdot \Zita(\phi_{1,\infty}\otimes\ch(\eta_{m+1}K_{m,n}))
\end{displaymath}
\fbox{$m= 0$} We can apply the same reasoning for all $v$ but for $v\equiv -1$ modulo $\ell$. For such $v$ we find $\ch(\eta_{1}^{(1+v)}K_{0,n})=\ch(K_{0,n})$, hence we get  
\begin{displaymath}
U'(\ell)\cdot \Zita(\phi_{1,\infty}\otimes \ch(K_{0,n}))=(\ell-1) \cdot \Zita(\phi_{1,\infty}\otimes\ch(\eta_{1}K_{0,n}))+\Zita(\phi_{1,\infty}\otimes\ch(K_{0,n}))
\end{displaymath}
\end{proof}
\end{prop}

We now want to go back to the case where $W$ is the smooth dual of  principal series representation and $\tau=I_H(\chi,\psi)$ in order to use the bijection of Proposition \ref{bij} and Theorem \ref{thmzita}. First of all let $K=G(\Z_{\ell})$. We assume that the Haar measures on $G(\Q_{\ell})$ and on $H(\Q_{\ell})$ are normalised to that $\vol(G(\Z_{\ell}))=\vol(H(\Z_{\ell}))=1$. We also recall the Siegel section map used above
\begin{align*}
\Ss(\Q_{\ell}^2 &,\C) \rightarrow I_H(\chi,\psi) \\& \phi \mapsto F_{\phi,\chi,\psi}:=f_{\hat{\phi},\chi,\psi},
\end{align*}
that is $H(\Q_{\ell})$ equivariant if $\chi,\psi$ are unramified.
\begin{cor}\label{THEcor}
Let $W=\sigma^{\vee}$ for $\sigma$ a principal series representation with central character $\chi_{\sigma}$. Let $\chi,\psi$ unramified characters such that
\begin{itemize}
\item $\chi=|\cdot |^{1/2+k+h} \tau$, for $\tau$ a finite order character (that may be ramified) and $k,h\geq 0$ integers;
\item $\psi=|\cdot |^{-1/2+h}$;
\item we assume that $\sigma$ satisfies $\chi\psi\cdot\chi_{\sigma}=1$ and $($\ref{inertmult1}$)$ holds.
\end{itemize}
Let $\Zita\in \X(\sigma^{\vee})$ and assume that it \underline{factors through the Siegel section map} for the above $\chi,\psi$, i.e.
\[
\begin{tikzcd}[column sep=1.5em]
\Ss(\Q_{\ell}^2,\C)\otimes \hh(G) \arrow{dr}{} \arrow{rr}{\Zita} && \sigma^{\vee} \\
& I_H(\chi,\psi)\otimes \hh(G) \arrow[ur,dotted]
\end{tikzcd}
\]
Then we have
\begin{displaymath}
\Zita(\phi_{1,\infty} \otimes (\ch(K)-\ch(\eta_1K)))=\tfrac{\ell}{\ell-1}L(\as(\sigma),h)^{-1}\cdot \Zita(\phi_0,\ch(K)),
\end{displaymath}
where $\phi_0$ is as in Definition \ref{defischwt}.
\begin{proof}
We write $\phi_{0,1}:=\ch(\ell\Z_{\ell} \times \Z_{\ell}^{\times})$ and
\begin{displaymath}
K_0:=\stab(\ell\Z_{\ell} \times \Z_{\ell}^{\times})=\lbrace \gamma \in H(\Z_{\ell}) : \gamma \equiv \SmallMatrix{*&*\\0&*} \ \text{mod } \ell\rbrace.
\end{displaymath}
We also recall that $\phi_{1,1}=\ch(\ell\Z_{\ell} \times (1+\ell\Z_{\ell}))$ and write $K_1:=K_{H,1}(\ell)=\stab(\ell\Z_{\ell} \times (1+\ell\Z_{\ell}))$. For every $\sigma\in H(\Z_{\ell}\subset K=G(\Z_{\ell})$, we have $\sigma\cdot \ch(K)=\ch(K)$ and hence
\begin{displaymath}
\Zita(\sigma\cdot \phi_{1,1}\times\ch(K))=\Zita(\sigma\cdot (\phi_{1,1}\times\ch(K)) )=\Zita(\phi_{1,1}\times\ch(K)).
\end{displaymath}
Applying this and writing $\phi_{0,1}=\sum_{\sigma\in K_0/K_1}\sigma \cdot \phi_{1,1}$, we obtain
\begin{equation}\label{eq1cor}
\begin{aligned}
\Zita(\phi_{1,\infty} \otimes \ch(K))&=\tfrac{1}{\vol(K_1)}\Zita(\phi_{1,1} \otimes \ch(K))=\tfrac{1}{[K_0:K_1]\vol(K_1)}\Zita(\phi_{0,1} \otimes \ch(K))\\&=\tfrac{1}{\vol(K_0)}\Zita(\phi_{0,1} \otimes \ch(K))=(\ell+1)\Zita(\phi_{0,1} \otimes \ch(K)),
\end{aligned}
\end{equation}
where in the last step we used the fact that $\vol(H(\Z_{\ell}))=1$ and $[H(\Z_{\ell}):K_0]=\ell+1$. 

Then applying the previous proposition we find, writing $K_{G,1}:=K_{0,1}$,
\begin{displaymath}
\Zita(\phi_{1,\infty}\otimes \ch(\eta_{1}K_{G,1}))=
\tfrac{1}{\ell -1}(U'(\ell)-1) \cdot \Zita(\phi_{1,\infty}\otimes \ch(K_{G,1})).
\end{displaymath}
Let $K_{G,0}$ be the subgroup of $G$ given by matrices congruent to $\SmallMatrix{*&*\\0&*}$ modulo $\ell$. Next we sum on both side of the last equality over representatives of $K/K_{G,1}$. On the left hand side we obtain $\Zita(\phi_{1,\infty}\otimes \eta_1\ch(K))$. On the right hand side, writing $K/K_{G,1}=K/K_{G,0}\cdot K_{G,0}/K_{G,1}$, and using the fact that $K_{G,0}/K_{G,1}$ commutes with the Hecke operator $U'(\ell)$, we obtain
\begin{displaymath}
\frac{1}{\ell -1} \sum_{\gamma\in K/K_{G,0}} \gamma \cdot (U'(\ell)-1)\Zita(\phi_{1,\infty}\otimes \ch(K_{G,0})).
\end{displaymath}
Moreover we can argue as before, using the fact that $\sigma\cdot \ch(K_{G,0})=\ch(K_{G,0})$ for $\sigma\in K_0\subset K_{G,0}$, we can rewrite $\Zita(\phi_{1,\infty}\otimes \ch(K_{G,0}))$ as $(\ell+1)\Zita(\phi_{0,1}\otimes \ch(K_{G,0}))$. Overall we have obtained
\begin{equation}\label{eq2cor}
\begin{aligned}
\Zita(\phi_{1,\infty}\otimes \eta_1\ch(K))&=\frac{\ell+1}{\ell-1}\sum_{\gamma\in K/K_{G,0}} \gamma \cdot (U'(\ell)-1)\Zita(\phi_{0,1}\otimes \ch(K_{G,0}))\\&=\frac{\ell+1}{\ell-1}\sum_{\gamma\in K/K_{G,0}} \gamma \cdot U'(\ell)\Zita(\phi_{0,1}\otimes \ch(K_{G,0}))   - \frac{\ell+1}{\ell-1}\Zita(\phi_{0,1}\otimes \ch(K)).
\end{aligned}
\end{equation}
Combining (\ref{eq1cor}) and (\ref{eq2cor}), one obtains
\begin{equation}\label{eq3cor}
\begin{aligned}
\Zita(\phi_{1,\infty} \otimes (\ch(K)-\ch(\eta_1K)))&=(\ell+1)(1+\tfrac{1}{\ell-1})\Zita(\phi_{0,1}\otimes \ch(K))\\&- \tfrac{\ell+1}{\ell-1}\sum_{\gamma\in K/K_{G,0}} \gamma \cdot U'(\ell)\Zita(\phi_{0,1}\otimes \ch(K_{G,0})).
\end{aligned}
\end{equation}

We finally use the assumption that $\Zita$ factors through the Siegel section. First we suppose that $\tau$ is ramified. Since both $\phi_0$ and $\phi_{0,1}$ are invariant under the action of matrices of the form $\SmallMatrix{a&*\\0&d}$ for $a,d\in\Z_{\ell}^{\times}$, we get
\begin{displaymath}
F_{\phi_0,\chi,\psi}=\chi(a)\cdot F_{\phi_0,\chi,\psi},
\end{displaymath}
and being $\chi$ ramified, this implies that $F_{\phi_0,\chi,\psi}=0$. Similarly $F_{\phi_{0,1},\chi,\psi}=0$ and the claimed equality reads $0=0$. So we can suppose $\tau$ unramified, so that we are able to apply Theorem \ref{thmzita} (where $\phi_1$ is our $\phi_{0,1}$). Using Remark \ref{rmkR'}, the two equalities of the theorem give us
\begin{displaymath}
\Zita(\phi_{0,1}\otimes \ch(K))=\tfrac{1}{(\ell +1)}\cdot \left( 1 - \tfrac{\ell^k}{\tau(\ell)}\right)\cdot \Zita(\phi_0\otimes \ch(K)),
\end{displaymath}
\begin{displaymath}
U'(\ell)\cdot\Zita(\phi_{0,1}\otimes \ch(K_{G,0}))=\tfrac{\ell}{(\ell +1)}\cdot \left[ \left( 1 - \tfrac{\ell^k}{\tau(\ell)}\right) - L(\as(\sigma), h)^{-1}\right] \cdot \Zita(\phi_{0}\otimes \ch(K_{G,0})).
\end{displaymath}
Hence we rewrite the two terms on the right hand side of (\ref{eq3cor}) as
\begin{align*}
(\ell+1)(1+\tfrac{1}{\ell-1})\Zita(\phi_{0,1}\otimes \ch(K))=\tfrac{\ell}{\ell-1}\left( 1 - \tfrac{\ell^k}{\tau(\ell)}\right)\cdot \Zita(\phi_0\otimes \ch(K)),
\end{align*}
\begin{align*}
\tfrac{\ell+1}{\ell-1}\sum_{\gamma\in K/K_{G,0}} \gamma \cdot U'(\ell)\Zita(\phi_{0,1}\otimes \ch(K_{G,0}))&=\tfrac{\ell}{(\ell -1)}\cdot \left[ \left( 1 - \tfrac{\ell^k}{\tau(\ell)}\right) - L(\as(\sigma), h)^{-1}\right]\sum_{\gamma\in K/K_{G,0}} \gamma  \cdot \Zita(\phi_{0}\otimes \ch(K_{G,0}))\\&=\tfrac{\ell}{(\ell -1)}\cdot \left[ \left( 1 - \tfrac{\ell^k}{\tau(\ell)}\right) - L(\as(\sigma), h)^{-1}\right]\Zita(\phi_{0}\otimes \ch(K))
\end{align*}
and get the claimed equality. 
\end{proof}
\end{cor}

\begin{rmk}[Towards Asai--Flach Euler system] As anticipated in the introduction, in order to (re)define the Euler system constructed in \cite{HMS}, we will define a special map $\AF$ for $k,k'\geq 0$ integers and $0\leq j \leq \min(k,k')$ with values in $W=H^3_{\mot}(Y_G,\mathcal{D}(2))$, where $Y_G$ is the Shimura variety associated to $G$ and $\mathcal{D}$ is a motivic sheaf depending on $k,k',j$. Such map will be of ``global nature'', more precisely it is a map 
\begin{displaymath}
\AF: \Ss(\mathbb{A}_f^2, \Q) \otimes \hh(G(\A_f),\Q) \longrightarrow H^3_{\mot}(Y_G,\mathcal{D}(2))
\end{displaymath}
satisfying conditions of $H(\A_f)\times G(\A_f)$-equivariance with actions defined as in Definition \ref{defZita}. The Asai--Flach classes will be defined by images via $\AF$ of very precise elements in $\Ss(\mathbb{A}_f^2, \Q) \otimes \hh(G(\A_f),\Q)$, whose local components will be the one we considered in this section. Proving norm relations (in motivic cohomology) will turn out to be equivalent to prove relations of such classes locally at a certain prime $q$, i.e. we will be looking at a map 
\begin{displaymath}
\Zita:=(\AF)_q: \Ss(\Q_{q}^2, \Q) \otimes \hh(G(\Q_q),\Q) \longrightarrow W= H^3_{\mot}(Y_G,\mathcal{D}(2))  \in \X(W).
\end{displaymath}
In order to prove norm relations of vertical type, we will be able to apply Proposition \ref{THEprop}. While for proving ``tame norm relatitons'' the input local data will be essentially the one in Corollary \ref{THEcor}, but we have the strong assumption on $W$. We will have to apply the \'{e}tale regulator and Hochschild–-Serre spectral sequence to pass to Galois cohomology and finally take the projection to an automorphic representation of $G$ associated to an Hilbert modular form of weight $(k+2,k'+2)$. As anticipated in Remark \ref{rmkautom}, the local component at a ``good prime'' $\ell$ of this representation will be a spherical principal series representation, so we will finally be
able to apply Corollary \ref{THEcor}.
\end{rmk}
\section{Eisenstein classes for \texorpdfstring{$H=\GL_2$}{GL2}}\label{eisenstein}

\subsection{Motivic cohomology}
Let $X$ be an object in the category $Sm$ of smooth variety over a field $k \subset \C$. Then Voevodsky defined motivic cohomology as homomorphisms in the triangulated category $\dm$ of motivic complexes. For a construction of this category see \cite{voevo}; he equips it with a functor $M: Sm\to\dm$ and with a Tate object $\textbf{Q}(1)$. 
\begin{defi} The \textit{motivic cohomology} of $X$ as above is defined by
\begin{displaymath}
H^i_{\mot}(X,\textbf{Q}(j)):=\Hom_{\dm}(M(X),\textbf{Q}(j)[i]).
\end{displaymath}
\end{defi}
He is then able to identify motivic cohhomology with hypercohomology with respect to the Zariski topology, more precisely
\begin{displaymath}
H^i_{\mot}(X,\textbf{Q}(j))\simeq \mathbb{H}^i_{Zar}(X,C_{\bullet}(\Z(j))),
\end{displaymath}
where $C_{\bullet}(\Z(j))$ is the Suslin complex of sheaves in the Zariski topology (see \cite{voevo} for more details).

The idea of motives and motivic cohomology is in some sense to collect together the information coming from all Weil cohomology theories $\T$. There are regulator maps
\begin{displaymath}
\reg_{\T}: H^i_{\mot}(X,\textbf{Q}(n)) \rightarrow H^i_{\T}(X,\textbf{Q}_{\T}(n)),
\end{displaymath}
all compatible with comparison isomorphisms. For this see \cite{huber}.

We can similarly construct motivic cohomology with ``non trivial coefficient sheaves'', using the formalism of relative Chow motives of \cite{deninger}. Consider $S$ a smooth, connected, quasiprojective $k$-variety and the category of relative Chow motives over $S$, denoted with $\operatorname{CHM}(S)_{\Q}$. It is a pseudo-abelian tensor category. For any field of characteristic zero one can similarly consider $\operatorname{CHM}(S)_{L}$, defined as the pseudo-abelian envelope of $\operatorname{CHM}(S)_{\Q}\otimes L$. The objects of such category are triples $(X,p,n)$, where $X$ is a smooth projective $S$-variety of relative dimension $m$, $p$ is an idempotent element of CH$^m(X\times_S X)$ and $n\in \Z$. This category comes equipped with a contravariant functor from the category $SmPr(S)$ of smooth projective $S$-schemes
\begin{displaymath}
M: SmPr(S)\to \text{CHM}(S)_{\Q}.
\end{displaymath}
One can take $\mathcal{F}_{\T}$ the realisation of an object in $\mathcal{F}\in\text{CHM}(S)_L$ in a cohomology theory $\T$ as above. This  takes value in the category of shaves on $S$ with extra structure depending on $\T$, this is why we will often, by abuse of notation, refer to an object $\mathcal{F}$ as \emph{motivic sheaf}. In particular, if $\T=\et$ is the $p$-adic \'{e}tale cohomology and $L$ is a $p$-adic field, $\mathcal{F}_{\et}$ is a lisse \'{e}tale $L$-sheaf over $S$. The sheaves $\mathcal{F}_{\T}$ are naturally graded objects, in particular $\mathcal{F}_{\et}=\oplus_{j}$Gr$^j\mathcal{F}_{\et}$ and, if $\mathcal{F}=M(X)$, then Gr$^j\mathcal{F}_{\et}=\mathcal{H}_{\et}^j(X/S)$ the relative \'{e}tale cohomology sheaf of $X/S$ of degree $j$.   

\begin{thm}[\cite{deninger}]\label{deningerthm} Let $A/S$ be an abelian variety. Then there is a canonical decomposition in the category of relative Chow motives over $S$
\begin{displaymath}
M(A)=\bigoplus_{i=0}^{2\dim A} M^i(A),
\end{displaymath}
such that for all the realisations $\operatorname{Gr}^jM^i(A)_{\T} =0$ if $i\neq j$.
\end{thm}

We are now ready to define motivic cohomology groups with \emph{coefficients in} $\mathcal{F}=(X,p,n)$, an object in the category $\text{CHM}(S)_{Q}$ as above. We assume that the realisations of $\mathcal{F}$ are non-zero only in one degree $r$ and let
\begin{displaymath}
H^i_{\mot}(S,\mathcal{F}(j)):=p^*H^{i+r+2n}_{\mot}(X,\mathbf{Q}(j+n)).
\end{displaymath}
As above, we find regulator maps
\begin{displaymath}
\reg_{\T}: H^i_{\mot}(S,\mathcal{F}(j)) \rightarrow H^i_{\T}(S,\mathcal{F}_{\T}(j)),
\end{displaymath}
and similarly when extending to a field extension $L$.

If $\iota: S\hookrightarrow T$ is a closed immersion of codimension $d$, there is a Gysin map
\begin{equation}\label{eqGY}
\iota_*: H^i_{\mot}(S,\iota^*\mathcal{F}(j))\to H^{i+2d}_{\mot}(T, \mathcal{F}(j+d)),
\end{equation}
where one uses the existence of a pullback functor
\begin{displaymath}
\iota^*: \operatorname{CHM}(T)_L\to \operatorname{CHM}(S)_L.
\end{displaymath}

We will be interested in sheaves over the modular curve and over the Hilbert modular surface arising from universal modular abelian varieties over them. We first fix some notation. Given a Shimura datum $(\mathcal{G},X)$, we write $Y_{\mathcal{G}}$ for the inverse limit over $K$, compact open subgroups of $\mathcal{G}(\A_f)$, of the varieties $Y_{\mathcal{G}}(K):={\mathcal{G}}(\Q)\backslash \mathcal{G}(\A_f)\times X\slash K$. Similarly every time we consider a cohomology group for $Y_{\mathcal{G}}$ we mean the limit of the cohomology groups of $Y_{\mathcal{G}}(K)$. We work with ${\mathcal{G}}=H,G^*,G$, for which $X$ is the Siegel plane in the first case and two copies of the Siegel plane in the second and third ones (with action of $G^*(\Q),G(\Q)$ given by the two real embeddings $\sigma_1,\sigma_2$ on the two copies). We obtain a smooth curve $Y_H$, which is the infinite level modular curve, and smooth surfaces $Y_{G^*},Y_G$ which are infinite level Hilbert modular surfaces. They are defined over $\Q$. 
 
If $\mathcal{G}=H,G^*$, then the corresponding finite level Shimura varieties are of PEL type and, using the functor of \cite[Theorem 8.6]{ancona}, one can associate to representations of $\mathcal{G}$ a relative Chow motive over $Y_{\mathcal{G}}(K)$ for any sufficiently small level $K$. We recall some details of this construction.

\subsubsection{Modular curves}
We consider $\mathcal{E}\rightarrow Y$, where $Y=Y_H(K)$ is the modular curve of level $K$, $K$ is a sufficiently small open compact of $H(\A_f)$ and $\mathcal{E}$ is the universal elliptic curve over $Y$. 
\begin{defi}
Let $k\geq 0$ be an integer. We define $\tsym^{k}\hh_{L}(\mathcal{E})$ to be the object of $\operatorname{CHM}(Y)_L$ given by the $k$-th symmetric power of $M^1(\mathcal{E})(1)$, where $M^1(\mathcal{E})$ is given by the decomposition of Theorem \ref{deningerthm} and $(1)$ denotes the twist by the Tate object $(Y,\operatorname{id},1)$.
\end{defi}

\subsubsection{Hilbert modular surfaces}\label{sechilbertsheaves}
Similarly, consider $\mathcal{A}\rightarrow Y^*$, where $Y^*=Y_{G^*}(K^*)$ is the Hilbert modular surface of level $K^*$, $K^*$ is a sufficiently small open compact of $G^*(\A_f)$ and $\mathcal{A}$ is the universal abelian surface over $Y^*$. Recall that $\Oo_F$ acts on $\mathcal{A}$ by endomorphisms. Consider the following object of $\operatorname{CHM}(Y^*)_L$
\begin{displaymath}
\hh_L(\mathcal{A})=M^3(\mathcal{A})(2).
\end{displaymath} 
By enlarging $L$ if necessary, we assume we have two non-zero embeddings $\theta_i:F\hookrightarrow L$. Then the object considered above decomposes as 
\begin{displaymath}
\hh_L(\mathcal{A})=\hh_L(\mathcal{A})^{(1)}\oplus \hh_L(\mathcal{A})^{(2)},
\end{displaymath}
where $\hh_L(\mathcal{A})^{(i)}$ is the direct summand where, for $x\in\Oo_F$, we have $[x]_*=\sigma_i(x)$ (see \cite[$\S$ 3.2b]{HMS} for more details).
\begin{defi}
Let $k,k'\geq 0$ be an integer. We define $\tsym^{[k,k']}\hh_{L}(\mathcal{A})$ to be the object of $\operatorname{CHM}(Y^*)_L$ given by 
\begin{displaymath}
\tsym^k(\hh_L(\mathcal{A})^{(1)})\otimes \tsym^{k'}(\hh_L(\mathcal{A})^{(2)})
\end{displaymath}
\end{defi}
One can similarly define $\sym^{(k,k')}\hh_{L}(\mathcal{A})$ and we have that its dual is $\tsym^{[k,k']}\hh_{L}(\mathcal{A})$.\\

From these sheaves, one constructs motivic sheaves for Hilbert modular surfaces $Y_G$ with respect to the larger group $G$ we were considering in the previous sections (and which is not of PEL type). The \'{e}tale cohomology of the \'{e}tale realisation of such sheaves will be the natural place where the Galois representations we are interested in will show up. Let us start by considering integers $k,k',t,t'$ such that $k,k'\geq 0$ and $k+2t=k'+2t'$. Write $\lambda$ for the quadruple $(k,k',t,t')$. Fix an open compact subgroup $U\subset G(\A_f)$ and consider $Y_G(U)$ and $Y_{G^*}(U\cap G^*(\A_f))$. One considers, with notation as above, the sheaf $\tilde{\hh}_L^{[\lambda]}$ over $Y_{G^*}(U\cap G^*(\A_f))$ defined by
\begin{displaymath}
\left[\tsym^k\left(\hh_L(\mathcal{A})^{(1)}\right)\otimes \det\left(\hh_L(\mathcal{A})^{(1)}\right)^t \right]\otimes \left[\tsym^{k'}\left(\hh_L(\mathcal{A})^{(2)}\right)\otimes \det\left(\hh_L(\mathcal{A})^{(2)}\right)^{t'}\right].
\end{displaymath}
Using the map $Y_{G^*}(U\cap G^*(\A_f)) \to Y_G(U)$ one can then define, pushforwarding the above sheaf and taking some component, a sheaf ${\hh}_L^{[\lambda]}$ over $Y_G(U)$. 
We refer to \cite[$\S$ 3.2c]{HMS} for more details. 

\begin{rmk}
One has similarly the dual sheaf ${\hh}_L^{(\lambda)}$ over $Y_G$. The \'{e}tale realisation of ${\hh}_L^{(\lambda)}$ is the lisse $\bar{\Q}_{\ell}$-sheaf associated to the representation of $\GL_2^{\Hom(F,\R)}$ given by 
\begin{displaymath}
\left(\sym^{k}(\std^{\vee})\otimes \text{det}^t \right) \otimes \left(\sym^{k'}(\std^{\vee})\otimes \text{det}^{t'} \right),
\end{displaymath}
where $\std^{\vee}$ is the dual of the standard two-dimensional representation of $\GL_2$. This is the sheaf $\mathscr{L}_{\xi,\ell}$ considered in \cite[$\S$5.5]{nekovar}.

The \'{e}tale realisation of ${\hh}_L^{[\lambda]}$ is the lisse $\bar{\Q}_{\ell}$-sheaf associated to the representation of $\GL_2^{\Hom(F,\R)}$ given by 
\begin{displaymath}
\left(\sym^{k}(\std)\otimes \text{det}^{-t} \right) \otimes \left(\sym^{k'}(\std)\otimes \text{det}^{-t'} \right),
\end{displaymath}
\end{rmk}

\subsubsection{Clebsch--Gordan map}
Write $Y_H$ and $Y_{G^*}$ for the Shimura curve and surface associated to $H$ and $G^*$. We have a closed embedding
\begin{displaymath}
\iota: Y_H \hookrightarrow Y_{G^*}.
\end{displaymath}
One has that the abelian variety $\iota^*(\mathcal{A})$ is canonically isomorphic to $\Oo_F\otimes_{\Z}\mathcal{E}$, compatibly with the $\Oo_F$ action. In particular both $\iota^*(\hh_L(\mathcal{A})^{(1)})$ and $\iota^*(\hh_L(\mathcal{A})^{(2)})$ can be identified with $\hh_L(\mathcal{E})$. Hence we obtain two maps
\begin{displaymath}
\tsym^{k+k'}\hh_L(\mathcal{E})\to \tsym^{k}\hh_L(\mathcal{E})\otimes \tsym^{k'}\hh_L(\mathcal{E})=\iota^*\left(\tsym^{[k,k']}\hh_L(\mathcal{A})\right),
\end{displaymath}
\begin{displaymath}
L(1)=\bigwedge^2_L \hh_L(\mathcal{E}) \to \hh_L(\mathcal{E})\otimes \hh_L(\mathcal{E}) = \iota^*\left(\tsym^{[1,1]}\hh_L(\mathcal{A})\right).
\end{displaymath}
Combining these two maps using multiplication in the symmetric tensor algebra, we find
\begin{prop}\label{CGmap}(\cite[Proposition 3.3.1]{HMS}). For any integers $k,k',j$ satisfying $0\leq j \leq \min(k,k')$, we have a morphism
\begin{displaymath}
CG_{\mot}^{[k,k',j]}: \tsym^{k+k'-2j}\hh(\mathcal{E}) \to \iota^*( \tsym^{[k,k']}\hh(\mathcal{A}))(-j).
\end{displaymath}
\end{prop} 
This is analogous to the map defined in \cite{kings} (see Corollary 5.2.2) for the $\GL_2\times\GL_2$ case. Moreover, writing $Y_G$ for the Shimura surface associated to $G$, one can use the fact that the pullback to $Y_{G^*}$ of the sheaf ${\hh}_L^{[\lambda]}$ over $Y_G$ is $\tsym^{[k,k']}\hh_{L}(\mathcal{A})(t+t')$ to find 
\begin{equation}\label{eqCG}
CG_{\mot}^{[k,k',j]}: \tsym^{k+k'-2j}\hh(\mathcal{E}) \to \iota_G^*( {\hh}_L^{[\lambda]}(-j-t-t')),
\end{equation}
where $\iota_G$ denotes the natural embedding $Y_H\xrightarrow{\iota} Y_{G^*} \to Y_G$.

\subsection{Eisenstein classes}
We now recall which are the elements in motivic cohomology that we are actually going to consider. This is \cite[$\S$ 7]{GSP4}.

Write $\Ss_0(\A_f^2,\Q)\subset\Ss(\A_f^2,\Q)$ for the subspace of functions $\phi$ satisfying $\phi(0,0)=0$. 
\begin{thm}[Eisenstein symbol maps] \

\begin{enumerate}
\item (\cite[Th\'{e}or\`{e}me 1.8]{colmez}) There is a canonical $H(\A_f)$-equivariant map
\begin{align*}
\Ss_0(\A_f^2&,\Q) \longrightarrow H^1_{\mot}(Y_H, \mathbf{Q}(1))=\Oo(Y_H)^{\times}\otimes \Q \\&\phi\longmapsto g_{\phi}
\end{align*}
characterised by the following: if $\phi=\ch((a,b)+N\hat{\Z})$ for some $N\geq , a,b\in\Q^2- N\Z^2$, then $g_{\phi}=g_{a/N,b/N}$, the Siegel unit in the notation of Kato (\cite[$\S$1.4]{Kato}).
\item (\cite[$\S$2]{poly}) Let $k\geq 1$. There is a $H(\A_f)$-equivariant map
\begin{align*}
\Ss(\A_f^2&,\Q) \longrightarrow H^1_{\mot}(Y_H, \tsym^k\hh_{\mathbf{Q}}(\mathcal{E})(1))\\&\phi\longmapsto \eis^k_{\phi},
\end{align*}
characterised by the following: the pullback of its de Rham realisation is the $\tsym^k\hh(\mathcal{E})$-valued differential 1-form $-F^{(k+2)}_{\phi}(\tau)(2\pi i dz)^k(2\pi i d\tau)$, where $F^{(k+2)}_{\phi}$ is the Eisenstein series defined as in \cite[Theorem 7.2.2]{GSP4}. 
\end{enumerate}
\end{thm}

\begin{rmk}
If $\phi=\ch((0,b)+N\hat{\Z})$, then $\eis^k_{\phi}$ is the class defined in \cite[Theorem 4.1.1]{kings}. Moreover, it is a consequence of Kronecker limit formula that if $\phi \in\Ss_0(\A_f^2,\Q)$, $d \log g_{\phi}$, which is the de Rham realisation of $g_{\phi}$, is equal to $-F_{\phi}^{(2)}(2\pi i d\tau)$. 
\end{rmk}

We will need a description of the target of these maps in terms of ``adelic induced representations''. The reader should have in mind, for the following discussion, that we are going to define classes using Eisenstein elements and, in order to apply the local results of the previous sections, it will be helpful to identify motivic cohomology with $H(\A_f)$ representations that locally look like $I_{H(\Q_{\ell})}(\chi,\psi)$. More precisely

\begin{defi} For $k\geq 0$ and $\eta$ a finite order character of $\A_f^{\times}/\Q^{\times +}$ such that $\eta(-1)=(-1)^k$, we define $I_k(\eta)$ to be the space of functions $f:H(\A_f)\to\C$ such that
\begin{displaymath}
f\left( \SmallMatrix{a&b\\0&d}g\right)=\parallel a \parallel^{k+1} \parallel d \parallel^{-1}\eta(a) f(g), \ \ \ \ \ \text{for every } g\in H(\A_f),a,d\in \A_f^{\times},b\in\A_f.
\end{displaymath}
We view it as a $H(\A_f)$ representation by right translation. For $k=0$ and $\eta=1$, we define $I^0_0(1)$ to be the subrepresentation which is the kernel of the integration over $H(A_f)/B(\A_f)$ on $I_0(1)$. 
\end{defi}

\begin{rmk}
Notice that restricting $f\in I_k(\eta)$ to $H(\Q_{\ell})$, we find an element $f_{\ell}$ in the space 
\begin{displaymath}
I_{H(\Q_{\ell})}(|\cdot|^{1/2+k}\eta_{\ell},|\cdot|^{-1/2}),
\end{displaymath}
with notation as in $\S$\ref{localrep}.
\end{rmk}

We finally relate motivic cohomology to these representations.
\begin{thm}\label{thmeis} With notation as above,
\begin{enumerate}
\item (\cite[Theorem 3]{scholl}) there is a $H(\A_f)$-equivariant isomorphism
\begin{displaymath}
\partial_0: \frac{\Oo^{\times}(Y)}{(\Q^{ab})^{\times}}\otimes \C \longrightarrow I_0^0(1) \oplus \bigoplus_{\eta\neq 1} I_0(\eta),
\end{displaymath}
characterised by the fact that $\partial_0(g)(1)$ is the order of vanishing of $g$ at the cusp $\infty$.
\item For $k\geq 1$, there is a surjective $H(\A_f)$-equivariant map
\begin{displaymath}
\partial_k: H^1_{\mot}(Y,\tsym^k\hh_{\mathbf{Q}}(\mathcal{E})(1)) \otimes \C \longrightarrow \bigoplus_{\eta} I_k(\eta),
\end{displaymath}
such that $\partial_k(x)$ is the residue at $\infty$ of the de Rham realisation of $x$. Moreover this map is an isomorphism on the image of the Eisenstein symbol. (See \cite[Theorem 7.4]{ss} and \cite[Lemma 4.3]{lemma}) 
\end{enumerate}
\end{thm}
Moreover, we have an explicit description of the image of the Eisenstein symbols via these maps. Write $\Ss(\A_f^2, \C)^{\eta}$ for the subspace of $\Ss(\A_f^2, \C)$ on which $\hat{\Z}^{\times}$ acts via the character $\eta$.
\begin{prop}\label{propeis}(\cite[Proposition 7.3.4]{GSP4}) Let $\phi\in \Ss(\A_f^2, \C)^{\eta}$ and write $\phi=\prod_{\ell}\phi_{\ell}$. If $k=0$ and $\eta=1$, assume that $\phi(0,0)=0$. Then we have
\begin{displaymath}
\partial_k(\eis^k_{\mot,\phi})=\frac{2(k+1)!L(k+2,\eta)}{(-2\pi i)^{k+2}}\prod_{\ell} f_{\widehat{\phi_{\ell}},|\cdot|^{1/2+k}\eta_{\ell},|\cdot|^{-1/2}},
\end{displaymath}
where the functions in the product are the Siegel sections of Proposition \ref{propdefF}.
\end{prop}

\section{Definition of Asai--Flach map and classes}\label{defi}
\subsection{Definition of the map}\label{defimap} We fix integers $k,k'\geq 0$ such that $k+2t=k'+2t'$ and write $\mathcal{D}^{k,k'}:={\hh}_L^{[\lambda]}(-t-t')$, as in \ref{sechilbertsheaves}. We will fix $j$ such that $0\leq j \leq \min(k,k')$. The goal of this section is to define a map 
\begin{displaymath}
\AF: \Ss(\A_f^2,\Q)\otimes \hh(G(\A_f),\Q)\longrightarrow H^3_{\mot}(Y_G,\mathcal{D}^{k,k'}(2-j))
\end{displaymath}
that is $H(\A_f)\times G(\A_f)$ equivariant, with actions given as follows
\begin{itemize}
\item $H(\A_f)$ acts trivially on the target and it acts on $\Ss(\A_f^2,\Q)\otimes \hh(G(\A_f),\Q)$ via
\begin{displaymath}
h\cdot (\phi\otimes \xi) = (h\cdot \phi) \otimes \xi(h^{-1}(-)).
\end{displaymath}
\item $G(\A_f)$ acts via the natural action on $H^3_{\mot}(Y_G,\mathcal{D}^{k,k'}(2-j))$ and on the source via 
\begin{displaymath}
g\cdot (\phi\otimes \xi) = \phi \otimes \xi((-)g).
\end{displaymath}
\end{itemize}

We will consider open compact subgroups $U\subset G(\A_f)$ such that the natural map
\begin{displaymath}
\iota_U:Y_H(U\cap H) \rightarrow Y_G(U)
\end{displaymath}
is a closed embedding. It is easy to check that this holds for $U$ sufficiently small. We then have that the Hecke algebra $\hh(G(\A_f),\Z)$ is generated as a $\Z$-module by the functions of the form $\ch(gU)$ where $g\in G(\A_f)$ and $U$ is as above. The following proposition is an adaptation of \cite[Proposition 8.2.3]{GSP4}.

\begin{prop}\label{Hsmall} Let $U$ be a sufficiently small subgroup of $G(\A_f)$ and $U'\leq U$. Write $V:=U\cap H(\A_f)$ and choose coset representative $(x_j)_{j\in J}$ of the double quotient $V\backslash U/U'$. Then write $U'_j:=x_jU'x_j^{-1}, V'_j:=H(\A_f)\cap U'_j$. Define the maps
\begin{displaymath}
\iota_j:  Y_H(V'_j) \to Y_G(U'_j) \xrightarrow{\cdot x_j} Y_G(U'),
\end{displaymath}
where the last arrow is induced by multiplication on the right by $x_j$. Then the $i_j$'s have disjoint image, hence $\sqcup Y_H(V'_j)\to Y_G(U')$ is a closed embedding. Moreover the following diagram is cartesian
\begin{displaymath}
\begin{tikzcd}
\sqcup_j Y_H(V'_j) \arrow{d}{} \arrow{r}{} & Y_G(U') \arrow{d}\\
Y_H(V) \arrow{r} & Y_G(U).
\end{tikzcd}
\end{displaymath}
\begin{proof}
First of all we notice that since $U'_j<U$, $U'_j$ is sufficiently small and hence each $i_j$ is a closed immersion. Now for each $j\in J$ we choose coset representatives $(v_i)_{i\in I_j}$ of $V/V'_j$ and get coset representatives $(v_ix_j)_{i\in I_j,j\in J}$ of $U/U'$. Hence the union of the images of the $\iota_j$'s is exactly the preimage of $Y_H(V)$. This implies that the diagram is cartesian and the top map is a closed immersion since the bottom one is. 
\end{proof}
\end{prop}

\begin{defi}
Fix an Haar measure on $H(\A_f)$ and let $V\subset H(\A_f)$ an open compact subgroup. We define a map $A_V:\Ss(\A_f^2,\Q)\rightarrow \Ss(\A_f^2,\Q)^V$ by
\begin{displaymath}
A_V(\phi):= \int_V h\cdot \phi dh = \vol(W)\cdot \sum_{v\in V \backslash W} v\cdot \phi,
\end{displaymath}
where $W$ is an open compact subgroup of $V$ fixing $\phi$.
\end{defi}
The following lemma is an immediate consequence of the definition.
\begin{lemma}\label{lemmaaverage}
If $V'\subset V$, we have $A_V(\phi)=\sum v\cdot A_{V'}(\phi)$, where $V=\bigsqcup vV'$.
\end{lemma}

Now let $x\in G(\A_f)$ and $U$ such that $xUx^{-1}$ is sufficiently small. Let
\begin{displaymath}
\xi=\ch(xU) \in \hh(G(\A_f),\Z), \ \ \ V=H\cap xUx^{-1}.
\end{displaymath}
We denote with $\iota_{xU}$ the closed embedding obtained by
\begin{displaymath}
\iota_{xU}: Y_H(V)\hookrightarrow Y_G(xUx^{-1}) \xrightarrow[\cdot x]{\simeq} Y_G(U).
\end{displaymath}
Moreover (\ref{eqCG}) gives a map
\begin{displaymath}
CG^{[k,k',j]}_{\mot}: H^i_{\mot}(Y_H(V), \tsym^{k+k'-2j}\hh_{L}(n))[j+t+t']\rightarrow  H^i_{\mot}(Y_H(V), \iota^*( {\hh}_L^{[\lambda]})(n-j-(t+t'))),
\end{displaymath}
where we added the twist by the $(j+t+t')$-th power of the determinant, meaning tensoring with the one dimensional representation on which $H(\A_f)$ acts as $(j+t+t')$-th power of the determinant.
One also has, as in (\ref{eqGY}), a pushforward map
\begin{displaymath}
(\iota_{xU})_*: H^i_{\mot}(Y_H(V), \iota_{xU}^*({\hh}_L^{[\lambda]})(n)) \rightarrow H^{i+2}_{\mot}(Y_G(U), {\hh}_L^{[\lambda]}(n+1)).
\end{displaymath}
Composing such morphisms for $i=1,n=1$ we obtain a map
\begin{displaymath}
\iota_{xU,*}^{[k,k',j]}: H^1_{\mot}(Y_H(V), \tsym^{k+k'-2j}\hh_{\Q}(1))[j+t+t']\longrightarrow H^{3}_{\mot}(Y_G(U), {\hh}_L^{[\lambda]}(2-j-(t+t'))).
\end{displaymath}
We also have, from the previous chapter, a $H(\A_f)$-equivariant map
\begin{align*}
\Ss(\A_f^2&,\Q) \rightarrow H^1_{\mot}(Y_H, \tsym^{k+k'-2j}\hh_{\Q}(1)) \\&\phi \longmapsto \eis^{k+k'-2j}_{\mot,\phi}.
\end{align*}
In particular if $\phi \in \Ss(\A_f^2,\Q)^V$ for some $V\subset H(\A_f)$, we have $\eis^{k+k'-2j}_{\mot,\phi}\in H^1_{\mot}(Y_H(V), \tsym^{k+k'-2j}\hh_{\Q}(1))$.
We can finally make the following definition:
\begin{defi}
The level $U$ motivic Asai--Flach map for $k,k',j$ and $U$ as above is defined by
\begin{align*}
\AF_{,U}: \Ss(\A_f^2,\Q)[j+t+t']\otimes \hh&(G(\A_f),\Z)\longrightarrow H^{3}_{\mot}(Y_G(U), {\hh}_L^{[\lambda]}(2-j-(t+t'))) \\ &\phi \otimes \xi \longmapsto \iota_{xU,*}^{[k,k',j]}(\eis^{k+k'-2j}_{\mot,A_V(\phi)}), 
\end{align*}
where $\xi=\ch(xU)$ as above and $V=H\cap xUx^{-1}$. Since the Hecke algebra is spanned by functions of this form, $\AF$ is defined extending by $\Z$ linearity.
\end{defi}

\begin{prop}
The above defined map satisfies
\begin{itemize}
\item[(a)] if $\xi'=g\cdot\xi$ for $g\in G(\A_f)$, then 
\begin{displaymath}
\AF_{,gUg^{-1}}(\phi\otimes\xi')=g\cdot \AF_{,U}(\phi \otimes \xi).
\end{displaymath}
\item[(b)] For every $h\in H(\A_f)$, one has
\begin{displaymath}
 \AF_{,U}(h\cdot(\phi \otimes \xi))= \AF_{,U}(\phi \otimes \xi).
\end{displaymath}
\item[(c)] If $U'\subset U$, writing $\pi:Y_G(U')\rightarrow Y_G(U)$ for the natural projection map, we find 
\begin{displaymath}
 \AF_{,U'}(\phi \otimes \xi)=\pi^* \AF_{,U}(\phi \otimes \xi).
\end{displaymath}
\end{itemize} 
\begin{proof} We prove all the statements for $\xi=\ch(xU)$, which is enough because they span the Hecke algebra.

(a) We find that $\xi'=\ch(xg^{-1}(gUg))$. Then the statement follows from the commutativity of the following diagram 
\begin{displaymath}
\begin{tikzcd}
Y_G(xUx^{-1}) \arrow{d}{\id} \arrow{r}{\cdot x} & Y_G(U) \arrow{d}{\cdot g^{-1}}\\
Y_G(xUx^{-1}) \arrow{r}{\cdot xg^{-1}} & Y_G(gUg^{-1}).
\end{tikzcd}
\end{displaymath}
together with the fact that the action of $g$ on cohomology is precisely given by the pushforward of the right vertical map. 

(b) We have, by definition of the action, $h\cdot(\phi\otimes \xi)=h\cdot \phi \otimes \ch(hxU)$. Writing $V=xU\cap H$ and using $A_{hVh^{-1}}(h\cdot \phi)=h\cdot A_V(\phi)$ we get the desired equality.

(c) We can write $\xi=\ch(xU)=\sum_{u\in U/U'} \ch(xuU')$. We then apply Proposition \ref{Hsmall} to find $\pi^*\circ\pi_*=\sum _{u\in U/U'} u$. Combining this with $\pi_*(\AF_{,U'}(\phi \otimes \ch(xU')))=\AF_{,U}(\phi \otimes \ch(xU))$ and $\sum _{u\in U/U'} u\cdot \ch(xU')=\ch(xU)$, we are done. 
\end{proof}
\end{prop}

\begin{defi} We define
\begin{displaymath}
\AF: \Ss(\A_f^2,\Q)[j+t+t']\otimes \hh(G(\A_f),\Z)\longrightarrow H^{3}_{\mot}(Y_G, {\hh}_L^{[\lambda]}(2-j-(t+t')))
\end{displaymath}
to be the direct limit $\varinjlim_{U}\AF_{,U}$. This is well defined thanks to (c) in the above Proposition and is $H(\A_f)\times G(\A_f)$-equivariant with respect to the action given above thanks to (a)-(b) in the above Proposition. 
\end{defi}
\subsection{Definition of the classes in motivic cohomology}\label{deficlass}
In order to define the Asai--Flach elements in motivic cohomology, we will specify the choice of an element in $\Ss(\A_f^2,\Q)\otimes \hh(G(\A_f),\Z)$ to which we will apply $\AF$.

We start by fixing a prime $p$, a finite set of primes $S$ not containing $p$. Our choice will also depend on  integers $m,M\geq 1$ with $M$ coprime to $S$ and $p$. We now will define $K_{M,m,n}\subset G(\A_f), W\subset H(\A_f)$ and $\phi_{M,m,n}\in \Ss(\A_f^2,\Q),\xi_{M,m,n}\in\hh(G(\A_f),\Z)$ satisfying certain properties and apply $\AF$ in order to define an element 
\begin{displaymath}
z_{M,m,n}^{[k,k',j]}:= \tfrac{1}{\vol(W)}\AF (\phi_{M,m,n}\otimes\xi_{M,m,n} ) \in H^{3}_{\mot}(Y_G(K_{M,m,n}), \hh_L^{[\lambda]}(2-j-(t+t'))).
\end{displaymath}

Every definition of such data will be given in term of local data. Writing $\Q_S=\prod_{\ell\in S}\Q_{\ell}$ we will define
\begin{itemize}
\item subgroups $K_S\subset G(\Q_S),K_{p,n}\subset G(\Q_p)$ and let
\begin{displaymath}
K_n:= K_S \times K_{p,n}\times \prod_{\ell \not\in S\cup \{ p\}} G(\Z_{\ell})\subset G(\A_f);
\end{displaymath}
\item A subgroup $K_{M,m,n}\subset K_n$, defined by $K_n \cap \det^{-1}(1+Mp^m\Oo_F)$;
\item functions $\phi_S\in \Ss(\Q_S^2,\Z), \phi_{\ell} \in\Ss(\Q_{\ell}^2,\Z)$ for $\ell\not\in S$ and let
\begin{displaymath}
\phi_{M,m,n}:=\phi_S\otimes \bigotimes_{\ell\not\in S}\phi_{\ell};
\end{displaymath}
\item elements $\xi_{\ell}\in \hh(G(\Q_{\ell}),\Z)$ for $\ell\not\in S$ and let
\begin{displaymath}
\xi_{M,m,n}:=\ch(K_S)\otimes  \bigotimes_{\ell\not\in S}\xi_{\ell};
\end{displaymath}
\item an open compact subgroup $W\subset H(\A_f)$ defined choosing $W_S\subset H(\Q_S)\cap K_S$ acting trivially on $\phi_S$ and $W_{\ell}\subset H(\Q_{\ell})$ for $\ell\not\in S$ and letting
\begin{displaymath}
W:= W_S \times \prod_{\ell\not\in S}W_{\ell}.
\end{displaymath}
\end{itemize}
We consider fixed the choices at $S$ and require that the global elements satisfy the following
\begin{itemize}
\item[(i)] $\xi_{M,m,n}$ is fixed by right translation of $K_{M,m,n}$,
\item[(ii)] $\xi_{M,m,n}$ is fixed by left translation of $W$, 
\item[(iii)] $\phi_{M,m,n}$ is stable under the action of $W$.
\end{itemize}
We first define the level subgroup $K_n$. We are only left with saying what is the choice at $p$. We let
\begin{displaymath}
K_{p,n}:= \lbrace \SmallMatrix{a&b\\c&d}\in G(\Z_p): c\equiv d-1\equiv 0 \text{ mod } p^n\rbrace.
\end{displaymath}
The desired subgroup $K_{M,m,n}$ will then be given at $p$ by 
\begin{displaymath}
\lbrace g\in K_{p,n}: \det g  \equiv 1 \text{ mod }Mp^m\hat{\Z}\rbrace.
\end{displaymath}
Write $K_{M,m,n}^*:=K_{M,m,n}\cap G^*(\A_f)\subset K_n^*:=K_n\cap G^*(\A_f)$. We then have (\textit{cf.} \cite[Proposition 5.4.2]{GSP4}) that the determinant map induces an isomorphism
\begin{equation}\label{eqDET}
Y_{G^*}(K^*_{M,m,n})\simeq Y_{G^*}(K^*_n)\times_{\Q} \mu_{Mp^m}. 
\end{equation}
We now define the local terms of $\xi_{M,m,n},\phi_{M,m,n},W$ at places $\ell\not\in S$, dividing the three cases $\ell\nmid Mp, \ell\mid M,\ell=p$. First of all we write, for $r\geq 0$,
\begin{displaymath}
\eta_{\ell,r}:=\begin{cases}
\left( \SmallMatrix{1&0\\0&1},\SmallMatrix{1&\tfrac{1}{\ell^r}\\0&1}\right)\in \GL_2(\Q_{\ell})\times \GL_2(\Q_{\ell}) &\text{if $\ell$ splits} \\
\SmallMatrix{1& \delta\tfrac{1}{\ell^r}\\0&1}\in \GL_2(F_{\ell}) &\text{if $\ell$ is inert. }
\end{cases}
\end{displaymath}
This is the element $\eta_{r}$ defined in Definition \ref{etadef}.

\underline{$\ell\nmid Mp$}: We let
\begin{displaymath}
\xi_{\ell}=\ch(G(\Z_{\ell})), \ \ \ W_{\ell}=H(\Z_{\ell}), \ \ \ \phi_{\ell}=\ch(\Z_{\ell}^{2}).
\end{displaymath}

\underline{$\ell\mid M$}: First we define $K_{\ell,1}=\{g\in G(\Z_{\ell}): \det g \equiv 1 \text{ mod } \ell\}$. We then let
\begin{displaymath}
\xi_{\ell}=\ch(K_{\ell,1})-\ch(\eta_{\ell,1}K_{\ell,1}),
\end{displaymath}
\begin{displaymath}
W_{\ell}=\lbrace h\in H(\Z_{\ell}): \det h \equiv 1 \text{ mod } \ell, h\equiv \SmallMatrix{*&*\\0&1} \text{ mod } \ell^2\rbrace,
\end{displaymath}
\begin{displaymath}
\phi_{\ell}=\ch(\ell^2\Z_{\ell} \times (1+\ell^2\Z_{\ell})).
\end{displaymath}

\underline{$\ell=p$}: We define $K_{p,m,n}=\lbrace g\in G(\Z_{p}): \det g\equiv 1 \text{ mod } p^m, g\equiv \SmallMatrix{*&*\\0&1} \text{ mod } p^n\rbrace$. Let
\begin{displaymath}
\xi_p=\ch(\eta_{p,m}K_{p,m,n}).
\end{displaymath}
We then choose an integer $t\geq 1$ big enough such that $W_p\subset \eta_{p,m}K_{p,m,n}\eta_{p,m}^{-1}$, where
\begin{displaymath}
W_p=\lbrace h\in H(\Z_{\ell}): \det h \equiv 1 \text{ mod } p^m, h\equiv \SmallMatrix{*&*\\0&1} \text{ mod } p^t \rbrace.
\end{displaymath}
Finally for such choice of $t$ we let
\begin{displaymath}
\phi_p=\ch(p^t\Z_{p} \times (1+p^t\Z_{p})).
\end{displaymath}

It follows easily from the definitions that conditions (i),(ii),(iii) above are satisfied. We finally can, as anticipated above, make the following definition:
\begin{defi}\label{motivicdefi} For $M,m,n$ and $W,\phi_{M,m,n},\xi_{M,m,n}$ as above, we define
\begin{displaymath}
z_{M,m,n}^{[k,k',j]}:= \tfrac{1}{\vol(W)}\AF (\phi_{M,m,n}\otimes\xi_{M,m,n} ) \in H^{3}_{\mot}(Y_G,\hh^{[\lambda]}(2-j-(t+t')))
\end{displaymath}
\end{defi}

\begin{lemma}
The above definition is independent on the choice of the Haar measure on $H(\A_f)$ and on the choice of $t$ at the place $p$. 
\begin{proof}
Writing $U:=K_{M,m,n}$, we have, from (i) and (ii) that $\xi_{M,m,n}\in \hh(W\backslash G /U)$. We rewrite it as
\begin{displaymath}
\xi_{M,m,n}=\sum \ch (x_i U),
\end{displaymath}
where $\ch (x_i U)$ is left invariant under $W$, i.e. $W\subset V_i:=H(\A_f)\cap x_iUx_i^{-1}$. Hence writing $\iota_i:=\iota_{x_iU}$ we have that by definition our classes are
\begin{displaymath}
\tfrac{1}{\vol(W)}\sum_i (\iota_i)_* \circ CG_{\mot}^{[k,k',j]} (\eis^{k+k'-2j}_{\mot,A_{V_i}(\phi)}),
\end{displaymath}
where $\phi=\phi_{M,m,n}$. Using (iii), the definition of the averaging map and the fact that $\eis^{k+k'-2j}_{\mot,-}$ is $H(\A_f)$-equivariant, we can write $\eis^{k+k'-2j}_{\mot,A_{V_i}(\phi)}=\vol(W)\sum_{v\in V_i/W}v\cdot \eis^{k+k'-2j}_{\mot,\phi}$, from which the independence on the Haar measure becomes clear. 

Write now $W$ for the subgroup defined by the condition at $p$ with a fixed choice of $t$ and $W^0$ for the subgroup defined with a different choice, say $t_0> t$. We similarly write $\phi_{M,m,n}$ and $\phi^0_{M,m,n}$.  We can write
\begin{displaymath}
\phi_{M,m,n}=\sum_{w\in W/W^0}w\cdot \phi^0_{M,m,n}.
\end{displaymath}
We obtain 
\begin{align*}
\tfrac{1}{\vol(W)}\AF(\phi_{M,m,n}\otimes\xi_{M,m,n})&=\tfrac{1}{\vol(W)}\sum_w \AF(w\cdot \phi^0_{M,m,n}\otimes\xi_{M,m,n})\\&=\tfrac{1}{\vol(W)}\sum_w \AF(w\cdot (\phi^0_{M,m,n}\otimes\xi_{M,m,n}))\\&=\tfrac{[W:W^0]}{\vol(W)}\AF(\phi^0_{M,m,n}\otimes\xi_{M,m,n})\\&=\tfrac{1}{\vol(W^0)}\AF(\phi^0_{M,m,n}\otimes\xi_{M,m,n}).
\end{align*}
In the second equality we used (ii), and in the third the fact that $\AF$ is $H(\A_f)$-equivariant. 
\end{proof}
\end{lemma}
\begin{rmk}\label{rmkHMS} Condition (i) together with the fact that $\AF$ is $G(\A_f)$-equivariant implies
\begin{displaymath}
z_{M,m,n}^{[k,k',j]}\in H^{3}_{\mot}(Y_G(K_{M,m,n}),\hh^{[\lambda]}(2-j-(t+t'))),
\end{displaymath}
as wanted. Pulling back these classes to $Y_{G^*}$, we find elements in 
\begin{displaymath}
H^{3}_{\mot}(Y_{G^*}(K_{n})\times \mu_{Mp^m},\tsym^{[k,k']}\hh(\mathcal{A})(2-j)).
\end{displaymath}
It is clear from the construction that, if $M=1,m=0$, these are the classes $\text{AF}_{\mot,p^n}^{[k,k',j]}$ of \cite[Definition 3.4.2]{HMS}. In \textit{loc.cit.}, the authors obtain classes $\text{AF}_{\mot,Mp^m,p^n,a}^{[k,k',j]}$ twisting the class $\text{AF}_{\mot,Mp^mp^n}^{[k,k',j]}$ by the matrix $\SmallMatrix{1&\tfrac{a}{Mp^m}\\0&1}$, for any $a\in \Oo_F/(Mp^m\Oo_F+\Z)$. The input elements in $\hh(G(\A_f))$ we are considering involve matrices of the same form for a specific choice of $a\in \Oo_F\otimes\hat{\Z}/(Mp^m+\hat{\Z})$. The places where our classes and the ones of \cite[Definition 3.5.1]{HMS} differ are the primes $\ell\mid M$, where there is a correction term given by $\ch(K_{\ell,1})$. Indeed the tame norm relations we will prove are not the same as the ones obtained in \cite[Theorem 3.5.3, Corollary 4.3.8]{HMS}, where the factor appearing is the local Euler factor $P(\sigma_{\ell}^{-1} \ell^{-1-j})$ summed  with a term divisible by $\ell-1$. One then applies \cite[Lemma 7.3.4]{RSCMF} (\cite[Lemma IX.6.1]{rubin}) to find Galois cohomology classes satisfying the correct relations. Thanks to the correction term at the primes dividing $M$, we in fact don't get this extra factor and we are able to show that our classes form an Euler system.
\end{rmk}
\section{The Asai--Flach Euler systems norm relations}\label{ES}
\subsection{Pushforward compatibilities in motivic cohomology}
We now prove that the classes just defined satisfy compatibility properties if we vary the level $K_n$ and if we vary the cyclotomic field in the $p$-direction. 
\begin{thm}\label{changelevel}
For $n\geq 1$ we have, writing $\pi_n:Y_G(K_{M,m,n+1})\rightarrow Y_G(K_{M,m,n})$ for the natural projection,
\begin{displaymath}
(\pi_n)_*(z_{M,m,n+1}^{[k,k',j]})=z_{M,m,n}^{[k,k',j]}.
\end{displaymath}
\begin{proof}
Going back to the definition of the local data in $\S$\ref{deficlass}, we see that the only place where these differ is $p$, where
\begin{displaymath}
\xi_p=\ch(\eta_{p,m}K_{p,m,n}),
\end{displaymath}
while we can choose the same $t$ sufficiently large, so that we have the same $W$ and the same $\phi_p$ in the definition of $z_{M,m,n+1}^{[k,k',j]}$ and $z_{M,m,n}^{[k,k',j]}$. So locally at $p$, we need to check that
\begin{displaymath} 
(\pi_n)_*((\AF)_{p}(\phi \otimes \ch(xK_{p,m,n+1})))=(\AF)_{p}(\phi \otimes \ch(xK_{p,m,n})).
\end{displaymath}
But this is true, since we can write $\ch(xK_{p,m,n})=\sum_{k\in K_{p,m,n}/K_{p,m,n+1}}\ch(xkK_{p,m,n+1})$ and the pushforward act on cohomology by multiplication of coset representatives $k\in K_{p,m,n}/K_{p,m,n+1}$.
\end{proof}
\end{thm}
The following theorem is essentially the proof of the vertical type Euler system norm relation for the classes we will obtain in Galois cohomology in the next section starting with the motivic input $ z_{M,m,n}^{[k,k',j]}$.
\begin{thm}\label{quasivert}
For $m\geq 1$ we have, writing $\pi_m:Y_G(K_{M,m+1,n})\rightarrow Y_G(K_{M,m,n})$ for the natural projection,
\begin{displaymath}
(\pi_m)_*(z_{M,m+1,n}^{[k,k',j]})=\begin{cases}
\tfrac{U'(p)}{p^j}  \\
\left( \tfrac{U'(p)}{p^j}-1\right)  
\end{cases} \cdot z_{M,m,n}^{[k,k',j]} \ \  \begin{matrix}
\text{if } m\geq 1\\
\text{if } m=0,
\end{matrix}
\end{displaymath}
where $U'(p)$ is the Hecke operator in $\hh(K_{p,m,n}\backslash G(\Q_p)/K_{p,m,n})$ given by the double coset of $\SmallMatrix{p^{-1}&0\\0&1}$.
\begin{proof}
This theorem follows from the choice of the local data and from Proposition \ref{THEprop}. As in the previous theorem, the elements $\phi_{M,m+1,n}\otimes \xi_{M,m+1,n}$ and $\phi_{M,m,n}\otimes \xi_{M,m,n}$ are the same at places different from $p$. Hence we are comparing two values of the $p$-part map
\begin{displaymath}
(\AF)_p: \Ss(\Q_p)\otimes \hh(G(\Q_p)) \to H^{3}_{\mot}(Y_G, {\hh}_L^{[\lambda]}(2-j-(t+t')))[-j],
\end{displaymath}
which is $H(\Q_p)\times G(\Q_p)$-equivariant. Since it is enough to check the equality after tensoring with $\C$, we can apply Proposition \ref{THEprop}. Indeed, reasoning as in the proof of the previous theorem, we have, that on the left hand side we have $\tfrac{1}{\vol(W_{p,m+1})}(\AF)_p(\phi_{p}\otimes \ch(\eta_{p,m+1}K_{p,m,n}))$, where $\phi_p=\ch(p^t\Z_{p} \times (1+p^t\Z_{p}))$ and $W_{p,m+1}=\lbrace h\in H(\Z_{\ell}): \det h \equiv 1 \text{ mod } p^{m+1}, h\equiv \SmallMatrix{*&*\\0&1} \text{ mod } p^t \rbrace$. Hence the classes we need to compare are
\begin{displaymath}
[K_{H,1}(p^t):W_{p,m+1}](\AF)_p(\phi_{1,\infty}\otimes \ch(\eta_{p,m+1}K_{p,m,n})),
\end{displaymath}
\begin{displaymath}
[K_{H,1}(p^t):W_{p,m}](\AF)_p(\phi_{1,\infty}\otimes \ch(\eta_{p,m}K_{p,m,n})).
\end{displaymath}
The proposition tells us that 
\begin{displaymath}
(\AF)_p(\phi_{1,\infty}\otimes \ch(\eta_{p,m+1}K_{p,m,n}))=\begin{cases}
\tfrac{1}{p}U'(p)  \\
\tfrac{1}{p -1}(U'(p)-1) 
\end{cases} \cdot (\AF)_p(\phi_{1,\infty}\otimes\ch(\eta_{p,m}K_{p,m,n})) \ \  \begin{matrix}
\text{if } m\geq 1\\
\text{if } m=0.
\end{matrix}
\end{displaymath}
The factors $\tfrac{1}{\ell},\tfrac{1}{\ell-1}$ cancels out since $[K_{H,1}(p^t):W_{p,m}]/[K_{H,1}(p^t):W_{p,m+1}]=[W_{p,m}:W_{p,m+1}]$ is equal to $\ell,\ell-1$ respectively. Moreover the twist $[j]$ makes the factor $p^j$ appear.
\end{proof}
\end{thm}

\subsection{Hilbert cuspforms and Galois representations}
In the previous chapter, we constructed some classes
\begin{displaymath}
z_{M,m,n}^{[k,k',j]}\in H^{3}_{\mot}(Y_G(K_{M,m,n}), \hh^{[\lambda]}(2-j-(t+t'))).
\end{displaymath} 
We now will realize these classes in \'{e}tale cohomology and use the Hochschild--Serre spectral sequence to find elements in Galois cohomology of the representation attached to a weight $(k+2,k'+2)$ Hilbert cuspform. We will then show that they satisfy the Euler system norm relations.

We consider $f$ a cuspidal Hilbert newform of weight $(k+2,k'+2)$ and of level $K_f\subset G(\A_f)$. We assume $k\equiv k'$ mod $2$ and we write $w=k+2+2t=k'+2+2t'$. Denote with $L$ a number field generated by the Hecke eigenvalues $\{\lambda_{\mathfrak{m}}\}_{\mathfrak{m}\subset \Oo_F}$ and fix a prime $p$. We fix an arbitrary place $v$ of $L$ dividing $p$. 

As in the case of classical modular forms, for every finite place $w$ of $L$, one has a Galois representation
\begin{displaymath}
\rho_{f,w}: \gal(\bar{F}/F)\rightarrow \GL_2(L_w)
\end{displaymath}
such that for all but finitely many prime (e.g. for $\mathfrak{p}\nmid \N\text{Nm}_{L/\Q}(w)$, if $f$ is of level $\N$), the representation $\rho_{f,w}$ is unramified at $\mathfrak{p}$ and Trace$(\frob_{\mathfrak{p}})=\lambda_{\mathfrak{p}}$. The more detailed characterization of this representation is a result due to Blasius, Rogawski and Taylor. One can then consider the \textit{classical $L$-function attached to $f$}. In \cite{modular}, it is proved that all elliptic curves over real quadratic fields are modular and hence the $L$-function of such an elliptic curve is equal to the $L$-function of $f$ for the Hilbert eigenform $f$ of weight $(2,2)$ corresponding to it.
 
We consider the $L_{v}$-linear representations of $G_{\Q}$ attached to $f$.   
The classes constructed above are not related to this Galois representation, but rather to a $G_{\Q}$-representation obtained from $\rho_{f,v}$.
One defines the \textit{Asai Galois representation attached to $f$} using the tensor induction by
\begin{displaymath}
\rhoas_{f,v} := \left(\bigotimes -\ind\right)_F^{\Q}(\rho_{f,v}) \otimes L_{v}(t+t'): \gal(\bar{\Q}/\Q) \rightarrow \GL_4(L_{v}).
\end{displaymath}
It is called the Asai Galois representation attached to $f$ because it was first considered by Asai in \cite{Asai}. We will denote it with $V_f^{\as}$.

\begin{rmk} 
In the degenerate case $F=\Q\oplus \Q$, we can think of this Asai representation as the representation attached to the Rankin--Selberg convolution of two modular forms. 
\end{rmk}

\begin{defi}\label{defpoli}
For $f$ as above, we define the local Euler factor for $\ell\neq p$ to be
\begin{displaymath}
P_{\ell}^{\as}(f,X):=\det(1-X\frob_{\ell}^{-1}|(V_f^{\as})^{I_{\ell}}),
\end{displaymath}
where $\frob_{\ell}$ is the arithmetic Frobenius at $\ell$ and $I_{\ell}$ is the inertia subgroup at $\ell$. The local Euler factor at $p$ is defined by the same polynomial acting on the Galois representation $\rhoas_{f,w}$ for some auxiliary $w$ such that $p\nmid w$. 

Then the \textit{Asai $L$-function} is defined by 
\begin{displaymath}
\Las(f,s):=\prod_{\ell}P_{\ell}^{\as}(f,\ell^{-s})
\end{displaymath}
\end{defi}
This product converges for Re$(s)>\frac{k+k'}{2}$ and it admits an analytic continuation to the whole complex plane. It also satisfies a functional equation relating the value at $s$ with the value at $k+k'-1-s$. We are interested in writing the Euler factors in terms of the automorphic representation attached to $f$. We will need the following result
\begin{prop}(\cite[Theorem 2]{Asai}, \cite[Proposition 4.3.4]{HMS})
If $f$ is of level $\N$ and $\ell\nmid \nm(\N) p$ then 
\begin{displaymath}
P_{\ell}^{\as}(f,\ell^{t+t'}X)= \begin{cases}
(1-\alpha_1\alpha_2 X)(1-\alpha_1\beta_2 X)(1-\beta_1\alpha_2 X)(1-\beta_1\beta_2 X), &\text{if $\ell=\mathfrak{l}_1\cdot {\mathfrak{l}_2}$ splits in $F$} \\
(1-\alpha X)(1-\beta X)(1-\alpha \beta X^2) &\text{if $\ell$ is inert in $F$}.
\end{cases}
\end{displaymath}
where $\alpha_i,\beta_i$ and $\alpha,\beta$ are the roots of $X^2-a_{\mathfrak{l}_i}(f)X+\ell^{w-1}\varepsilon_{\mathfrak{l}_i}(f)$ and of $X^2-a_{\ell}(f)X+\ell^{2(w-1)}\varepsilon_{\ell}(f)$ respectively, where $w=k+2+2t=k'+2+2t'$. 
\end{prop}
Now, recall that from the action of $G(\A_f)$ on $H^2_{\et}\left((Y_G)_{\bar{\Q}},\hh_{L_v}^{(\lambda)}(t+t')\right)$ we obtain the finite part of the automorphic representation corresponding to $f$. We will denote it with $\Pi_{f}=\otimes'_{\ell}\Pi_{\ell}$, where $\Pi_{\ell}$ is a $G(\Q_{\ell})$-representation and it is spherical for all but finitely many primes $\ell$. We can describe these $\Pi_{\ell}$ and relate the local $L$-factor with the Euler factor at $\ell$ using the previous proposition.  
\begin{cor}\label{corpoli} For $\ell$ as above, let 
\begin{displaymath}
\sigma=\begin{cases}
I_{G(\Q_{\ell})}(\underline{\chi},\underline{\psi})  &\text{if $\ell$ splits} \\
I_{G(\Q_{\ell})}(\chi,\psi) &\text{if $\ell$ is inert},
\end{cases}
\end{displaymath}
where
\begin{displaymath}
\begin{matrix}
\chi_1(\ell)=\alpha_1 \ell^{-1/2}, & \psi_1(\ell)=\beta_1 \ell^{-1/2}, \\ \chi_2(\ell)=\alpha_2 \ell^{-1/2}, & \psi_2(\ell)=\beta_2 \ell^{-1/2}.
\end{matrix}
\ \ \ \ \text{and} \ \ \ \ \chi(\ell)=\alpha \ell^{-1},  \psi(\ell)=\beta \ell^{-1}.
\end{displaymath}
We then find that $\Pi_{\ell}\simeq \sigma$ and 
\begin{displaymath}
P_{\ell}^{\as}(f,\ell^{-1-s+t+t'})=\begin{cases}
L(\sigma, s) \\
L(\as(\sigma),s)
\end{cases}
\end{displaymath}
\begin{proof}
This follows from the above Proposition and by applying Theorem \ref{sphericalprinc} and Remark \ref{sphericalprincrmk}. First we deal with the split prime case. At a place $\ell$ as above the spherical representation is determined by the values $\chi_i(\ell),\psi_i(\ell)$ being the roots of $X^2-\ell^{-1/2}\lambda_i X +\mu_i$, where $\lambda_i,\mu_i$ are the eigenvalues of $T(\mathfrak{l}_i),R(\mathfrak{l}_i)$. Since $f$ is a newform we have $\lambda_i=a_{\mathfrak{l}_i}(f)$ and  $\mu_i=\ell^{w-2} \varepsilon_{\mathfrak{\ell}_i}(f)$. Hence we need to solve 
\begin{displaymath}
\begin{matrix}
a_{\mathfrak{l}_1}(f)=\alpha_1+\beta_1=\ell^{1/2}(\chi_1(\ell)+\psi_1(\ell)), & \ell^{w-2} \varepsilon_{\mathfrak{l}_1}(f)=\ell^{-1}\alpha_1\beta_1=\chi_1(\ell)\psi_1(\ell);\\
a_{\mathfrak{l}_2}(f)=\alpha_2+\beta_2=\ell^{1/2}(\chi_2(\ell)+\psi_2(\ell)), & \ell^{w-2} \varepsilon_{\mathfrak{l}_2}(f)=\ell^{-1}\alpha_2\beta_2=\chi_2(\ell)\psi_2(\ell).
\end{matrix}
\end{displaymath}
From where we find the claimed values of $\chi_i(\ell),\psi_i(\ell)$.

For the inert prime case we proceed similarly, finding $\chi(\ell),\psi(\ell)$ to be roots of $X^2-(\ell^2)^{-1/2}\lambda X +\mu$, where $\lambda,\mu$ are the eigenvalues of $T(\ell),R(\ell)$. Now $\lambda=a_{\ell}(f)$ and  $\mu=\ell^{2(w-2)} \varepsilon_{\ell}(f)$. Hence from
\begin{displaymath}
a_{\ell}(f)=\alpha+\beta=\ell(\chi(\ell)+\psi(\ell)), \ \  \ell^{2(w-2)} \varepsilon_{\ell}(f)=\ell^{-2}\alpha\beta=\chi(\ell)\psi(\ell)
\end{displaymath}
we find the claimed values of $\chi(\ell),\psi(\ell)$.
\end{proof}
\end{cor}

We use the characterisation of the local components of $\Pi=\Pi_f$ obtained in the previous corollary to prove that if the Hilbert modular form is not a base change lift of a modular form of $\GL_2/\Q$, then a certain $\Hom$-space is zero. We denote with $\omega_{\Pi}$ the Hecke character of $F$ given by the central character of $\Pi$ and we let $\chi_{\Pi_{\ell}}$ be the character of $\Q_{\ell}^{\times}$ given by the restriction of $\Pi_{\ell}$ to the center of $H(\Q_{\ell})$.

\begin{prop}\label{propbasechange}
Let $\tau$ be the representation of $H(\A_f)$ given by $\gamma(\det)$, where $\gamma$ is a character of the id\`{e}les of $\Q$ such that $\gamma_{\ell}^2$ is equal to $\chi_{\Pi_{\ell}}$ for every $\ell$. If $\Pi$ is not a twist of a base change lift of a cuspidal representation of $H(\A_f)$, then
\begin{displaymath}
\Hom_{H(\A_f)}(\Pi,\tau)=0.
\end{displaymath}
\begin{proof}
We will assume for simplicity that $\gamma$ is trivial. If $\Hom_{H(\A_f)}(\Pi,\tau)\neq 0$ then $\Hom_{H(\Q_{\ell})}(\Pi_{\ell},\mathbf{1})\neq 0$ for every $\ell$. In particular for all primes $\ell$ as above which split in $F$, we have 
\begin{displaymath}
\Hom_{H(\Q_{\ell})}(I(\chi_1,\psi_1)\otimes I(\chi_2,\psi_2),\mathbf{1})=\Hom_{H(\Q_{\ell})}(I(\chi_1,\psi_1),I(\chi_2^{-1},\psi_2^{-1}))\neq 0. 
\end{displaymath}
Hence $\Pi_{\ell}=\Pi_{\lambda}\otimes \Pi_{\bar{\lambda}}$ is of the form
\begin{displaymath}
I(\chi_1,\psi_1)\otimes I(\chi_1^{-1},\psi_1^{-1}) \ \ \ \text{ or } \ \ \ I(\chi_1,\psi_1)\otimes I(\psi_1^{-1},\chi_1^{-1}).
\end{displaymath}
Hence $\Pi_{\lambda}\simeq \Pi_{\bar{\lambda}}\otimes \chi_1\psi_1$. 
Letting $\sigma$ be the non trivial automorphism of $F/\Q$ and $\sigma(\Pi)_{\lambda}=\Pi_{\sigma(\lambda)}$, we hence found that the representations $\Pi$ and $\sigma(\Pi) \otimes \omega_{\Pi}$ are isomorphic at all but finitely many primes. Moreover $\omega_{\Pi}$ restricted to the id\`{e}les of $\Q$ is trivial.
We can then apply \cite[Theorem 2(a)]{twists}, which implies that $\Pi$ is a twist of a base change lift of a cuspidal representation of $\GL_2/\Q$ and reach the desired contradiction. 
\end{proof}
\end{prop}

We now see that the Asai representation appears in the parabolic \'{e}tale cohomology of $Y_G$.  Write $\lambda=(k,k',t,t')$. 
\begin{defi}
We define $\hh_{L_v}^{[\lambda]}$ to be the \'{e}tale sheaf of $L_v$-vector spaces on $Y_G$, for $U$ sufficiently small, which is the \'{e}tale realisation of the motivic sheaf $\hh_{L}^{[\lambda]}$ of $\S$\ref{sechilbertsheaves}. We denote with $\hh_{L_v}^{(\lambda)}$ its dual.
\end{defi}

For simplicity let $\mathcal{L}:=\hh^{(\lambda)}_{L_v}(t+t')$. We consider parabolic \'{e}tale cohomology: let $Y_G^{BB}$ be the Bailey-Borel compactification of $Y_G$ and write $j:Y_G \to Y_G^{BB}$ for the natural open embedding. Then \textit{parabolic cohomology} is defined by
\begin{displaymath}
H^i_{\et,!}(Y_{G,\bar{\Q}},\mathcal{L})=\lim_{\to K}H^i_{\et}((Y_{G}(K)^{BB})_{\bar{\Q}},j_{!*}\mathcal{L}).
\end{displaymath}

These cohomology groups have both a $G_{\Q}$ and a $G(\A_f)$ action. 

\begin{thm}[\cite{nekovar},\cite{BL}]\label{thmdec} Let $\mathcal{L}$ be as above, with $\lambda=(k,k',t,t')$ where $k+2t=k'+2t'$. There is a $G_{\Q}\times G(\A_f)$-equivariant decomposition
\begin{displaymath}
H^2_{\et,!}(Y_{G,\bar{\Q}},\mathcal{L})=\bigoplus_{\Pi}V_{\Pi}\otimes \Pi^{\vee},
\end{displaymath}
where $\Pi$ runs over the finite part of cuspidal automorphic representations $\Pi\otimes \Pi_{\infty}$ of $G$ where $\Pi_{\infty}$ is a discrete series of weight $(k+2,k'+2)$. We denote with $\Pi^{\vee}$ its dual $G(\A_f)$-representation and $V_{\Pi}$ is the $G_{\Q}$-representation defined by the tensor induction of $\rho_{\Pi}$ twisted by $t+t'$, where $JL(\rho_{\Pi})=\Pi$. In other words, if $\Pi$ is the automorphic representation generated by a Hilbert cuspidal eigenform $f$, $\rho_{\Pi}=\rho_{f,v}$ and $V_{\Pi}=V_f^{\as}$.
\end{thm}

Taking the dual (as $G_{\Q}$-module) of the cohomology group in the theorem, we get a $G_{\Q}\times G(\A_f)$-equivariant decomposition
\begin{displaymath}
H^2_{\et,!}(Y_{G,\bar{\Q}},\hh^{[\lambda]}_{L_v}(2-(t+t')))=\bigoplus_{\Pi}V_{\Pi}^*\otimes \Pi^{\vee}.
\end{displaymath}

Let's now fix an automorphic representation $\Pi$. We have the following
\begin{prop}(see \cite[$\S$ 4.4]{HMS}).\label{localis}
Let $K\subset G(\A_f)$ be a level such that $\Pi^K\neq 0$ and $T$ a set of primes including the ones at which $K$ is ramified. Let $I$ be the maximal ideal of the Hecke algebra away from $T$ given by the kernel of the action on $\Pi^K$. Then the localisation at  $I$ of $H^i_{\et}(Y_G(K)_{\bar{\Q}},\mathcal{L})$ is zero for $i\neq 2$ is 0 and is equal to the localisation of parabolic cohomology for $i=2$.

Moreover such localisation is given by 
\begin{displaymath}
\left(H^2_{\et}(Y_G(K)_{\bar{\Q}},\mathcal{L})\right)_{I}=\left(H^2_{\et,!}(Y_G(K)_{\bar{\Q}},\mathcal{L})\right)_{I}=\bigoplus_{\Pi'\sim \Pi}V_{\Pi'}\otimes ((\Pi')^{\vee}[t+t'])^K,
\end{displaymath}
where the sum runs over automorphic representations as in the theorem such that $\Pi'_v\simeq \Pi_v$ for almost all $v$. In particular the localisation is independent on $T$.
\end{prop}

Now recall that the target of our map $\AF$ is $ H^{3}_{\mot}(Y_G, {\hh}_L^{[\lambda]}(2-j-(t+t')))$. Let $f$ be a fixed Hilbert eigenform of weight $(k+2,k'+2)$ and $\Pi$ the corresponding $G(\A_f)$-representation, so that $V_{\Pi}^*=(V_f^{\as})^*$. In order to find classes in Galois cohomology of $(V_f^{\as})^*$ we will, roughly, use continuous \'{e}tale realisation map and then apply the above proposition together with Hochschild--Serre spectral sequence. We will find a $G(\A_f)$-equivariant map
\begin{displaymath}
pr_{\Pi}: H^3_{\mot}(Y_G,\hh^{[\lambda]}_{L}(2-j-(t+t')))\longrightarrow H^1(\Q,(V_{\Pi})^*(-j))\otimes \Pi^{\vee}.
\end{displaymath}
We work for any $K$ level subgroup of $G(\A_f)$. 
\begin{itemize}
\item We have (see \cite{huber}) a realisation functor for continuous \'{e}tale cohomology (as defined in \cite{jan}) for varieties defined over $\Q$
\begin{displaymath}
r_{\et}: H^3_{\mot}(Y_G(K),\hh^{[\lambda]}_{L}(2-j-(t+t')))\longrightarrow H^3_{\et}(Y_G(K),\hh^{[\lambda]}_{L_v}(2-j-(t+t'))).
\end{displaymath}
\item There is an Hochschild--Serre spectral sequence (see again \cite{jan}) relating continuous \'{e}tale cohomology for varieties over $\Q$ with \'{e}tale cohomology of the base change over $\bar{\Q}$
\begin{displaymath}
E_2^{p,q} = H^p(\Q,H^q_{\et}(Y_G(K),\mathcal{D})) \Rightarrow H^{p+q}_{\et}(Y_G(K)_{\bar{\Q}},\mathcal{D}).
\end{displaymath}
From this, one gets a map from the kernel of  the map $H^i_{\et}(Y_G(K),\mathcal{D}) \to H^i_{\et}(Y_G(K)_{\bar{\Q}}, \mathcal{D})^{G_{\Q}}$ to $H^1\left(\Q, H^{i-1}_{\et}(Y_G(K)_{\bar{\Q}}, \mathcal{D})\right)$. In particular, for $i=3$, since Artin vanishing theorem tells us that $H^i_{\et}(Y_G(K)_{\bar{\Q}},\mathcal{D})=0$ being $i>\dim(Y_G(K))=2$, we obtain a map 
\begin{displaymath}
HS: H^3_{\et}(Y_G(K),\hh^{[\lambda]}_{L_v}(2-j-(t+t'))) \longrightarrow H^1\left(\Q, H^{2}_{\et}(Y_G(K)_{\bar{\Q}}, \hh^{[\lambda]}_{L_v}(2-j-(t+t')))\right).
\end{displaymath}
\item We now localise at the maximal ideal $I$ given by the kernel of the Hecke algebra acting on $\Pi^K$ as in Proposition \ref{localis}. Applying such proposition and projecting to the $\Pi$-isotypic part we find
\begin{displaymath}
(H^3_{\et}(Y_G(K),\hh^{[\lambda]}_{L_v}(2-j-(t+t'))))_I \longrightarrow H^1\left(\Q, (V_f^{\as})^*(-j)\right) \otimes (\Pi^{\vee})^K.
\end{displaymath}
\end{itemize}

Since all these maps are compatible with respect to changing $K$ and since, by Proposition \ref{localis}, the localisation is independent on the choice of the set of primes $T$ (which may vary changing $K$), we can construct a map of $G(\A_f)$-representation.

\begin{defi}
We define $pr_{\Pi}$ to be the $G(\A_f)$-equivariant map
\begin{displaymath}
pr_{\Pi}: H^3_{\mot}(Y_G,\hh^{[\lambda]}_{L}(2-j-(t+t'))) \longrightarrow H^1\left(\Q, (V_f^{\as})^*(-j)\right) \otimes \Pi^{\vee}
\end{displaymath}
obtained by the previous steps and taking the limit with respect to $K$. 
\end{defi}

In order to define classes in Galois cohomology, we need to take a ``projection'' to $H^1\left(\Q, (V_f^{\as})^*(-j)\right)$ from the target of the map in the previous definition. To do that we assume that $\Pi$ is ordinary at $p$. This means that, letting 
\begin{displaymath}
K_0(p)=\{\gamma\in G(\Z_{p}): \gamma \equiv  \SmallMatrix{*&*\\0&*}  \text{ mod } p\},
\end{displaymath}
the eigenvalue of the Hecke operator $U(p)$ acting on $\Pi_p^{K_0(p)}$ (which is equal to the one of $U'(p)$ acting on $(\Pi_p^{\vee})^{K_0(p)}$) is a $p$-adic unit; we will denote it with $\alpha$. We fix a finite set of primes $S$ to be set of primes outside which $\Pi_{\ell}$ is a spherical representation. We now fix the local data as in $\S$\ref{deficlass}. Write 
\begin{displaymath}
K':=K_S \times \prod_{\ell\nmid pS}G(\Z_{\ell})\times K_0(p),
\end{displaymath}
where $K_S$ is chosen so that $\Pi^{K'}\neq 0$. Choose a vector $v_{\alpha}\in\Pi^{K'}$ in the $U(p)=\alpha$ eigenspace.  This gives a homomorphism
\begin{displaymath}
v_{\alpha}: (\Pi^{\vee})^{K'}\longrightarrow L_v.
\end{displaymath}

What we are going to do is to consider the image of the $K_{M,m,n}$-invariant classes defined in $\S$\ref{deficlass}, take the image via the $G(\A_f)$-equivariant map $pr_{\Pi}$ and then apply $v_{\alpha}$. For $W,\phi_{M,m,n},\xi_{M,m,n}$ as in $\S$\ref{deficlass}, we consider $z_{M,m,n}^{[k,k',j]}$ as in Definition \ref{motivicdefi},
\begin{displaymath}
z_{M,m,n}^{[k,k',j]}= \tfrac{1}{\vol(W)}\AF (\phi_{M,m,n}\otimes\xi_{M,m,n} ) \in H^{3}_{\mot}(Y_G,\hh^{[\lambda]}(2-j-(t+t'))).
\end{displaymath}
Since these elements actually lied in the $K_{M,m,n}$-invariant subspace of the motivic cohomology group, when we apply the \'{e}tale regulator and the map obtained via Hochschild--Serre we obtain classes in  
\begin{displaymath}
H^1\left(\Q, H^{2}_{\et}(Y_G(K_{M,m,n})_{\bar{\Q}}, \hh^{[\lambda]}_{L_v}(2-j-(t+t')))\right).
\end{displaymath}

Recall from (\ref{eqDET}) that, restricting to $G^*$, we find 
\begin{displaymath}
Y_{G^*}(K^*_{M,m,n})\simeq Y_{G^*}(K^*_n)\times_{\Q} \mu_{Mp^m}.
\end{displaymath}
We now recall a result that will be useful to use the above isomorphism to land in Galois cohomology over cyclotomic extensions.

\begin{prop}\cite[Corollary 5.8]{nekovar}. Let $U\subset G(\A_f)$ be the stabiliser of $Y_{G^*}$. We have a $G_{\Q}\times G(\A_f)$ isomorphism
\begin{displaymath}
H^i_{\et}(Y_{G,\bar{\Q}},\mathcal{L})\simeq \ind_U^{G(\A_f)}H^i_{\et}(Y_{G^*,\bar{\Q}},\iota^*\mathcal{L}),
\end{displaymath}
where the natural embedding $\iota: Y_{G^*}\hookrightarrow Y_G$ is an open immersion. 
\end{prop}

We also recall that for any variety $X$ over $\Q$ we naturally have the following isomorphism of $G_{\Q}$-modules
\begin{displaymath}
H^i_{\et}((X\times_{\Q} \mu_N)_{\bar{\Q}}, \mathcal{L})\simeq \ind_{G_{\Q(\mu_N)}}^{G_{\Q}} H^i_{\et}(X_{\bar{\Q}},\mathcal{L}).
\end{displaymath}
Moreover, by Shapiro's lemma we have
\begin{displaymath}
H^1({\Q},  \ind_{G_{\Q(\mu_N)}}^{G_{\Q}} V)=H^1(\Q(\mu_N),V).
\end{displaymath}
Applying the above proposition and these isomorphisms for $N=Mp^m$ and for the $G_{\Q}$-module $H^{2}_{\et}(Y_{G^*}(K_{M,m,n}^*)_{\bar{\Q}}, \tsym^{[k,k']}\hh_{L_v}(\mathcal{A})(2-j))$, we can give the following
\begin{defi}\label{galoisclasses}
For $m\geq 0$ we define a class
\begin{displaymath}
z^{\Pi,j}_{Mp^m,\alpha}\in H^{1}(\Q(\mu_{Mp^m}),(V_f^{\as})^*(-j))
\end{displaymath}
by letting
\begin{displaymath}
\tfrac{1}{M}\cdot \begin{cases}
\left( \tfrac{p^{j}\sigma_p}{\alpha}\right)^m  \\
\left( 1-\tfrac{p^{j}\sigma_p}{\alpha}\right)  
\end{cases} \cdot (v_{\alpha}\circ \text{pr}_{\Pi^{\vee}})\left(z_{M,m,n}^{[k,k',j]}\right) \ \  \begin{matrix}
\text{if } m\geq 1\\
\text{if } m=0,
\end{matrix}
\end{displaymath}
where $\sigma_p$ is the arithmetic Frobenius at $p$ in $\gal(\Q(\mu_M)/\Q)$. 
\end{defi}

\subsection{Norm relations in Galois cohomology}
These are the classes that, as we are going to show, form an Euler system for $V_f^{\as}(j+1)$. 
\begin{thm}[Vertical norm relations]\label{verticalthm}
Let $j\leq \min(k,k')$ and $k,k'\geq 0$. We have
\begin{displaymath}
\cores_{\Q(\mu_{Mp^m})}^{\Q(\mu_{Mp^{m+1}})} \left( z^{\Pi,j}_{Mp^{m+1},\alpha}\right)=z^{\Pi,j}_{Mp^m,\alpha}.
\end{displaymath}
\begin{proof}
Since the pushforward by
\begin{displaymath}
\pi_m:H^2_{\et}(Y_{G^*}(K_n^*)\times\mu_{Mp^{m+1}}, \tsym^{[k,k']}\hh_{L_v}(\mathcal{A}))\to H^2_{\et}(Y_{G^*}(K_n^*)\times\mu_{Mp^m}, \tsym^{[k,k']}\hh_{L_v}(\mathcal{A}))
\end{displaymath}
induces corestriction in Galois cohomology, the result is an immediate corollary of Theorem \ref{quasivert}. This indeed can be rewritten as
\begin{displaymath}
(\pi_m)_*(z_{M,m+1,n}^{[k,k',j]})=\begin{cases}
\tfrac{U'(p)}{p^j }  \\
\left( \tfrac{U'(p)}{p^j}-1\right)  
\end{cases} \cdot z_{M,m,n}^{[k,k',j]} \ \  \begin{matrix}
\text{if } m\geq 1\\
\text{if } m=0,
\end{matrix}
\end{displaymath}
seen as elements in $H^2_{\et}(Y_{G^*}(K_{M,m,n}^*), \tsym^{[k,k']}\hh_{L_v}(\mathcal{A}))$. Here $U'(p)$ is the Hecke operator given by the double coset of $\SmallMatrix{p^{-1}&0\\0&1}$ in $\hh(K_{M,m,n}^*\backslash G^*(\A_f)/K_{M,m,n}^*)$, and we used the fact that, 
as explained in the proof \cite[Proposition 4.3.4]{HMS}\footnote{The pullback of the projection from $Y_{G^*}$ to $Y_G$ intertwines $U'(p)$ on the cohomology of $Y_{G^*}$ with $p^{-(t+t')}\mathcal{U}'(p)$ on the cohomology of $Y_G$, where $\mathcal{U}'(p)$ is the normalised Hecke operator given by $p^{t+t'}U'(p)$.}, the Hecke operator of Theorem \ref{quasivert}, acts on $H^2_{\et}(Y_{G^*}(K_{M,m,n}^*)_{\bar{\Q}},\tsym^{[k,k']}(\hh_{L_v}(\mathcal{A})))$ as $U'(p)$. The isomorphism (\ref{eqDET}) intertwines $U'(p)$ with  $U'(p)\times \sigma_p^{-1}$, where $U'(p)\in \hh(K_{n}^*\backslash G^*(\A_f)/K_{n}^*)$ and $\sigma^{-1}_{p}$ is the arithmetic Frobenius at $p$ in $\gal(\Q(\mu_M)/\Q)$. Since $v_{\alpha}$ projects to the $U'(p)= \alpha$ eigenspace, the theorem follows. 
\end{proof}
\end{thm}
 
\begin{thm}[Tame norm relations]\label{tame}
Let $j\leq \min(k,k')$ and $k,k'\geq 0$. We assume that $\Pi$ is not a twist of a base change lift of a cuspidal representation of $\GL_2/\Q$. For any $\ell\nmid Mp$, $\ell\not\in S$, we have
\begin{displaymath}
\cores_{\Q(\mu_{Mp^m})}^{\Q(\mu_{\ell Mp^{m}})} \left( z^{\Pi,j}_{\ell Mp^{m},\alpha}\right)=Q(\sigma_{\ell}^{-1})z^{\Pi,j}_{Mp^m,\alpha}, 
\end{displaymath}
where $Q(X)=\det(1-X\frob_{\ell}^{-1}|V_f^{\as}(1+j))$, i.e. $Q(\sigma_{\ell}^{-1})=P_{\ell}(\ell^{-1-j}\sigma_{\ell}^{-1})$ for $P_{\ell}(X)=\det(1-X\frob_{\ell}^{-1}|V_f^{\as})$ as in Definition \ref{defpoli}.
\begin{proof}
First of all we notice as above that the corestriction map is induced by pushforward under the projection $\pi:Y_G(K_{\ell M,m,n})\rightarrow Y_G(K_{M,m,n})$. The class on the left is then obtained in motivic cohomology by applying 
\begin{displaymath}
\Ss(\A_f^2,\Z)[j+t+t']\otimes \hh(G(\A_f)) \xrightarrow{\AF} H^{3}_{\mot}(Y_G(K_{\ell M,m,n}), \mathcal{D}(2-j))\xrightarrow{\pi_{*}}H^{3}_{\mot}(Y_G(K_{M,m,n}), \mathcal{D}(2-j))
\end{displaymath}
We have $\pi_*\circ\AF(\phi\otimes\xi)=\sum_{k}k\cdot \AF(\phi\otimes\xi)=\sum_{k}\AF(\phi\otimes (k\cdot \xi))$, where $k$ runs over coset representatives of $K_{\ell M,m,n}/ K_{M,m,n}$. In particular we find that
\begin{displaymath}
\sum_{k}k\cdot \xi_{\ell M,m,n}=\xi'_{\ell M,m,n}, 
\end{displaymath}
where $\xi'_{\ell M,m,n}$ is equal to $\xi_{\ell M,m,n}$ at every component but at $\ell$ where we find 
\begin{displaymath}
\sum_{k}k_{\ell}\cdot\left( \ch(K_{\ell,1})-\ch(\eta_{\ell,1}K_{\ell,1}) \right) =\ch(G(\Z_{\ell}))-\ch(\eta_{\ell,1}G(\Z_{\ell})).
\end{displaymath}
Hence both the left hand side and the right hand side of the claimed equality are obtained as image of the same map $v_{\alpha}\circ \text{pr}_{\Pi^{\vee}}\circ \AF$
\begin{displaymath}
\Ss(\A_f^2,\Z)[j+t+t']\otimes \hh(G(A_f),\Z) \longrightarrow H^{1}({\Q(\mu_{Mp^m})},(V_f^{\as})^*(-j)).
\end{displaymath}
They are obtained as image of elements that are the same at every component different from $\ell$, where the right hand side is the image of $\tfrac{\ell -1}{\ell}\left( (\ch(\ell^2\Z_{\ell} \times (1+\ell^2\Z_{\ell})) \otimes (\ch(G(\Z_{\ell}))-\ch(\eta_{\ell,1}G(\Z_{\ell})))\right)$ and the left hand side the image of $\ch(\Z_{\ell}^2)\otimes \ch(G(\Z_{\ell}))$. The factor $\ell -1$ appears comparing $\vol(W)$ for the two different motivic classes, while $\tfrac{1}{\ell}$ comes from the $\tfrac{1}{M}$-factor in the definition of the Galois cohomology classes. So it is enough to compare the image of these two elements via the component at $\ell$ of the above map. We will first compare the images through each of the the maps
\begin{equation}\label{eqmapproof}
\Zita: \Ss(\Q_{\ell}^2,\Z)[j+t+t']\otimes \hh(G(\Q_{\ell}),\Z) \longrightarrow H^{1}({\Q(\mu_{Mp^m})},(V_f^{\as})^*(-j)) \otimes \Pi_{\ell}^{\vee}\longrightarrow \Pi_{\ell}^{\vee},
\end{equation}
where the last map is is a $G(\A_f)$-equivariant projection to $\Pi_{\ell}^{\vee}$, obtained by choosing a basis element of $H^{1}({\Q(\mu_{Mp^m})},(V_f^{\as})^*(-j))$. Note that this Galois cohomology group is a priori infinite dimensional, but since it is actually equal to the Galois cohomology of some maximal unramified (outside a finite set of places) extension, we are reduced to take this projection map for a finite number of basis elements. 
First we assume that $k+k'-2j\neq 0$. By definition and by Theorem \ref{thmeis}, Proposition \ref{propeis} and Theorem \ref{thmdec}, we find that $\Zita$ satisfies the condition  of Corollary \ref{THEcor}, for integers $k$ and $h$ equal to $k+k'-2j$ and $j+t+t'$ respectively. Condition (\ref{inertmult1}) follows from purity. 

If $M=1$ we can apply then Corollary \ref{THEcor} with $\sigma=\Pi_{\ell}$. The factor of discrepancy is then $L(\sigma,h)^{-1}$. We then apply Corollary \ref{corpoli} to get
\begin{displaymath}
L(\sigma,h)^{-1}=L(\Pi_{\ell}, j+t+t')^{-1}=P_{\ell}(\ell^{-1-j}).
\end{displaymath}
The multiplication by such scalar is carried when we take the projection via $v_{\alpha}$ into Galois cohomology and this is precisely what we were looking for (since $\sigma_{\ell}$ is trivial in this situation).

If $M>1$, we apply this to every twist by Dirichlet characters modulo $M$ and apply Shapiro's lemma. First of all we suppose $m=1$. This is no problem since Theorem \ref{verticalthm} allows us to reduce to this case. We are now comparing classes in $H^1({\Q(\mu_M)}, (V_f^{\as})^*(-j))$. Since $\rho:=(V_f^{\as})^*(-j)$ is a $G_{\Q}$-module, we have $\ind^{G_{\Q(\mu_M)}}_{G_{\Q}}(\rho)=\bigoplus_{\eta}\rho \otimes \eta$, where $\eta$ varies over all characters of the quotient $G_{\Q}/G_{\Q(\mu_M)}=\gal(\Q(\mu_M)/\Q)\simeq (\Z/M\Z)^{\times}$. We hence find
\begin{displaymath}
H^1({\Q(\mu_M)}, (V_f^{\as})^*) = \bigoplus_{\eta}H^1({\Q},(V_f^{\as})^*\otimes \eta).
\end{displaymath}
Since $\sigma_{\ell}$ is the image of $\ell^{-1}$ in $\gal(\Q(\mu_M)/\Q)$, if we write $z\in H^1({\Q(\mu_M)}, (V_f^{\as})^*(-j))$ as $(z_{\eta})_{\eta}$, we have that $\sigma_{\ell}^{-1}\cdot z = (\eta(\ell)\cdot z_{\eta}){_\eta}$. Hence we have reduced to prove that the $\eta$-components of the classes we are considering differ by the factor $P_{\ell}(\ell^{-1-j}\eta(\ell))$. We are then again in the case $M=1$. The character $\eta$ can be seen as an unramified character of $\Q_{\ell}^{\times}$ for $\ell\nmid M$ via class field theory and it then defines a one dimensional representation of $G^*(\Q_{\ell})$ and of $H(\Q_{\ell})$ via the determinant map. Hence the classes $z_{\eta}$ we are considering are locally at $\ell$ images of the map (\ref{eqmapproof}) with the action of $H(\Q_{\ell})$ twisted by $\eta$. The space of such maps factoring through the Siegel section will now be isomorphic, via the bijection of Proposition \ref{bij}, to a space of the form $\Hom_H(I_H(\chi\eta,\psi\eta)\otimes \Pi_{\ell},\eta)$, where $\chi\psi\cdot\chi_{\Pi_{\ell}}=1$. Theorem \ref{thmmult1} implies that this space is again one dimensional, and the construction of a basis carries through as in Section $\S$ \ref{towards}, where in the choice of the auxiliary character in Definition \ref{defizetino} $\psi$ is replaced by $\psi\eta$. We obtain the same results, but with $L(\sigma\otimes \eta,h)$ in place of $L(\sigma,h)$. We then find, as we wanted,
\begin{displaymath}
L(\Pi_{\ell}\otimes \eta, j+t+t')^{-1}=P_{\ell}(\ell^{-1-j}\eta(\ell)).
\end{displaymath}

We are left with the case $k+k'-2j = 0$. The issue here is that the divisor map from $\Oo^{\times}(Y)\otimes \C$ in (1) of Theorem \ref{thmeis} has a kernel. It consists of non-generic representations of $H(\A_f)$. For any such representation $\tau$ we have that $\Hom_{H(\A_f)}(\tau\otimes \Pi,\C)= 0$ thanks to the assumption that $\Pi$ is not a base change lift from $\GL_2/\Q$ and Proposition \ref{propbasechange}. Hence the local map factors through the Siegel section also in this case and the proof follows as above.
\end{proof}
\end{thm}

\begin{rmk}
These classes hence satisfy the Euler system norm relations (\ref{eqNR}) as stated in the Introduction. In particular we proved the tame norm relations for all primes $\ell\not\in S$. In \cite{HMS} these were proved only for $\ell$ inert in $F$ or $\ell$ split with the condition of the two primes ideal in $F$ above it being narrowly principal. 
\end{rmk}

\begin{rmk}[Integral classes] In fact, one is interested in ``integral classes'': fixing a $G_{\Q}$-stable lattice $T\subset (V_f^{\as})^*(-j)$, we would like to have classes in $H^1(\Q({\mu_m}),T)$ satisfying the same norm relations. To do that one works with integral Eisenstein classes, applies the map $\AF$ and slightly modifies the projection map $v_{\alpha}\circ pr_{\Pi}$ by choosing an appropriate Hecke operator that will define a lattice as above. This is explained in details for the case $G=$GSp$_4$ in \cite{GSP4} in the discussion following Proposition 10.5.2.
\end{rmk}

\subsection{A remark on Beilison--Flach Euler system}
It should now be clear to the reader that, proceeding in a completely analogous way, one can reprove Euler system norm relations for Beilinson--Flach classes. These elements were constructed in \cite{RSCMF} and \cite{kings} and lay in Galois cohomology of the representation attached to the Rankin--Selberg convolution of two modular forms $f,g$ of weight $k+2,k'+2$ respectively. This means that in this case one works with $\Pi=\Pi_f\otimes \Pi_g$, where $\Pi_f,\Pi_g$ are automorphic representations of $\GL_2(\A_f)$. Hence we have, at all but finitely many places, a spherical representation $\Pi_{\ell}$ of $G(\Q_{\ell})$ as in Definition \ref{sphgl2gl2}, where now $G=\GL_2\times \GL_2$. Using $\S$\ref{splitzeta}, one can restate all the results of $\S$\ref{towards} for $G$; everything is already there, since we are in the degenerate case where all primes split. One then defines a map
\begin{displaymath}
\BF: \Ss(\A_f^2,\Q)[j]\otimes \hh(G(\A_f),\Z)\longrightarrow H^{3}_{\mot}(Y_{G}, \tsym^{[k,k']}\hh(\mathcal{E})(2-j)),
\end{displaymath}
similarly as in $\S$\ref{defimap}, where now $\tsym^{[k,k']}\hh_{L}(\mathcal{E})$ is a motivic sheaf over the $\GL_2\times \GL_2$ Shimura variety and the considered embedding $\iota$  at the level of algebraic groups is the diagonal embedding $\GL_2\hookrightarrow \GL_2\times \GL_2$. In this case one uses
\begin{displaymath}
CG_{\mot}^{[k,k',j]}: \tsym^{k+k'-2j}\hh(\mathcal{E}) \to \iota^*( \tsym^{[k,k']}\hh(\mathcal{E}))(-j), 
\end{displaymath}
as in \cite[Corollary 5.2.2]{kings}. The local input is then the same as in $\S$\ref{deficlass}, again in the ``all split primes'' case. The proofs of all results in $\S$\ref{ES} carry over, where in this setting the Galois representation $(V_f\otimes V_g)^*$ appears in $H^{2}_{\et}(Y_{G,\bar{\Q}}, \tsym^{[k,k']}\hh(\mathcal{E})(2))$. 

Both in this case and in the Asai--Flach one, the obtained classes are not exactly the ones obtained pushing forward Eisenstein classes via ``perturbed embeddings''. The classes explicitly defined this way in \cite{RSCMF} and \cite{HMS} satisfy the expected tame norm relations at $\ell$ only modulo $(\ell -1)$; one obtains an Euler system thanks to a result by Rubin stating that these relations are enough to ``lift'' such classes to an Euler system. This error term does not appear in this setting because at primes $\ell\mid M$ we already add a correction term in the definition of the local Hecke algebra element $\xi_{M,m,n}$ (see Remark \ref{rmkHMS}). This can be seen to be the right choice from the local computation of Corollary \ref{THEcor}.

\section*{Funding}
This work was supported by the Engineering and Physical Sciences Research Council [EP/L015234/1], the EPSRC Centre for Doctoral Training in Geometry and Number Theory (The London School of Geometry and Number Theory), University College London.

\bibliographystyle{alpha}
\bibliography{biblio}

\Addresses

\end{document}